\input amstex
\documentstyle{amsppt}

\magnification=1200 
 
\TagsOnRight 
\NoRunningHeads
\NoBlackBoxes
 
\def\C{{\Bbb C}} 
\def\R{{\Bbb R}} 
\def\Z{{\Bbb Z}} 
\def\gl{\frak {gl}}

\def\p{\partial} 

\def\det{\operatorname{det}}
\def\diag{\operatorname{diag}}
\def\dim{\operatorname{dim}}
\def\per{\operatorname{per}} 

\def\deg{\operatorname{deg}} 
\def\gr{\operatorname{gr}} 
\def\Hom{\operatorname{Hom}}
\def\M{\operatorname{M}}
\def\Mnm{\operatorname{M}(n,m)} 
 
\def\RTab{\operatorname{RTab}}
\def\tr{\operatorname{tr}}
\def\const{\operatorname{const}} 
\def\sgn{\operatorname{sgn}}

\def\End{\operatorname{End}}
\def\per{\operatorname{per}}

\def\mi{\kern0.8pt|\kern0.8pt}

\def\L{\Lambda}
\def\Ls{\Lambda^*}
\def\gln{\frak{gl}(n)}
\def\glm{\frak{gl}(m)}
\def\Ugln{\Cal U(\gln)}
\def\Uglm{\Cal U(\glm)}

\def\U{\Cal U}

\def\pc{\psi_{C}}

\def\l{\lambda}

\def\si{\sigma}
\def\chm{\chi^\mu}

\def\tht{\thetag}
\def\dt{\delta}
\def\f{\downharpoonright}
\def\u{\upharpoonright}

\def\a{\alpha}
\def\dimn{\dim_{GL(n)}}
\def\ov{\,/\,}
\def\S{\Bbb S}
\def\la{\langle}
\def\ra{\rangle}
\def\gl{\frak{gl}}
\def\glN{\frak{gl}(N)}
\def\Zgln{\frak Z(\gln)}
\def\ZglN{\frak Z(\glN)}
\def\UglN{\Cal U(\glN)}
\def\Zglm{\frak Z(\glm)}
\def\Dmnm{\Delta_\mu^{(n,m)}}
\def\DmNM{\Delta_\mu^{(N,M)}}
\def\Ind{\operatorname{Ind}}
\def\Res{\operatorname{Res}}

\def\ps{p^*}
\def\Ps{P^*}
\def\pc{\check p}
\def\ph{\hat p}
\def\o{\omega}

\def\s{s^*}

\def\tht{\thetag}
\def\dt{\delta}
\def\f{\downharpoonright}
\def\u{\upharpoonright}

\def\a{\alpha}
\def\dimn{\dim_{GL(n)}}
\def\ov{\,/\,}
\def\S{\Bbb S}
\def\Sm{\S_\mu}
\def\Smn{\S_{\mu|n}}

\def\Smnt{\widetilde{\S}_{\mu|n}}
\def\SmNt{\widetilde{\S}_{\mu|N}}

\def\la{\langle}
\def\ra{\rangle}
\def\gl{\frak{gl}}
\def\glN{\frak{gl}(N)}
\def\Zgln{\frak Z(\gln)}
\def\ZglN{\frak Z(\glN)}
\def\UglN{\Cal U(\glN)}
\def\Zglm{\frak Z(\glm)}
\def\Dm{\Delta_\mu}
\def\Ind{\operatorname{Ind}}
\def\Res{\operatorname{Res}}
\def\Av{\operatorname{Avr}_{nN}}

\def\ps{p^*}
\def\Ps{P^*}
\def\pc{\check p}
\def\ph{\hat p}
\def\o{\omega}

\def\1{\widehat{\,1\,}}
\def\tM{\widetilde{M}}

\def\li #1 #2 #3 {\vrule height #1  width #2  depth #3 }
% \li draws a black box

\def\hl #1 #2 #3 {\rlap{ \kern #1 pc \raise -#2 pc 
\hbox{ \li {0.1pt} {#3pc} 0pt }}}
% \hl draws a horizontal line starting in the point (#1,#2)
% with length #3 and width #4

\def\vl #1 #2 #3 {\rlap{ \kern #1 pc \raise -#2 pc
\hbox{ \li {#3pc} {0.1pt} 0pt }}}
% \vl draws a vertical line starting in the point (#1,#2)
% with length #3 and width #4

\def\bhl #1 #2 #3 {\rlap{ \kern #1 pc \raise -#2 pc
\hbox{ \li {1pt} {#3pc} 0pt }}}
% \hl draws a horizontal line starting in the point (#1,#2)
% with length #3 and width #4

\def\bvl #1 #2 #3 {\rlap{ \kern #1 pc \raise -#2 pc
\hbox{ \li {#3pc} {1pt} 0pt }}}
% \vl draws a vertical line starting in the point (#1,#2)
% with length #3 and width #4

\def\wr #1 #2 #3 {\rlap{ \kern #1 pc \raise -#2 pc \hbox{\tenrm #3 }}}
% \wr writes #3 starting from the point (#1,#2)

\def\boxit #1 #2 #3 {\rlap{ \kern #1 pc \raise -#2 pc
\vbox to 2 pc{
\vfill \hbox to 2pc{\hfill {\tenrm #3} \hfill } \vfill }}}

\def\botit #1 #2 #3 {\rlap{ \kern #1 pc \raise -#2 pc 
\hbox{\raise 4 pt
\hbox to 2pc{\hfill {\tenrm #3} \hfill }}}}

\def\picture #1 #2 #3 \endpicture
   {\smallskip \centerline{\hbox to #1 pc 
   {\nullfont #3 \hss}} \nobreak \smallskip \centerline{#2} \smallskip}
% creates a picture with width #1 and name #2 

\topmatter
\title Shifted Schur functions\endtitle
\thanks
The authors were supported by the International Science
Foundation (under grants MBI300 and MQV000, respectively) and by Russian Foundation for Basic Research
(grant 95-01-00814)
\endthanks
\endtopmatter
\bigskip

\centerline{\smc Andrei Okounkov
\footnote{current address: School of Mathematics,
Institute for Advanced Study, Princeton, NJ 08540,
e-mail: okounkov\@math.ias.edu} and Grigori Olshanski}
\bigskip
\centerline{Institute for Problems of Information Transmission}
\centerline{Bolshoy Karetny 19, 101447 Moscow GSP--4, Russia}
\medskip
\centerline{E-mail: okounkov$\@$ippi.ac.msk.su,
olsh$\@$ippi.ac.msk.su} 
\bigskip
\bigskip
\bigskip
\bigskip
\bigskip

\head Abstract \endhead
The classical algebra $\Lambda$ of symmetric functions has a remarkable
deformation $\Lambda^*$, which we call the algebra of shifted symmetric
functions. In the latter algebra, there is a distinguished basis formed by
shifted Schur functions $s^*_\mu$, where $\mu$ ranges over the set of
all partitions. The main significance of the shifted Schur functions is
that they determine a natural basis in $Z(\frak{gl}(n))$, the center of the
universal enveloping algebra $U(\frak{gl}(n))$, $n=1,2,\ldots$. 

The functions $s^*_\mu$ are closely related to the factorial Schur
functions introduced by Biedenharn and Louck and further studied by
Macdonald and other authors.

A part of our results about the functions $s^*_\mu$ has natural classical
analogues (combinatorial presentation, generating series, Jacobi--Trudi
identity, Pieri formula). Other results are of different nature
(connection with the binomial formula for characters of $GL(n)$, an
explicit expression for the dimension of skew shapes $\lambda/\mu$,
Capelli--type identities, a characterization of the functions $s^*_\mu$ by
their vanishing properties, `coherence property', special
symmetrization map $S(\frak{gl}(n))\to U(\frak{gl}(n))$.

The main application that we have in mind is the asymptotic
character theory for the unitary groups $U(n)$ and 
symmetric groups $S(n)$ as $n\to\infty$.

\newpage
\head Contents\endhead

Introduction\par
\item{1.} The algebra of shifted symmetric functions
\item{2.} Quantum immanants
\item{3.} Characterization of $s^*$-functions
\item{4.} Duality
\item{5.} Binomial Theorem
\item{6.} Eigenvalues of higher Capelli operators
\item{7.} Capelli--type identity for Schur--Weyl duality
\item{8.} Dimension of skew Young diagrams
\item{9.} Pieri--type formula for $s^*$-functions
\item{10.} Coherence property of quantum immanants and shifted Schur
polynomials 
\item{11.} Combinatorial formula for $s^*$-functions
\item{12.} Generating series for $h^*$- and $e^*$-functions
\item{13.} Jacobi--Trudi formula for $s^*$-functions
\item{14.} Special symmetrization
\item{15.} Concluding remarks

References

\newpage

\head Introduction \endhead

\subhead 1. Shifted Schur polynomials and factorial Schur polynomials
\endsubhead 
Recall that the {\it Schur function\/} (or {\it Schur polynomial\/}) in
$n$ variables can be defined as ratio of two $n\times n$ determinants
$$
s_\mu(x_1,\ldots,x_n)=\dfrac
{\det[x_i^{\mu_j+n-j}]}{\det[x_i^{n-j}]},\tag0.1
$$
where $\mu$, the parameter of the polynomial, is an arbitrary partition
$\mu_1\ge\mu_2\ge\ldots\ge\mu_n\ge0$ of length $\le n$.

Denote by the symbol $(x\f k)$ the $k$-th {\it falling
factorial power\/} of a variable $x$,
$$
(x\f k)=\cases
x(x-1)\ldots(x-k+1),&\text{if $k=1,2,\ldots$,}\\
1,&\text{if $k=0$.}
\endcases\tag0.2
$$

In the present paper we study the following Schur--type polynomials:
$$
s^*_\mu(x_1,\ldots,x_n)=\dfrac
{\det[(x_i+n-i\f \mu_j+n-j)]}{\det[(x_i+n-i\f n-j)]}.\tag0.3
$$
We call them the {\it shifted Schur polynomials\/}.

These new polynomials differ by a shift of arguments only from the {\it
factorial Schur polynomials\/}
$$
t_\mu(x_1,\ldots,x_n)=\dfrac
{\det[(x_i\f \mu_j+n-j)]}{\det[(x_i\f n-j)]}.\tag0.4
$$
The polynomials $t_\mu$ were introduced by Biedenharn and Louck
\cite{BL1}, \cite{BL2}. \footnote{ Note that the original definition of
these authors was different from \tht{0.4}. The fact that their definition
can be written in the form \tht{0.4} was established later by Macdonald
\cite{M2}.} Then they were studied in Chen--Louck \cite{CL}, Goulden--Greene
\cite{GG}, Goulden--Hamel \cite{GH}, and Macdonald \cite{M2} (see also
Macdonald \cite{M1}, 2nd edition, Ch. I, section 3, Examples 20--21). In
particular, Macdonald developed a theory of more general `factorial'
polynomials including as special cases the polynomials $s_\mu$ and
$t_\mu$. 

In these works it was shown that several important facts about the ordinary
Schur polynomials (e.g., the Jacobi--Trudy identity and the combinatorial
presentation) can be transferred to the factorial polynomials. 

Our results about the shifted Schur
polynomials $s^*_\mu$ can be restated as certain new results about the
factorial polynomials $t_\mu$. However, as it will be shown, the use of the
shifted polynomials has many advantages and provides a new insight.

\subhead 2. Stability and the algebra $\Ls$ of shifted Schur functions
\endsubhead 
Recall that the ordinary Schur polynomials are {\it stable\/} in the
following sense:
$$
s_\mu(x_1,\ldots,x_n,0)=s_\mu(x_1,\ldots,x_n).\tag0.5
$$
This stability property holds also for the shifted polynomials,
$$
s^*_\mu(x_1,\ldots,x_n,0)=s^*_\mu(x_1,\ldots,x_n),\tag0.6
$$
but {\it fails\/} for the factorial polynomials $t_\mu$.

The stability property \tht{0.6} allows us to introduce the {\it
functions\/} $s^*_\mu$ in infinitely many variables --- just in the same
way as for the classical Schur functions. These functions $s^*_\mu$ form a
distinguished basis in a certain new algebra which we denote by $\Ls$ and
call the {\it algebra of shifted symmetric functions\/}. 

As in the classical context of symmetric functions, elements of the
algebra $\Ls$ may be viewed as functions $f(x_1,x_2,\ldots)$ on
infinite sequences of arguments such that $x_i=0$ for $i$ large enough. But
the ordinary symmetry is replaced by the `shifted symmetry':
$$
f(x_1,\ldots,x_i,x_{i+1},\ldots)=
f(x_1,\ldots,x_{i+1}-1,x_i+1,\ldots),\qquad i=1,2,\ldots.\tag0.7
$$

As examples of shifted symmetric functions one can take the {\it complete
shifted functions\/} $h^*_r=s^*_{(r)}$ and the {\it elementary shifted
functions\/} $e^*_r=s^*_{(1^r)}$, $r=1,2,\ldots$:
$$
\gather
h^*_r(x_1,x_2,\ldots)=\sum_{1\le i_1\le\ldots\le i_r<\infty}
(x_{i_1}-r+1)(x_{i_2}-r+2)\ldots x_{i_r},\tag0.8\\
e^*_r(x_1,x_2,\ldots)=\sum_{1\le i_1<\ldots< i_r<\infty}
(x_{i_1}+r-1)(x_{i_2}+r-2)\ldots x_{i_r}.\tag0.9
\endgather
$$

The algebra $\Ls$ (but yet without the functions $s^*_\mu$) appeared in
Olshanski's papers \cite{O1}, \cite{O2} in connection with a construction
of Laplace operators on the infinite--dimensional classical groups. Then
this algebra was studied in the note \cite{KO} by Kerov and Olshanski. 
We will often call the functions $s^*_\mu\in\Ls$ the $s^*$-{\it
functions\/}.  

\subhead 3. The $s^*$-functions and the center of the universal
enveloping algebra $\Ugln$ \endsubhead
Let $\Ls(n)$ denote for the algebra of shifted symmetric polynomials in $n$
variables, and note that the polynomials $s^*_\mu$ form a basis in
$\Ls(n)$. Let $\Zgln$ denote the center of $\Ugln$, the universal
enveloping algebra of $\gln=\gl(n,\C)$.

There exists a canonical algebra
isomorphism 
$$
\Zgln\;\to\;\Ls(n),\qquad A\mapsto f_A,\tag0.10
$$
which plays the key role in our paper. 
This isomorphism is simply the well--known Harish--Chandra isomorphism; it
takes a central element $A\in\Zgln$ to its eigenvalue $f_A(\l_1,\ldots,\l_n)$
in a highest weight module $V_\l$, where $\l$ varies over $\C^n$.

By taking the preimage of the polynomials $s^*_\mu\in\Ls(n)$ under the
isomorphism \tht{0.10} we obtain a certain distinguished basis $\{\Sm\}$ in
the center $\Zgln$. It will be shown that the central elements $\Sm$
possess many remarkable properties.

In Okounkov's paper \cite{Ok1} an explicit formula
expressing the elements $\Sm$ in terms of the standard generators
$E_{ij}\in\gln$ was found; then it was improved in Nazarov \cite {N2} and
in \cite{Ok2}. 
\footnote{ The results of the present paper also provide us with an
expression for $\Sm$ but it is less satisfactory.} In the particular case
of $\mu=(1^k)$, $k=1,\ldots,n$, this formula turns into a classical
formula occurring in the well--known Capelli identity (see 
Howe \cite{H}, Howe--Umeda \cite{HU}). We would like to note that the paper
\cite{HU} was one of the starting points of our work.

Note that a kind of isomorphism \tht{0.10} also exists for the algebra
$\Ls$ although the naive $\infty$-dimensional analogue of the algebras
$\Ugln$, the inductive limit algebra $\varinjlim\Ugln$, has trivial center
(see \cite{O1}, \cite{O2}).

\subhead 4. The $s^*$-functions and Vershik--Kerov's asymptotic theory of
characters \endsubhead
Suppose we have an infinite increasing chain
$$
G(1)\subset G(2)\subset \ldots \tag0.11
$$
of finite or compact groups. In the asymptotic theory of characters (see
Vershik--Kerov \cite{VK1}, \cite{VK2}, \cite{VK3}) a central place is
taken by the problem of the limit behavior of the expression
$$
\chi_\pi(g)=\frac{\tr\pi(g)}{\dim\pi}\;,\tag0.12
$$
where $g$ is an arbitrary but fixed element of a group $G(k)$ while
$\pi=\pi_n$ is an irreducible representation of $G(n)$ (where $n\ge k$)
which varies as $n\to\infty$.

When the groups $G(n)$ are the symmetric groups $S(n)$ or the unitary
groups $U(n)$, Vershik and Kerov found necessary and sufficient conditions
on sequences 
$\{\pi_n\}$ under which the limit of the expression \tht{0.12} exists for
any $g\in\varinjlim G(k)$. 

It turns out that in the symmetric group case, for the normalized
character \tht{0.12} there exists an explicit formula in terms of
$s^*$-functions; this is a corollary of a new formula for dimension of
skew Young diagrams. In the case of $G(n)=U(n)$ an explicit formula in
terms of $s^*$-functions also exists provided group elements $g\in U(k)$
are replaced by elements of the algebra $\U(\gl(k))$. 
\footnote{Note that the idea 
to replace in \tht{0.12} the compact groups by the universal
enveloping algebras is present, in an implicit form, in Vershik--Kerov's note
\cite{VK2}.} 

The present paper is much obliged to Vershik--Kerov's work: in fact, our
initial aim was to analyze the sketch of proof of the main theorem in
\cite{VK2} about the characters of $U(\infty)$.

We plan to present a detailed exposition of this fundamental theorem (as
well as its generalization to other classical groups), based on the
machinery of $s^*$-functions, in next papers. 

\subhead 5. Main results \endsubhead

{\bf I. Binomial formula.} Let $\l=(\l\ge\ldots\ge\l_n)\in\Z^n$ be a
dominant weight for the group $GL(n,\C)$ and let $gl(n)_\l(z_1,\ldots,z_n)$
stand for the corresponding irreducible character, viewed as a function on
the subgroup of diagonal matrices. Then we have the following expansion
(Theorem 5.1) 
$$
\frac{gl(n)_\l(1+x_1,\ldots, 1+x_n)}
{gl(n)_\l(1,\ldots,1)}=
\sum_\mu\frac{s^*_\mu(\l_1,\ldots,\l_n)s_\mu(x_1,\ldots,x_n)}
{c(n,\mu)},\tag0.13
$$
where $\mu$ ranges over the set of partitions of length $\le n$ and
$c(n,\mu)$ are certain number factors not depending on $\l$. We call
\tht{0.13} the {\it binomial formula\/} for the (normalized) characters of
the group $GL(n)$.

In a very different form, the expansion \tht{0.13} can be found in Macdonald
\cite{M1}, Ch.I, Section 3, Example 10 (this example is due to Lascoux).
But the new point here is 
the observation that the 
binomial formula turns out to be related to the $s^*$-functions. 

The binomial formula plays an essential role in applications to the
asymptotic character theory. In the next paper we shall present similar
formulas for other classical groups. 

{\bf II. Dimension of skew Young diagrams.} Let $\mu\vdash k$ and
$\l\vdash n$ be two partitions, also viewed as Young diagrams. Let us
assume $k\le n$ and $\mu\subset\l$, and denote by $\dim\l/\mu$ the number
of standard tableaux of shape $\l/\mu$; in particular,
$\dim\l=\dim\l/\varnothing$. Recall that for $\dim\l$ nice explicit
formulas are known. Now we have
$$
\frac{\dim\l/\mu}{\dim\l}=
\frac{s^*_\mu(\l)}{n(n-1)\ldots(n-k+1)}\;,\tag0.14
$$
where $s^*_\mu(\l)=s^*_\mu(\l_1,\l_2,\ldots)$.

To our knowledge, this is a new formula, which is a quite surprising fact.

Formula \tht{0.14} can be applied to the asymptotic character theory of
the symmetric groups.

{\bf III. Higher Capelli identities.} We shall consider differential
operators with polynomial coefficients on the space $\Mnm$ of $n\times m$
matrices. Let $x_{ij}$ denote the natural coordinates in $\Mnm$ and let
$\p_{ij}$ be the corresponding partial derivatives. Let $\mu\vdash k$ be
an arbitrary partition of length $\le\min(n,m)$. We define the {\it higher
Capelli operator\/}, indexed by $\mu$, as the following differential
operator on $\Mnm$:
$$
\Dmnm=(k!)^{-1}
\sum_{i_1,\dots,i_k=1\,}^n 
\sum_{j_1,\dots,j_k=1\,}^m
\sum_{s\in S(k)} \chm(s)\cdot
x_{i_1 j_1} \dots x_{i_k j_k}
\p_{i_{s(1)} j_1} \dots \p_{i_{s(k)} j_k}\;,\tag0.15
$$
where $\chm(s)$ is the value of the irreducible character of $S(k)$,
indexed by $\mu$, at the permutation $s\in S(k)$.

These operators are invariant with respect to the action of the
groups $GL(n)$ and $GL(m)$ by left and right multiplications on $\Mnm$,
respectively. In the particular case $n=m$, $\mu=(1^n)$, the operator
\tht{0.15} reduces to the well--known Capelli operator
$$
\det\bmatrix x_{11}&\ldots&x_{1n}\\
\vdots&\phantom{\ldots}&\vdots\\
x_{n1}&\ldots&x_{nn}
\endbmatrix
\det\bmatrix \p_{11}&\ldots&\p_{1n}\\
\vdots&\phantom{\ldots}&\vdots\\
\p_{n1}&\ldots&\p_{nn}
\endbmatrix  \tag0.16
$$
Other Capelli operators are obtained when $n$ and $m$ are arbitrary and
$\mu=(1^k)$. (About Capelli operators and Capelli identities see  
Howe \cite{H}, Howe--Umeda \cite{HU}.) When $\mu=(k)$, we obtain
Capelli--type operators found by Nazarov \cite{N1}.

The action of the groups $GL(n)$ and $GL(m)$ on the space $\Mnm$ induces
two homomorphisms, $L$ and $R$, of the universal enveloping
algebras $\Ugln$ and $\Uglm$, respectively, in the algebra of differential
operators on $\Mnm$ with polynomial coefficients. Let us use a more detailed
notation $\Smn$ for the central elements $\S_\mu$ defined above. Then we
have the following result (Corollary 6.8)
$$
L(\S_{\mu|n})=R(\S_{\mu|m})=\Dmnm \tag0.17
$$
for all partitions $\mu$ of length $\le\min(n,m)$.

Together with the explicit formula for the quantum immanants $\Smn$,
obtained in  \cite{Ok1} (see also \cite{N2} and \cite{Ok2}) the relations
\tht{0.17} provide  {\it higher analogues\/} of the classical Capelli
identity. 

{\bf IV. Characterization Theorem for the $s^*$-functions.} An important
idea (already used in \cite{KO}) is to interpret elements of the algebra
$\Ls$ as functions $f(\l)=f(\l_1,\l_2,\ldots)$ on the set of partitions. 

\proclaim{\smc Characterization Theorem \cite{Ok1}} {\rm(See also Theorems
3.3 and 3.4 below.)} Fix a partition $\mu$. Then $s^*_\mu$ is the unique,
within a scalar factor, element of the algebra $\Ls$ such that
$s^*_\mu(\l)=0$ for all $\l\ne\mu$ with $|\l|\le|\mu|$.
\endproclaim

This characterization of $s^*$-functions turns out to be an efficient tool
for proving various results about the $s^*$-functions. Its role is
especially important in the proof of the identities \tht{0.17}.

{\bf V. Combinatorial presentation of $s^*$-functions.} Recall that the
ordinary Schur function $s_\mu$ admits a nice combinatorial description in
terms of tableaux. There exists a similar description for
the shifted Schur functions (Theorem 11.1):
$$
s^*_\mu(x_1,x_2,\ldots)=\sum_T\sum_{\a\in\mu}(x_{T(\a)}-c(\a)),\tag0.18
$$
summed over all {\it reverse tableaux\/} $T$ of shape $\mu$ and over all
boxes $\a$ of $\mu$, where $c(\a)$ is the content of the box $\a$ and in a
reverse tableau, in contrast to the conventional one, the entries {\it
decrease\/} left to right along each row (weakly) and down each column
(strictly). 

Note that \tht{0.8} and \tht{0.9} are particular cases of \tht{0.18}.

Formula \tht{0.18} can be easily derived from the combinatorial
presentation of the polynomials $t_\mu$ (see Chen--Louck \cite{CL} and
Macdonald \cite{M2}), but we also give an independent proof, based
on the Characterization Theorem.

{\bf VI. The coherence property of shifted Schur polynomials.} There is
one more stability property of the shifted Schur polynomials, which is
best stated in terms of the central elements $\Smn\in\Zgln$.

Note that for each $m=1,2,\ldots$ there exists a canonical projection
$$
\Uglm\to\Zglm,\tag0.17
$$
commuting with the adjoint representation. It turns out that if we apply
to $\Smn$ the $m$-th projection \tht{0.17}, where $m>n$, then the result
will be proportional to $\S_{\mu|m}$ (Theorem 10.1). We call this the {\it
coherence property\/}. When expressed in terms of the polynomials $s^*_\mu$
the coherence property leads to an interesting identity, see \tht{10.30} 
below.

We think the coherence property is an important argument in favor of the
thesis that the elements $\Smn$ constitute a distinguished basis of the
center. 

{\bf VII. Generating series for $h^*$- and $e^*$-functions.} The are nice
generating series for the complete symmetric functions $h_1,h_2,\ldots$
and elementary symmetric functions $e_1,e_2,\ldots$. In Theorem 12.1 we
present their analogues for $h^*_1,h^*_2,\ldots$ and $e^*_1,e^*_2,\ldots$.
An interesting feature of these new series is that they involve 
{\it inverse factorial powers\/} of the formal parameter.

{\bf VIII. Jacobi--Trudi formula.} For factorial Schur polynomials $t_\mu$
a Jacobi--Trudi--type formula was given in  Chen--Louck
\cite{CL}, Goulden--Hamel \cite{GH}, and Macdonald \cite{M2}, \cite{M1},
Ch. I, section 3, Examples 20--21.
However, this formula, being rewritten in terms of the polynomials
$s^*_{\mu|n}$, turns out to be not stable as $n\to\infty$ and so makes no
sense in the algebra $\Ls$. In Theorem 13.1
we obtain a different formula, which is stable and so expresses the
$s^*$-functions in terms of $h^*_1,h^*_2,\ldots$. There is also a dual
formula expressing $s^*_\mu$ through $e^*_1,e^*_2,\ldots$. Then a general
result due to Macdonald implies that there is a determinantal expression
of the $s^*$-functions in terms of the `hook' $s^*$-functions which is
just the same as in the classical Giambelli formula.

\subhead 6. Notes\endsubhead  First observations about the
$s^*$-functions were made  by Olshanski \cite{O3}. The 
results of the present work were the
subject of a talk at the 7-th Conference on formal power series and
algebraic combinatorics (Universit\'e Marne--la--Vall\'ee, May 29--June 2,
1995).  

\subhead 7. Acknowledgements \endsubhead We are grateful to Sergei Kerov,
Alain Lascoux, Alexander Molev, and Maxim Nazarov for discussions. We are
especially indebted to Ian Macdonald for his letter \cite{M3}, which contained
an important idea (see section 15 below), and to Alexander Postnikov for
the first version of the inversion formula \tht{14.9}.

The both authors were supported by the International Science Foundation
(under grants MBI300 and MQV000, respectively) and by the Russian Foundation 
for Basic Research under grant 95--01--00814.

\head 1. The algebra of shifted symmetric functions \endhead

Throughout the paper we assume that the ground field is $\C$ (although
most results hold over any field of characteristic zero).

Recall the definition of the algebra $\L$
of symmetric functions, see Macdonald \cite{M1}. Let $\L(n)$
denote the algebra of symmetric polynomials
in $x_1,\dots,x_n$. This algebra is graded
by degree of polynomials. The specialization
$x_{n+1}=0$ is a morphism of graded algebras
$$
\L(n+1)\to\L(n)\,.\tag1.1
$$
By definition $\L$ is the projective limit 
$$
\L=\varprojlim\L(n), \quad n\to\infty\,,
$$
in the category of graded algebras, taken with respect to morphisms 
\tht{1.1}. 
An element $f\in\L$ is by definition a
sequence $(f_n)_{n\ge1}$ such that:
\roster
\item $f_n\in\L(n)$, $n=1,2,\dots$,
\item $f_{n+1}(x_1,\dots,x_n,0)=f_{n}(x_1,\dots,x_n)$ 
(the stability condition),
\item $\sup_n\deg f_n<\infty$\,.
\endroster

Now let us denote by $\Ls(n)$ the algebra of polynomials in 
$x_1,\dots,x_n$ that become symmetric in new
variables
$$
x'_i=x_i-i+\const,\quad i=1,\dots,n \,. \tag1.2
$$
Here `$\const$' is an arbitrary fixed number;
note that the definition does not depend on
its choice. We call such polynomials
{\it shifted symmetric\/}. The algebra $\Ls(n)$
is filtered by degree of polynomials, and
the specialization 
$x_{n+1}=0$ is a morphism of filtered algebras
$$
\Ls(n+1)\to\Ls(n)\,.\tag1.3
$$

\definition{\smc Definition 1.1} Let
$$
\Ls=\varprojlim\Ls(n), \quad n\to\infty\,,\tag1.4
$$
be the projective limit in the category of filtered algebras, 
taken with respect to morphisms \tht{1.3}. We call $\Ls$ the {\it algebra 
of shifted symmetric functions}.
\enddefinition

Throughout this paper we use the notation
$$
(x\f k)=\cases x(x-1)\ldots (-k+1),&\text{if $k=1,2,\ldots$,}\\
1,&\text{if $k=0$,} \endcases \tag1.5
$$
for the $k$-th {\it falling factorial power\/} of a variable $x$. By
$\ell(\mu)$ we denote the {\it length\/} of a partition $\mu$.

\example{\smc Definition 1.2} Let $\mu=(\mu_1,\ldots,\mu_n)$ be a
partition, $\ell(\mu)\le n$. We define the {\it shifted Schur polynomial in
$n$ variables, indexed by $\mu$\/} as ratio of two $n\times n$
determinants, 
$$
s^*_\mu(x_1,\ldots,x_n)=\frac{\det[(x_i+n-i\f \mu_j+n-j)]}
{\det[(x_i+n-i\f n-j)]}, \tag1.6
$$
where $1\le i,j\le n$.
\endexample

Note that the denominator in \tht{1.6} equals the Vandermonde determinant in
the variables \tht{1.2}. Since the numerator is 
skew--symmetric in these variables, the ratio is indeed a
polynomial. Sometimes we will denote the polynomial
$s^*_\mu(x_1,\ldots,x_n)$ by $s^*_{\mu|n}$. Let us agree that
$s^*_{\mu|n}=0$ if 
$\mu$ is a partition with $\ell(\mu)>n$.

Let us show that the sequence $\{s^*_{\mu|n}\}$ defines an 
element of the algebra $\Ls$. It is
clear that $s^*_{\mu|n}$ is shifted symmetric and 
$\deg s^*_{\mu|n} = |\mu|$ for $n\ge\ell(\mu)$, so the degree is bounded
as $n\to\infty$. Let us verify the stability condition:

\proclaim{\smc Proposition 1.3} For each partition $\mu$ and each $n$
$$
s^*_\mu(x_1,\dots,x_n,0)=s^*_\mu(x_1,\dots,x_n)\,. \tag1.7
$$
\endproclaim
\demo{Proof}
By definition, 
$$
\gathered
s^*_\mu(x_1,\dots,x_{n+1})=
\frac
{\det\big[(x_i+n+1-i\f \mu_j+n+1-j)\big]}
{\det\big[(x_i+n+1-i\f n+1-j)\big]}, \\
\text{where $i,j=1,2,\ldots,n+1$.}
\endgathered\tag1.8
$$
First suppose $\ell(\mu)\le n$;  that is $\mu_{n+1}=0$.
Then substituting $x_{n+1}=0$ in the numerator of \tht{1.8}
we obtain 
$$
\det
\left[
\matrix
(x_1+n\f \mu_1+n)&\hdots&(x_1+n\f \mu_n+1)&1\\
\vdots&&\vdots&\vdots\\
(x_n+1\f \mu_1+n)&\hdots&(x_n+1\f \mu_n+1)&1\\
0&\hdots&0&1
\endmatrix
\right] \,.
$$
Clearly this determinant equals
$$
(x_1+n)\dots(x_n+1)
\det
\left[
\matrix
(x_1+n-1\f \mu_1+n-1)&\hdots&(x_1+n-1\f \mu_n)\\
\vdots&&\vdots\\
(x_n\f \mu_1+n-1)&\hdots&(x_n\f \mu_n)
\endmatrix
\right] \,.
$$
Likewise, the denominator of \tht{1.8} becomes
$$
(x_1+n)\dots(x_n+1)
\det
\big[
(x_i+n-i\f n-j)
\big],\quad i,j=1\dots n \,.
$$
This yields \tht{1.7} provided $\ell(\mu)\le n$. 
If $\ell(\mu)=n+1$ then 
substituting $x_{n+1}=0$ in the numerator of \tht{1.8}
we obtain 
$$
\det
\left[
\matrix
(x_1+n\f \mu_1+n)&\hdots&(x_1+n\f \mu_n+1)&0\\
\vdots&&\vdots&\vdots\\
(x_n+1\f \mu_1+n)&\hdots&(x_n+1\f \mu_n+1)&0\\
0&\hdots&0&0
\endmatrix
\right] \,.
$$
This determinant clearly vanishes
and hence the left--hand side of \tht{1.7} vanishes.
On the other hand the right--hand side equals zero
by definition. If $\ell(\mu)>n+1$ then
both sides equal zero by definition.
This completes the proof. \qed
\enddemo

\definition {\smc Definition 1.4} By Proposition 1.3, for each partition
$\mu$ the sequence \tht{1.6}, where $n\to\infty$,  
defines an element of the algebra $\Ls$. We denote it by $s^*_\mu$ and 
call it the {\it shifted Schur
function}, indexed by $\mu$.  These functions will also be called  
$s^*$-{\it functions\/} for short.
\enddefinition

Note that the shift of variables
\tht{1.2} establishes an isomorphism between
$\Ls(n)$ and $\L(n)$ but no shift of
variables maps $\Ls$ to $\L$. However
$\Ls$ and $\L$ are related in the following
way:

\proclaim{\smc Proposition 1.5} The graded algebra $\gr\Lambda^*$
corresponding to the filtered algebra $\Lambda^*$ is canonically
isomorphic to the algebra $\Lambda$.
\endproclaim

\demo{Proof} Note that if $f\in\Lambda^*(n)$ is of degree $\le d$ then its
$d$-th homogeneous component is a symmetric polynomial. It follows that the
algebras $\gr\Lambda^*(n)$ and $\Lambda(n)$ are canonically isomorphic for
each $n$. Since the isomorphisms $\gr\Lambda^*(n)\to\Lambda(n)$ are
compatible with the specialization maps $x_n=0$, we obtain a canonical
algebra isomorphism $\gr\Lambda^*\to\Lambda$.\qed
\enddemo

If $f$ is an element of $\Ls$ (resp., of $\Ls(n)$) of degree $d$ then its
image in the $d$-th homogeneous component of the graded algebra $\L$
(resp., $\L(n)$) will be called the {\it highest term\/} of $f$. 

Note that the highest term of $s^*_\mu\in\Lambda^*$ is the ordinary Schur
function $s_\mu\in\Lambda$ and the highest term of $s^*_{\mu|n}\in\Ls(n)$
is $s_{\mu|n}\in\L(n)$, the Schur polynomial in $n$ variables. This
follows at once from the comparison of \tht{1.6} with the well--known
expression for $s_{\mu|n}$,
$$
s_\mu(x_1,\ldots,x_n)=\frac{\det[x_i^{\mu_j+n-j}]}
{\det[x_i^{n-j}]}.\tag1.9
$$

By definition, put
$$
\align
h^*_k&=s^*_{(k)}, \quad k=1,2,\dots \tag1.10\\
e^*_k&=s^*_{(1^k)}, \quad k=1,2,\dots \,.\tag1.11
\endalign
$$
These are shifted analogues of the complete homogeneous symmetric
functions and the elementary symmetric functions. 

Put also
$$
\ps_k = \sum_i\big((x_i-i)^k-(-i)^k\big)\,.\tag1.12
$$
These are certain analogues of Newton power sums, which appeared in  
Olshanski's papers \cite{O1} and \cite{O2}.

\proclaim{\smc Corollary 1.6} 
\roster
\item The shifted Schur functions $\{s^*_\mu\}$ form a linear
basis in $\Ls$.
\item The algebra $\Lambda^*$ is the algebra of
polynomials in $h^*_1,h^*_2,\dots$ or in $e^*_1,e^*_2,\dots$.
\item The algebra $\Lambda^*$ is the algebra of
polynomials in $\ps_1,\ps_2,\dots$.
\endroster
\endproclaim

\demo{Proof} Immediately follows from Proposition 1.5  and the similar
well--known claims for the algebra $\Lambda$. \qed
\enddemo

Note also that for any fixed $n$ the shifted Schur polynomials
$s^*_{\mu|n}$ form a linear basis in $\Ls(n)$.
\example{\smc Remark 1.7} The algebra $\Lambda^*$ also may be regarded as a
deformation of the algebra $\Lambda$. Indeed, let $\theta$ be
a number parameter. For each $n=1,2,\dots$ let
$\Lambda^*_\theta(n)$ be the algebra of polynomials in
$x_1,\dots,x_n$ which become symmetric in variables
$$
x_i'=x_i+\const-i\theta,\quad i=1,\dots,n\,,\tag1.13
$$
and define $\Lambda^*_\theta=\varprojlim\Lambda^*_\theta(n)$.
Then $\Lambda^*_1=\Ls$ and $\Lambda^*_0=\L$.
Note that the scaling
$$
x_i\mapsto x_i/\theta\tag1.14
$$
establishes an isomorphism $\L^*_\theta\cong\Ls$
for all nonzero $\theta$.
\endexample

\head 
2. Quantum immanants
\endhead

Throughout the paper we will use the following notation:

$\gln$ is the general linear Lie algebra $\gl(n,\C)$,

$E_{ij}$ are the standard generators of $\gl(n)$ --- the matrix units,

$\Ugln$ is the universal enveloping algebra of $\gln$,

$\Zgln$ is the center of $\Ugln$,

$S(\gln)$ is the symmetric algebra of $\gln$,

$GL(n)$ is the general linear group $GL(n,\C)$.

\medskip

Recall the construction of the Harish-Chandra isomorphism for the case of 
$\gln$, see, e.g., Bourbaki \cite{Bou}, Ch. VIII, 8.5, or Dixmier
\cite{D}, 7.3.

Suppose $X\in\Zgln$. Given $\l=(\l_1,\ldots,\l_n)\in\C^n$, we consider an
arbitrary (cyclic) highest weight $\gln$-module with highest weight $\l$
(relative to the upper triangular Borel subalgebra). Then $X$ acts in
this module as a scalar operator, say $f_X(\l)\operatorname{id}$. (Note
that $f_X(\l)$ does not depend on the choice of the module.) The
assignment 
$$
X\mapsto f_X(\cdot)
$$
is an isomorphism of the algebra $\Zgln$ onto the algebra of polynomials
in $\l\in\C^n$ that are symmetric in the coordinates of 
$\l'=\l+\rho$, where
$\rho$ stands for the half--sum of positive roots. This isomorphism is
called the {\it Harish--Chandra isomorphism\/}. 

Equivalently, the Harish--Chandra isomorphism can be defined by making use
of the projection 
$$
\multline
\Ugln=(\frak{n}_-\Ugln+\Ugln\frak{n}_+)\oplus\Cal{U}(\frak{h})\mapsto\\
\mapsto\Cal{U}(\frak{h})=S(\frak{h})=\C[\l_1,\ldots,\l_n],
\endmultline \tag2.1
$$
where $\frak{n}_+$ and $\frak{n}_-$ are the upper and lower triangular
nilpotent subalgebras of $\gln$ and $\frak{h}$ is the diagonal subalgebra.

Note also that
$$
\deg X = \deg f_X
$$
where $\deg X$ is the degree of $X$ with respect to
the natural filtration in $\Ugln$.

\proclaim{\smc Proposition 2.1} The Harish--Chandra isomorphism is an
algebra isomorphism
$$
\Zgln \to \Ls(n).\tag2.2
$$
\endproclaim

\demo{Proof} This follows at once from the definitions, because the shift
$\l\mapsto\l'=\l+\rho$ is of the form \tht{1.2}. (We recall that 
$$
\rho=(\tfrac{n-1}2,\tfrac{n-3}2,\ldots, -\tfrac{n-3}2,-\tfrac{n-1}2).
$$
 \qed
\enddemo

\example{\smc Example 2.2} Take as $X$ the Casimir element
$$
C=\sum_{i,j} E_{ij} E_{ji} \in \Zgln \,.
$$
Since
$$
C=\sum_i E_{ii}^2 + 2\sum_{i>j} E_{ij} E_{ji}
+\sum_{i<j} (E_{ii}-E_{jj}),
$$
its image under the projection \tht{2.1} is
$$
\sum_i E_{ii}^2 + \sum_{i<j} (E_{ii}-E_{jj}),
$$
whence
$$
\gather
f_C(\l)=\sum\l_i^2+\sum_{i<j}(\l_i-\l_j)=
\sum(\l_i^2+(n+1-2i)\l_i)\\
=\sum((\l_i+\tfrac{n+1}2 -i)^2-(\tfrac{n+1}2 -i)^2),
\endgather
$$
which is a shifted symmetric polynomial in $\l$.
\endexample

\definition{\smc Definition 2.3} By virtue of Proposition 2.1, for each
partition $\mu$ with $\ell(\mu)\le n$ there exists a central element
$$
\Sm=\Smn\in\Zgln \tag2.3
$$
corresponding to the shifted Schur polynomial $s^*_{\mu|n}\in\Ls(n)$, i.e.,

$\Smn$ is defined by 
$$
f_{\Sm}(\l_1,\ldots,\l_n)=s^*_\mu(\l_1,\ldots,l_n),
\qquad (\l_1,\ldots,\l_n)\in\C^n.\tag2.4
$$
We will call $\Sm$ the {\it quantum $\mu$-immanant\/} (an explanation of
this term will be given below in Remark 2.6).
\enddefinition

Let us calculate the highest term of $\Sm$ with respect to
the natural filtration in $\Ugln$. Recall that there is a canonical
isomorphism
$$
\gr\,\Ugln\,\cong\,S(\gln)\,,\tag2.5
$$
and denote by $I(\gln)$ the subalgebra of invariants in $S(\gln)$ under
the adjoint action of the group $GL(n)$. We have
$$
\gr\Zgln\,\cong\,I(\gln)\,.\tag2.6
$$
By means of the basis $\{E_{ij}\}$ of $\gln$ we can identify $\gln$ with its
dual space (under this identification each $E_{ij}$ becomes the coordinate
function $x_{ij}$). Then $S(\gln)$ can be identified with the algebra
$\C[\gln]$ of polynomial functions on $\gln$ and $I(\gln)$ turns into the
algebra of invariant polynomial functions.

Next, remark that the algebra $I(\gln)$ is isomorphic to the algebra
$\L(n)$: the isomorphism
$$
I(\gln)\,\cong\,\L(n) \tag2.7
$$
is simply the restriction of invariant polynomials
to the diagonal subalgebra of $\gln$ (the Chevalley restriction map, see
\cite{Bou}, Ch. VIII, 8.3, or \cite{D}, 7.3). It is well--known that the
Harish--Chandra homomorphism 
and the Chevalley restriction map are compatible within lower terms,
i.e., the following
diagram is commutative
$$
\CD
\gr\Zgln @>>> I(\gln) \\
@VVV @VVV \\
\gr\Ls(n) @>>> \L(n) \,, 
\endCD \tag2.8
$$
where the bottom arrow in \tht{2.8} is the canonical isomorphism mentioned
in Proposition 1.5.

Finally, let us introduce the element
$$
S_\mu=S_{\mu|n}\in I(\gln)\subset S(\gln), \tag2.9
$$
which corresponds to the Schur polynomial $s_{\mu|n}$ under the
isomorphism \tht{2.7}. If we consider $S_\mu$ as a polynomial function
$S_\mu(X)$ 
on $\gln$ then $S_\mu(X)$ is simply the Schur polynomial $s_\mu$ in the
eigenvalues of the matrix $X$. 

Now look at the commutative diagram \tht{2.8}. Since the highest term of
$s^*_\mu$ equals $s_\mu$ we see that $S_\mu$ is the highest term of $\Sm$.

In the next proposition we present an explicit formula for $S_\mu$. Note
that this is essentially the well--known definition of the Schur function
via the characteristic map, see Macdonald \cite{M1}, Ch. I, section 7.

\proclaim{\smc Proposition 2.4} Let $\mu\vdash k$ be a partition,
$\ell(\mu)\le n$, and let $\chm$ denote the corresponding irreducible
character of the symmetric group $S(k)$. Then we have
$$
S_\mu=(k!)^{-1}\sum_{\,\,i_1,\dots,i_k=1\,\,}^n \sum_{s\in S(k)}
\chm(s)\, E_{i_1,i_{s(1)}} \dots E_{i_k,i_{s(k)}} \, \tag2.10 
$$
or, as a polynomial invariant,
$$
S_\mu(X)=(k!)^{-1}\sum_{\,\,i_1,\dots,i_k=1\,\,}^n \sum_{s\in S(k)}
\chm(s)\, x_{i_1,i_{s(1)}} \dots x_{i_k,i_{s(k)}} \,, \tag2.11 
$$
where $X$ is a $n \times n$ matrix and $x_{ij}$ are its entries.
\endproclaim

\proclaim{\smc Corollary 2.5} We have
$$
\Sm=(k!)^{-1}\sum_{\,\,i_1,\dots,i_k=1\,\,}^n \sum_{s\in S(k)}
\chm(s)\, E_{i_1,i_{s(1)}} \dots E_{i_k,i_{s(k)}} 
+\text{\rm lower terms} \tag2.12
$$
\endproclaim

\demo{Proof of Proposition 2.4} Let $V$ denote the space $\C^n$ considered
as a natural 
$\gln$-module and let $V_\mu$ be the irreducible
polynomial $\gln$-module, indexed by $\mu$. Both $V$ and $V_\mu$ may also
be regarded as modules over the semigroup $\M(n)$ of all $n\times n$ matrices.
Since the Schur polynomials coincide with the characters of the
irreducible polynomial modules, we have
$$
S_\mu(X)=\tr_{V_\mu}(X), \qquad X\in\M(n).\tag2.13
$$

On the other hand, by the Schur--Weyl duality, $V_\mu$ occurs in the
decomposition of the module $V^{\otimes k}$ with multiplicity $\dim\mu$,
where $\dim\mu=\chm(e)$ is the dimension of the character $\chm$. It
follows 
$$
S_\mu(X)=\frac 1{\dim\mu}\tr_{V^{\otimes k}}(X^{\otimes k}\cdot P_\mu),
\tag2.14
$$
where $P_\mu$ stands for the projection in the tensor space $V^{\otimes
k}$ onto the isotypic component of $V_\mu$. This projection is given
by the central idempotent
$$
(1/k!)\dim\mu\sum_{s\in S(k)} \chm(s)\cdot s\in \C[S(k)],\tag2.15
$$
whence
$$
\align
s_\mu(X)
&=\tr (X^{\otimes k}\cdot\chm/k!)   \\
&=(1/k!)\sum_{\,\,i_1,\dots,i_k=1\,\,}^n \sum_{s\in S(k)}
\chm(s)\, x_{i_1,i_{s(1)}} \dots x_{i_k,i_{s(k)}} \,,
\endalign
$$
which proves \tht{2.11}. \qed
\enddemo

\example{\smc Remark 2.6}
Given a matrix $A=[a_{ij}]$, $i,j=1,\dots,k$,
the number
$$
\sum_{s\in S(k)} \chm(s) \, a_{1,s(1)} \dots a_{k,s(k)}
$$
is called the {\it $\mu$-immanant} of the matrix A, see Littlewood 
\cite{L}, 6.1.
If $\mu=(1^k),(k)$ then the $\mu$-immanant turns into
determinant and permanent respectively.
Note that \tht{2.11} expresses the invariant polynomial $S_\mu(X)$ 
as sum of $\mu$-immanants of principal
$k$-submatrices (possibly 
with repeated rows and columns) of the matrix $X$. This
was one of the reasons why the central element $\Sm$ (which may be viewed as
a `quantum analogue' of $S_\mu$) was called the `quantum immanant'.
Another reason is that in the special case $\mu=(1^k)$, $k=1,2,\ldots,n$,
the corresponding central elements appear when expanding the
so--called quantum determinant for the Yangian $Y(\gln)$ (about this, see 
Nazarov \cite{N1} and Molev--Nazarov--Olshanski \cite{MNO}).
\endexample

To find an explicit formula for $\Sm$ is a much less trivial problem.
Such a formula was obtained in \cite{Ok1}, \cite{N2}, and \cite{Ok2}; 
see also Remark 14.6 and section 15.

By virtue of the Harish-Chandra isomorphism, properties
of $s^*$-functions can be interpreted as properties of 
quantum immanants and vice versa. For example, Proposition 1.3 
asserts a kind of stability for quantum
immanants. This stability will be discussed below in
Remark 6.9.

\head 3. Characterization of $s^*$-functions \endhead

For any element $f\in\Ls$ its value $f(x_1,x_2,\ldots)$ on an arbitrary 
infinite
sequence $(x_1,x_2,\ldots)$ is well--defined provided  $x_i=0$
as $i$ is large enough. In particular, the value $f(\l)$ exists for any
partition $\l$; moreover, $f$ is uniquely determined by its values on all
partitions. Thus, the algebra $\Ls$ can be realized as a certain function
algebra on the set of partitions. 
\footnote{ This realization has been considered by Kerov and Olshanski
\cite{KO}.}
This new point of view (usually
partitions appear rather as parameters than arguments) turns out to be
extremely fruitful.

In this section we review some results of Okounkov's paper \cite{Ok1} 
about the $s^*$-functions viewed as functions on partitions. For
completeness we give proofs. 

Let $\mu$ be a partition, also viewed as a
Young diagram.
Denote by $H(\mu)$ the product of the hook lengths of
all boxes of $\mu$,
$$
H(\mu)=\prod_{\alpha\in\mu} h(\alpha) \,.\tag3.1
$$
Let $\l$ be another partition. Write
$\mu\subset\l$ if $\mu_i\le\l_i$ for all $i=1,2,\dots$. 

\proclaim{\smc Theorem 3.1 (Vanishing Theorem \cite{Ok1})} We have
$$
\align
s^*_\mu(\l)&=0\quad\text{unless}\quad\mu\subset\lambda\,,\tag3.2\\
s^*_\mu(\mu)&=H(\mu)\,.\tag3.3
\endalign
$$
\endproclaim

\demo{Proof}
Observe that
$$
(a\f b)=0\quad\text{if}\quad a,b\in\Z_+, b>a\,.
$$
Suppose $\l_l<\mu_l$ for some $l$ and choose an arbitrary
$n\ge\max(\ell(\mu),\ell(\l))$. Then in the $n\times n$ matrix
$$
\big[(\l_i+n-i\f \mu_j+n-j)\big]\tag3.4
$$
all entries with $i\ge l$ and $j\le l$ vanish. Let us expand 
the determinant of the matrix \tht{3.4} into an alternating sum of
monomials in the entries of \tht{3.4};
 then  each monomial will contain at least one entry with 
$i\ge l$ and $j\le l$. Hence the determinant of \tht{3.4} vanishes. Since the 
denominator in \tht{1.6} does not vanish when $(x_1,\ldots,x_n)$ 
is a partition
$(\l_1,\ldots,\l_n)$, \tht{3.2} follows.

Now take $\l=\mu$, and let $n\ge\ell(\mu)$. Then in the matrix
$$
\big[(\mu_i+n-i\f \mu_j+n-j)\big]\tag3.5
$$
all entries with $i>j$ vanish. Hence the determinant of \tht{3.5} equals
$$
\prod_i(\mu_i+n-i)! \,.
$$
Therefore
$$
s^*_\mu(\mu)=
\prod_i(\mu_i+n-i)! \,/\,
\prod_{i<j}(\mu_i-\mu_j-i+j) \,. \tag3.6
$$
Recall that there are two formulas for $\dim\mu$ (the dimension
of the irreducible representation of the symmetric group 
indexed by $\mu$),
$$
\align
\dim\mu&=|\mu|!\,/\,H(\mu) \quad\text{(the hook formula)}\tag3.7\\
&=|\mu|!
\prod_{i<j}(\mu_i-\mu_j-i+j) \,/\,
\prod_i(\mu_i+n-i)! \,.\tag3.8
\endalign
$$
Thus \tht{3.6} equals $H(\mu)$. \qed
\enddemo

\proclaim{\smc Theorem 3.2 (Characterization Theorem I\,\cite{Ok1})}
The function $s^*_\mu$ is the unique element of $\Ls$
such that $\deg s^*_\mu\le|\mu|$ and
$$
s^*_\mu(\l)=\dt_{\mu\l} H(\mu)
$$
for all $\l$  such that $|\l|\le|\mu|$.
\endproclaim

This follows from a more strong claim:

\proclaim{\smc Theorem 3.3 (Characterization Theorem I${}'$\,\cite{Ok1})}
Suppose $\ell(\mu)\le n$.
Then the polynomial $s^*_{\mu|n}$ is the unique element of $\Ls(n)$
such that $\deg s^*_{\mu|n}\le|\mu|$ and
$$
s^*_{\mu|n}(\l)=\dt_{\mu\l} H(\mu)
$$
for all $\l$ such that $|\l|\le|\mu|$ and $\ell(\l)\le n$.
\endproclaim

\demo{Proof}
We have to prove that if $f\in\Ls(n)$, $\deg f\le|\mu|$,
and $f(\l)=0$ for all $\l$ such that 
$|\l|\le|\mu|$, $\ell(\l)\le n$, then $f=0$. Put
$k=|\mu|$. The polynomials $\{s^*_{\l|n}\}$, where $|\l|\le k$,
$\ell(\l)\le n,$ constitute a linear basis in the space of shifted symmetric
polynomials in $n$ variables of degree $\le k$. Hence
$$
f=\sum c_\l s^*_{\l|n}, \qquad |\l|\le k,\quad\ell(\l)\le n\,, \tag3.9
$$
for some coefficients $c_\l$. Let us show that $c_\l=0$ for
all $\l$. Suppose there are partitions $\nu$ such that $c_\nu\ne0$.
Choose such a partition $\nu$ so that
$c_\nu\ne0$ and $c_\eta=0$ for all $\eta$, $|\eta|<|\nu|$, and 
evaluate \tht{3.9} at $\nu$. By the Vanishing Theorem we obtain
$$
0=c_\nu H(\nu) \,.
$$
Thus $c_\nu=0$, which leads to contradiction. \qed
\enddemo
 
There exists a slightly different version of the Characterization Theorem: 

\proclaim{\smc Theorem 3.4 (Characterization theorem II \,\cite{Ok1})}
The function $s^*_\mu$ is the unique element of $\Ls$
such that the highest term of $s^*_\mu$ is the 
ordinary Schur function $s_\mu$  and
$$
s^*_\mu(\l)=0 
$$
for all $\l$ such that $|\l|<|\mu|$.
\endproclaim
\demo{Proof}
Suppose there are two such elements $f_1$ and $f_2$ of $\Ls$.
Then $\deg(f_1-f_2)<|\mu|$ and
$(f_1-f_2)(\l)=0$ for all $\l$ such that $|\l|<|\mu|$.
By Characterization Theorem I we have $f_1-f_2=0$. \qed
\enddemo

Note that there is an obvious analogue of Theorem 3.4 for the polynomials
$s^*_{\mu|n}$. 

\example{\smc Example 3.5} As a first application of the Characterization
Theorem let us give a new proof of Proposition 1.3.

Suppose $\ell(\mu)\le n+1$ and consider the polynomial
$$
f=s^*_\mu(x_1,\dots,x_n,0)\,. \tag3.10
$$
Clearly it is shifted symmetric. By the Vanishing Theorem 
$$
\alignat2
f(\l)&=0&\quad&\text{  for all }\l\text{ such that }
\mu\not\subset\l,\ell(\l)\le n\,, \\
f(\mu)&=H(\mu)&\quad&\text{  if }\ell(\mu)\le n \,.
\endalignat 
$$ 
Hence by Characterization Theorem I$'$
$$
f=\cases
s^*_\mu(x_1,\dots,x_n), &\ell(\mu)\le n ,\\
0, &\ell(\mu) = n+1\,.
\endcases
$$
\endexample

\example{\smc Remark 3.6} The Vanishing and Characterization Theorems are
a way to control lower terms of the inhomogeneous polynomials
$s^*_{\mu|n}$. In particular cases similar arguments were used by several
people (see, e.g. Howe--Umeda \cite{HU}). In full generality these
arguments were developed by Sahi \cite{Sa1} (we learned about this
important paper 
after work on the present paper was completed). For an arbitrary 
fixed decreasing sequence of real numbers $\rho=(\rho_1,\ldots,\rho_n)$,
Sahi established existence and uniqueness of symmetric polynomials 
$p_\mu(y_1,\ldots,y_n)$ such that
$\operatorname{deg}p_\mu\le|\mu|$, $p_\mu(\lambda+\rho)=0$ 
when $|\lambda|\le|\mu|$, $\lambda\ne\mu$, and $p_\mu(\mu+\rho)\ne0$.
Sahi also described an inductive procedure to construct the 
polynomials $p_\mu$.

Theorem 3.3 shows that in the particular case $\rho_i-\rho_{i+1}=1$
Sahi's polynomials $p_\mu$ reduce to the factorial Schur polynomials
$t_\mu(y_1,\ldots,y_n)$ and so admit a nice closed expression. (It seems
difficult to see this directly from Sahi's general construction.)
\endexample

\head 4. Duality \endhead

Given a partition $\mu$, we denote by $\mu'$ the dual partition (i.e., the
transposed Young diagram). Recall that in the algebra $\L$ of symmetric
functions there exists an involutive automorphism $\o$ such that
$$
\o(s_\mu)=s_{\mu'} \qquad\text{for all $\mu$},\tag4.1
$$
and
$$
\o(p_k)=(-1)^{k-1}p_k \qquad\text{for $k=1,2,\ldots$,} \tag4.2
$$
see \cite{M1}, Ch. I, (2.13) and (3.8).

In this section section we aim to construct a similar automorphism for the
algebra $\Ls$. As we will see it admits a nice interpretation 
when the algebra
$\Ls$ is realized as an algebra of functions on the set of partitions.

Consider the elements $\ps_k$ introduced in \tht{1.12} and combine them into
the following generating series
$$
\Ps(u)=\sum_{k>0} \ps_k u^{-(k+1)}\,. \tag4.3
$$
We have
$$
\align
\Ps(u)
&=\sum_i \sum_{k>0} \big((x_i-i)^k-(-i)^k)\big) u^{-(k+1)}\\
&=\sum_i 
\left(
\frac{1}{u-x_i+i}
-\frac{1}{u+i}
\right) \\
&=
\frac{d}{du}\log
\prod_i
\frac{u-x_i+i}{u+i} \,.\tag4.4
\endalign
$$
Define an endomorphism $\o$ of the algebra $\Ls$ on
the generators $\ps_k$ of $\Ls$ by
$$
\o\Ps(u) = \Ps(-u-1)\,, \tag4.5
$$
where  
$$
\o\Ps(u)=\sum_{k>0} \o(\ps_k) u^{-(k+1)}\,. \tag4.6
$$
Remark that \tht{4.5} makes sense because a formal power
series in $(u+1)^{-1}$ can be reexpanded as a 
formal power series in $u^{-1}$. The mapping
$$
u\mapsto -u-1
$$
is involutive. Hence $\o$ is an involutive
automorphism of $\Ls$. We have
$$
\align
\o\Ps(u)
&=\sum_i 
\left(
\frac{1}{-u-(x_i+1-i)}
-\frac{1}{-u-(1-i)}
\right) \\
&=\sum_i \sum_{k>0} 
(-1)^{k-1}\big((x_i+1-i)^k-(1-i)^k)\big) u^{-(k+1)} \,,
\endalign
$$
whence
$$
\o(\ps_k)=
(-1)^{k-1}
\sum_i \big((x_i+1-i)^k-(1-i)^k)\big). \tag4.7
$$
This is an analogue of \tht{4.2}.

\proclaim{\smc Theorem 4.1} Suppose $\l$ is a partition.
Then
$$
[\o(f)](\l)= f(\l') \tag4.8
$$
for all $f\in\Ls$.
\endproclaim
\demo{Proof}
Introduce two factorial analogues of the functions $\ps_k$
$$
\ph_k(x)=\sum_i \big((x_i-i\u k)-(-i\u k)\big)\,, \tag4.9
$$
where $(x\u k)=x(x+1)\ldots(x+k-1)$ stands for the $k$-th {\it raising\/}
factorial power of $x$, and
$$
\pc_k(x)=\sum_i \big((x_i+1-i\f k)-(1-i\f k)\big)\,.\tag4.10
$$
Let us expand $(x\u k)$ in ordinary powers of $x$,
$$
(x\u k) = \sum_{l\le k} c_{kl}\, x^l\,, \quad c_{kl}\in\Z. \tag4.11
$$
It is clear that
$$
(x\f k) = \sum_{l\le k} (-1)^{k-l} c_{kl}\, x^l\,. \tag4.12
$$
The numbers $(-1)^{k-l} c_{kl}$  are known as Stirling numbers
of the first kind. We have
$$
\align
\o(\ph_k)
&=\sum_l c_{kl}\, \o(\ps_l) \\
&=\sum_l c_{kl} 
(-1)^{l-1} \sum_i \big((x_i+1-i)^l-(1-i)^l)\big)\\
&=(-1)^{k-1} \sum_l (-1)^{k-l} c_{kl} 
\sum_i \big((x_i+1-i)^l-(1-i)^l)\big)\\
&=(-1)^{k-1} \pc_k \,. \tag4.13
\endalign
$$
Since the functions $\ph_k$ are generators of $\Ls$
it suffices to check that
$$
\pc_k(\l)=(-1)^{k-1}\ph_k(\l') \tag4.14
$$
for all $\l$ and all $k>0$. In order to do this let us calculate the sum
$$
k \sum_{(i,j)\in\l} (j-i\f k-1) \, \tag4.15
$$
in two ways, by making use of the following elementary fact 
$$
k \sum_{l=a}^b (l\f k-1) = (b+1\f k) - (a\f k)\,. \tag4.16
$$

First sum \tht{4.15} along the rows. Then we obtain
$$
\sum_i \big((\l_i+1-i\f k)-(1-i\f k)\big) = \pc_k(\l) \,.
$$
Next sum \tht{4.15} along the columns. Then we obtain
$$
\multline
\sum_j \big((j\f k)-(j-\l'_j\f k)\big)  \\
=(-1)^{k-1}\sum_j \big(-(-j\u k)+(\l'_j-j\u k)\big) = 
(-1)^{k-1}\ph_k(\l')\,.
\endmultline
$$
This proves \tht{4.14} and the theorem. \qed
\enddemo

Note that the trick with two way summation of \tht{4.15}
was used in \cite{KO} in the proof of another
identity involving $\pc(\l)$\,.

\proclaim{\smc Theorem 4.2} For each partition $\mu$
$$
\o(s^*_\mu)=s^*_{\mu'}\,. \tag4.17
$$
\endproclaim

\demo{Proof}
This immediately follows from Theorem 4.1, Theorem 3.2, 
and the following obvious equality
$$
H(\mu)=H(\mu') \,. \qed
$$
\enddemo

This is an exact analogue of \tht{4.1}. As for ordinary symmetric functions,
we have
$$
\o(h^*_k)=e^*_k,\quad \o(e^*_k)=h^*_k\qquad \text{for $k=1,2,\ldots$.}
\tag4.18 
$$

\head 5. Binomial Theorem \endhead

The aim of this section is to write a Taylor--type expansion of the
irreducible characters of $GL(n)$ at $1\in GL(n)$ in terms of 
$s^*$-functions.

Introduce some notation. Let $\mu$ be a Young diagram and $\a=(i,j)$ be a
box of $\mu$. 
The number
$$
c(\a)=j-i \tag5.1
$$
is called the {\it content} of $\a$. Put
$$
(n\u\mu)=\prod_{\a\in\mu}(n+c(\a))\,. \tag5.2
$$
This is a generalization both of raising and falling factorial powers.
Indeed, 
$$
(n\u(k))=n(n+1)\ldots(n-k+1)\quad \text{and} \quad
(n\u (1^k))=(n\f k).
$$
Clearly
$$
(n\u\mu)=\prod_{i}(\mu_i+n-i)!\,/\,(n-i)!\,.\tag5.3
$$

By a {\it signature\/} we mean an arbitrary highest weight for the group
$GL(n)$, i.e., a $n$-tuple
$$
\l=(\l_1,\ldots,\l_n)\in\Z^n,
\qquad\text{such that $\l_1\ge\ldots\ge\l_n$.} 
$$
To each signature corresponds an irreducible finite--dimensional
$GL(n)$-module; let $gl(n)_\l(X)$ denote its character and
$\dimn\l$ its dimension. We also consider  $gl(n)_\l$ 
as a (rational) function
$gl(n)_\l(z_1,\ldots,z_n)$ on the torus of diagonal matrices in $GL(n)$.
When $\l_n\ge0$, the character $gl(n)_\l$ coincides with the polynomial
function $S_\l$ on $n\times n$ matrices, and $gl(n)_\l(z_1,\ldots,z_n)$
becomes the Schur polynomial
$s_{\l|n}$ (in the general case $gl(n)_\l(z_1,\ldots,z_n)$ is often called
the rational Schur 
function). The character $gl(n)_\l(z_1,\ldots,z_n)$ is given by exactly 
the same formula as
$s_{\l|n}$:
$$
gl(n)_\l(z_1,\ldots,z_n)=
\frac{\det[z_i^{\l_j+n-j}]}
{\det[z_i^{n-j}]},\quad 1\le i,j\le n; \tag5.4
$$
this is Weyl's character formula for $GL(n)$.
For the dimension $\dimn\l$ there are two useful expressions:
$$
\align
\dimn\l&=(n\u\l)\ov H(\l) \tag5.5\\
&=\prod_{i<j}(\l_i-\l_j+j-i)\ov\prod_{i}(n-i)!  \tag5.6
\endalign
$$
(recall that $H(\l)$ is the product of hook lengths in $\l$, see
\tht{3.1}). 
The first expression is called the {\it hook formula\/} for $GL(n)$, 
and the second one is Weyl's
dimension formula.

\proclaim{\smc Theorem 5.1 (Binomial Theorem for $GL(n)$)} Let
$\l=(\l_1,\ldots,\l_n)$ be a signature for $GL(n)$. Then
$$
\dfrac{gl_\l(1+x_1,\ldots,1+x_n)}{\dimn\l}= \sum_\mu
\frac1{(n\u\mu)}\, s^*_\mu(\l_1,\ldots,\l_n) \, s_\mu(x_1,\ldots,x_n)\,
\tag5.7
$$
\endproclaim

Recall that the symbol $(n\u\mu)$ is defined by \tht{5.2} and \tht{5.1}

\example{\smc Comments to \tht{5.7}}
\roster
\item If $n=1$ then \tht{5.7} turns into
$$
(1+x)^k=\sum_{m\ge 0} \frac1{m!}\,(k\f m)\,x^m \,. \tag5.8
$$
This is why \tht{5.7} is called the binomial formula. 
\item Let us regard the character $gl(n)_\l$ as a function on the whole
group $GL(n)$. Then \tht{5.7} may be written as 
$$
gl_\l(1+X)=\dimn\l \sum_\mu
\frac1{(n\u\mu)}\, s^*_\mu(\l) \, S_\mu(X)\,, \tag5.9
$$
where $X$ is a $n\times n$ matrix. This may be viewed as a Taylor--type
expansion  of the character at $1\in GL(n)$. 
\item The binomial formula \tht{5.7} is not new. In a different form 
it can be found 
in \cite{M1}, Ch. I, section 3, Example 10. 
What is new is the observation that the 
coefficients of the expansion are essentially the $s^*$-functions 
in $\l$. 
\item Note a remarkable symmetry  between
$(\l_1,\ldots,\l_n)$ 
and $(x_1,\ldots,x_n)$ in the right--hand side of \tht{5.7}.
\endroster
\endexample

Now we will give two different proofs of Theorem 5.1: the first proof
follows the arguments sketched in Macdonald \cite{M1}, Ch. I, section 3,
Example 10, while the second
proof is based on the Characterization Theorem.

\demo{First Proof} Put
$$
l_i=\l_i+n-i,\qquad i=1,\ldots,n. \tag5.10
$$
Then by \tht{5.8}
$$
\align
\det\big[(1+x_i)^{\l_j+n-j}]&=
\det\left[\sum_{m_i\ge0}
\frac{(l_j\f m_i)}{m_i!}\,x_i^{m_i}\right]
\,,\quad 1\le i,j\le n,\tag5.11 \\
 &=\sum_{m_1,\ldots,m_n\ge0}
\frac{\det[(l_j\f m_i)]}
{m_1!\ldots m_n!}
x_1^{m_1}\ldots x_n^{m_n}\;. 
\tag5.12
\endalign
$$
Remark that the determinant in \tht{5.12} vanishes unless the numbers
$m_1,\ldots,m_n$ are pairwise distinct; moreover
it is antisymmetric with respect to
permutations of these numbers. Hence 
$$
\det[(1+x_i)^{\l_j+n-j}]=\sum_{m_1>\ldots>m_n\ge0}
\frac{\det[(l_j\f m_i)]}
{m_1!\ldots m_n!}
\det[x_j^{m_i}] \,.
\tag5.13
$$
Put $\mu_i=m_i-n+i$, $i=1,\dots,n$. We obtain from \tht{5.4}, \tht{5.13}, 
and \tht{1.9} 
$$ 
gl(n)_\l (1+x_1,\ldots,1+x_n)
=\sum_{\mu}
\frac{\det[(l_j\f m_i)]}
{m_1!\ldots m_n!}
s_\mu(x_1,\ldots,x_n) 
\tag5.14
$$
In particular, if $x_1=\ldots=x_n=0$ then all summands
vanish except those corresponding to
$\mu=0$, and we obtain Weyl's dimension formula \tht{5.6}
$$
\dimn\l=\frac{\det[(\l_j+n-j\f n-i)]}
{\prod_{i}(n-i)!}
=\frac
{\prod_{i<j}(\l_i-\l_j+j-i)}{\prod_{i}(n-i)!}\,.\tag5.15
$$
It follows from \tht{5.14} and \tht{5.15} that
$$
gl(n)_\lambda(1+x_1,\ldots,1+x_n)
=\dimn\l\sum_\mu
\prod_{i=1}^n\frac{(n-i)!}{(\mu_i+n-i)!}\,
s^*_\mu(\l)\,
s_\mu(x_1,\ldots,x_n)\,.\tag5.16
$$
By virtue of \tht{5.3} this coincides with the desired formula.\qed 
\enddemo  

For the second proof we need some notation which will also be used later.

Let $\C[GL(n)]$ be the algebra of regular functions on the affine complex
algebraic group $GL(n)$. There is a canonical pairing
$$
\la\cdot,\cdot\ra:\quad 
\Ugln\otimes \C[GL(n)]\to\C,\tag5.17
$$
which arises when one looks at elements of $\Ugln$ as distributions
supported at the unity.
\footnote{ Since $GL(n)$ is a complex Lie group and distributions are usually
considered on real manifolds, one could interpret $\Ugln$ as the algebra of
distributions on the real form $U(n)\subset GL(n)$, supported at $1\in
U(n)$. Similarly, the elements of $\C[GL(n)]$ can  be viewed as
functions on  $U(n)$. Finally, as the space of test functions one could 
take, instead of
$\C[GL(n)]$, the larger space of smooth functions on the group $U(n)$.} 

Note that $\C[GL(n)]$ is formed by matrix coefficients of
finite--dimensional $GL(n)$-modules. If $X$ is an element of $GL(n)$, $V$
is a finite--dimensional $GL(n)$-module (also viewed as a $\gln$-module),
$\xi\in V$, $\eta\in V^*$, and $f_{\xi\eta}=((\cdot)\xi,\eta)$ is the
corresponding matrix coefficient, then
$$
\la X,f_{\xi\eta}\ra=(X\xi,\eta),\tag5.18
$$
where $X$ in the right--hand side is viewed as an operator in $V$.

By $I(GL(n))$ we denote the subspace of central functions in $\C[GL(n)]$;
it is spanned by the irreducible characters $gl(n)_\l$, where $\l$ is an
arbitrary signature. Note that there exists a canonical projection
$$
\C[GL(n)]\to I(GL(n)),\tag5.19
$$ 
the unique projection commuting with the action of $GL(n)$ by conjugations.

Recall that by $\C[\gln]$ we denote the algebra of polynomial functions on
$\gln$, and by $I(\gln)$ the subspace of invariant functions. This
subspace is spanned by the polynomials $S_{\mu|n}$. We have
again a canonical projection (introduced in section 2)
$$
\C[\gln]\to I(\gln).\tag5.20
$$

Note that the both projections, \tht{5.19} and \tht{5.20}, can be defined as
averaging with respect to the action of the compact group $U(n)$ by
conjugations.

\demo{Second Proof of Theorem 5.1} Suppose
$\mu$ ranges over the set of partitions with $\ell(\mu)\le n$. There exist
central distributions $\psi_\mu\in \Zgln$ such that for any $F\in I(GL(n))$ 
$$
F(1+X)=\sum_\mu \la \psi_\mu,F\ra S_\mu(X), \qquad X\in\gln.\tag5.21
$$
Indeed, to derive \tht{5.21} we write the Taylor expansion of $F(1+X)$ and
then `average' it using the projection \tht{5.20}. 

Let $\widetilde\psi_\mu \in\Ls(n)$ correspond to
$\psi_\mu$ under the Harish--Chandra isomorphism $\Zgln\cong \Ls(n)$, see
section 2. It follows from \tht{5.21} that $\deg\psi_\mu\le|\mu|$, hence
$$
\deg\widetilde\psi_\mu\le |\mu|.\tag5.22
$$

Let $\l$ be a signature, let $V_\l$ be the corresponding irreducible
$GL(n)$-module (also viewed as a $\gln$-module), and let $gl(n)_\l(X)$ be
the character of $V_\l$. Take $F=gl(n)_\l$ and remark that
$$
\la\psi_\mu,gl(n)_\l\ra=\tr_{V_\l}(\psi_\mu)=
\dimn\l\cdot\widetilde\psi_\mu(\l).\tag5.23
$$
It follows
$$
gl(n)_\l(1+X)=\dimn\l\cdot\sum_\mu 
\widetilde\psi_\mu(\l)S_\mu(X).\tag5.24
$$ 

Now it suffices to prove that
$$
\widetilde\psi_\mu=\frac 1{(n\u\mu)}s^*_{\mu|n}.\tag5.25
$$
To do this suppose $\l$ is a partition ($\ell(\l)\le n$), write $S_\l$
instead of $gl(n)_\l$ and rewrite \tht{5.24} as
$$
S_\l(1+X)=\dimn\l\cdot\sum_\mu \widetilde\psi_\mu(\l)S_\mu(X).\tag5.26
$$ 

On the other hand, observe that 
$$
S_\l(1+X)=S_\l(X)\;+\;\text{terms of degree strictly less than
$|\l|$}.\tag5.27 
$$
By comparing \tht{5.26} and \tht{5.27} we see that
$$
\gather
\widetilde\psi_\mu(\l)=0\qquad \text{if $|\l|\le|\mu|$, $\l\ne\mu$},\\
\dimn\mu\cdot\widetilde\psi_\mu(\mu)=1.
\endgather
$$
By using \tht{5.22} and applying Characterization Theorem I$'$ 
(Theorem 3.3) we conclude that
$$
\widetilde\psi_\mu=\frac{H(\mu)}{\dimn\mu}s^*_{\mu|n}=
\frac 1{(n\u\mu)}s^*_{\mu|n}
$$
(here we have used \tht{5.5}). This proves \tht{5.25} and completes 
the proof. \qed
\enddemo

\head 6. Eigenvalues of higher Capelli operators\endhead

Let $\Mnm$ denote the space of $n\times m$ matrices over $\C$ and let
$\C[\Mnm]$ denote the algebra of polynomial functions on $\Mnm$. The groups
$GL(n)$ and $GL(m)$ act on $\Mnm$ by left and right multiplications,
respectively. Thus, $\C[\Mnm]$ becomes a bi--module over $(GL(n),GL(m))$
and also over $(\Ugln,$ $\Uglm)$ (we may arrange the actions so that this
module would decompose on {\it polynomial\/} irreducible submodules over
$\Ugln$ and $\Uglm$).

Let $D(\Mnm)$ stand for the algebra of (analytic) differential operators 
on $\Mnm$ with polynomial coefficients. The space $\C[\Mnm]$ has a natural
structure of a $D(\Mnm)$-module, and the actions of $\Ugln$ and $\Uglm$
in the same space can be presented as algebra morphisms
$$
\gather
L:\;\Ugln\to\; D(\Mnm),\tag6.1\\
R:\;\Uglm\to\; D(\Mnm);\tag6.2
\endgather
$$
$L$ and $R$ are defined on the generators of both algebras (the
matrix units) as follows:
$$
\align
L(E_{pq})&=\sum_{j=1}^m x_{pj}\p_{qj},\qquad 1\le p,q\le n,\tag6.3\\
R(E_{rs})&=\sum_{i=1}^n x_{ir}\p_{is},\qquad 1\le r,s\le m. \tag6.4
\endalign
$$
(Here and below we denote by $x_{ij}$ the natural coordinate functions on
$\Mnm$ and by $\p_{ij}$ --- the corresponding partial derivatives.)

The aim of this section is to produce a certain differential operator
$\Delta_\mu^{(n,m)}$, indexed by an arbitrary partition $\mu$ with
$\ell(\mu)\le\min(n,m)$, and to prove the following 
identity: 
$$
L(\S_{\mu|n})=R(\S_{\mu|m})=\Delta_\mu^{(n,m)}.\tag6.5
$$

Introduce two formal $n\times m$ matrices $X=[x_{ij}]$ and $D=[\p_{ij}]$
and denote by $D'$ the transposed $m\times n$ matrix. To each partition
$\mu\vdash k$, $\ell(\mu)\le\min(n,m)$, we associate the following
differential operator 
$$
\align
\Dmnm&=\tr(X^{\otimes k}\cdot (D')^{\otimes k} 
\cdot \chm / k!) \\
&=(1/k!)
\sum_{i_1,\dots,i_k=1\,}^n 
\sum_{j_1,\dots,j_k=1\,}^m
\sum_{s\in S(k)} \chm(s)\cdot
x_{i_1 j_1} \dots x_{i_k j_k}
\p_{i_{s(1)} j_1} \dots \p_{i_{s(k)} j_k},\tag6.6
\endalign
$$
where $\chm$ is the irreducible character of the symmetric group $S(k)$,
indexed by $\mu$. We call $\Dmnm$ the {\it higher Capelli operator\/}.

A few comments to this formula: we can write
$$
\align
X&=\sum x_{ij}\otimes E_{ij}\in D(\Mnm)\otimes\Hom(\C^m,\C^n),\tag6.7\\
D&=\sum \p_{ij}\otimes E_{ij}\in D(\Mnm)\otimes\Hom(\C^m,\C^n),\tag6.8\\
D'&=\sum \p_{ij}\otimes E_{ji}\in D(\Mnm)\otimes\Hom(\C^n,\C^m),\tag6.9
\endalign
$$
and
$$
\align
X^{\otimes k}&\in 
D(\Mnm)\otimes\Hom((\C^m)^{\otimes k},(\C^n)^{\otimes k}), \tag6.10\\
(D')^{\otimes k}&\in 
D(\Mnm)\otimes\Hom((\C^n)^{\otimes k},(\C^m)^{\otimes k}), \tag6.11
\endalign
$$
so that
$$
X^{\otimes k}\cdot (D')^{\otimes k}\in D(\Mnm)\otimes 
\Hom((\C^n)^{\otimes k},(\C^n)^{\otimes k}).\tag6.12
$$
Further, the character
$$
\chm=\sum_{s\in S(k)}\chm(s)\cdot s\in\C[S(k)]
$$
is identified with its image in 
$\Hom((\C^n)^{\otimes k},(\C^n)^{\otimes k})$ under the action of $S(k)$ in
the tensor space $(\C^n)^{\otimes k}$. Thus, the product 
$$
X^{\otimes k}\cdot (D')^{\otimes k} \cdot \chm 
$$
is well-defined as an element of 
$$
D(\Mnm)\otimes 
\Hom((\C^n)^{\otimes k},(\C^n)^{\otimes k}).
$$
Finally, after taking the trace with respect to the space
$(\C^n)^{\otimes k}$ we obtain an element of $D(\Mnm)$.

\proclaim{\smc Proposition 6.1}
The operator $\Dmnm$ is $GL(n)\times GL(m)$-invariant.
\endproclaim

\demo{Proof}
Under the action of an element $(g,h)\in GL(n)\times GL(m)$ the matrices $X$
and $D'$ are transformed as follows
$$
X@>(g,h)>>gXh',\qquad D@>(g,h)>>(g')^{-1}Dh^{-1} \,.
$$

Therefore the action of $(g,h)$ maps $\Dmnm$  to
$$
\multline
\tr \big(g^{\otimes n} \cdot X^{\otimes n}
\cdot (h')^{\otimes n} \cdot ((h')^{-1})^{\otimes n}
 \cdot (D')^{\otimes n} 
\cdot (g^{-1})^{\otimes n})
\cdot \chm/k! \big) \\ 
= 
\tr \big(X^{\otimes n} \cdot (D')^{\otimes n}
\cdot \chm/k! \big) \,.   \qed
\endmultline
$$
\enddemo

\example{\smc Example 6.2}
Suppose $\mu=(1^k)$, where $k=1,\ldots,\min(n,m)$, and let the indices
$i_1,\ldots,i_k$ range over $\{1,\ldots,n\}$ while the the indices
$j_1,\ldots,j_k$ range over $\{1,\ldots,m\}$. Then $\Dmnm$ turns into the
$k$-th {\it classical Capelli operator\/} (see Howe--Umeda \cite{HU},
Nazarov \cite{N1})
$$
\align
\Delta_{(1^k)}^{(n.m)}
&=(1/k!)
\sum_{i_1,\dots,i_k\,} 
\sum_{j_1,\dots,j_k\,}
\sum_{s\in S(k)} \sgn(s)\cdot
x_{i_1 j_1} \dots x_{i_k j_k}
\p_{i_{s(1)} j_1} \dots \p_{i_{s(k)} j_k} \\
&=(1/k!)^2
\sum_{i_1,\dots,i_k\,} 
\sum_{j_1,\dots,j_k\,}
\sum_{s,t\in S(k)} \sgn(s t)\cdot
x_{i_1 j_1} \dots x_{i_k j_k}
\p_{i_{s t(1)} j_1} \dots \p_{i_{s t(k)} j_k} \\
&=(1/k!)^2
\sum_{i_1,\dots,i_k\,} 
\sum_{j_1,\dots,j_k\,}
\sum_{s,t\in S(k)} \sgn(s)\sgn(t)\cdot
x_{i_{t(1)} j_1} \dots x_{i_{t(k)} j_k}
\p_{i_{s(1)} j_1} \dots \p_{i_{s(k)} j_k} \\
&=(1/k!)^2
\sum_{i_1,\dots,i_k\,} 
\sum_{j_1,\dots,j_k\,}
\det 
\left[
\matrix
x_{i_1 j_1}&\hdots&x_{i_1 j_k}\\
\vdots&&\vdots\\
x_{i_k j_1}&\hdots&x_{i_k j_k}
\endmatrix \right]
\det 
\left[
\matrix
\p_{i_1 j_1}&\hdots&\p_{i_1 j_k}\\
\vdots&&\vdots\\
\p_{i_k j_1}&\hdots&\p_{i_k j_k}
\endmatrix \right] \\
&=
\sum_{i_1<\dots<i_k\,} 
\sum_{j_1<\dots<j_k\,}
\det 
\left[
\matrix
x_{i_1 j_1}&\hdots&x_{i_1 j_k}\\
\vdots&&\vdots\\
x_{i_k j_1}&\hdots&x_{i_k j_k}
\endmatrix \right]
\det 
\left[
\matrix
\p_{i_1 j_1}&\hdots&\p_{i_1 j_k}\\
\vdots&&\vdots\\
\p_{i_k j_1}&\hdots&\p_{i_k j_k}
\endmatrix \right]
\endalign
$$
This is why the operators $\Dmnm$ are called
the {\it higher Capelli operators}.
\endexample

\example{\smc Example 6.3} Suppose $\mu=(k)$, where $k=1,2,\ldots$. The
operator $\Delta_{(k)}^{(n,m)}$ can be written as follows
$$
\gather
\Delta_{(k)}^{(n,m)}=
\sum_{1\le i_1\le\dots\le i_k\le n\,} 
\sum_{1\le j_1\le\dots\le j_k\le m\,}
\frac 1{p_1!p_2!\ldots q_1!q_2!\ldots}\\
\cdot\per \left[
\matrix
x_{i_1 j_1}&\hdots&x_{i_1 j_k}\\
\vdots&&\vdots\\
x_{i_k j_1}&\hdots&x_{i_k j_k}
\endmatrix \right]
\per 
\left[
\matrix
\p_{i_1 j_1}&\hdots&\p_{i_1 j_k}\\
\vdots&&\vdots\\
\p_{i_k j_1}&\hdots&\p_{i_k j_k}
\endmatrix \right],\tag6.13
\endgather
$$
where $p_1,p_2,\ldots$ and $q_1,\ldots,q_2,\ldots$ stand for the
multiplicities in the multisets $\{i_1,\ldots,i_k\}$ and
$\{j_1,\ldots,j_k\}$, respectively, and `$\per$' means `permanent'.
The differential operators \tht{6.13} first appeared in Nazarov's work
\cite{N1}. 
\endexample

>From the very definition of the operators $\Dmnm$ it
readily follows that their dependence on $n$ and $m$ is not
essential. Suppose $N\ge n$ and $M\ge m$. We have the natural
projection
$$
\M(N,M)\to\Mnm \tag6.14
$$
(excision of the upper left corner of size $n\times m$ from a $N\times M$
matrix) and the corresponding inclusion (lifting)
$$
\C[\Mnm]\to \C[\M(N,M)] \,. \tag6.15
$$

\proclaim{\smc Proposition 6.4}
The subspace $\C[\Mnm]$ is an invariant subspace
of the operator $\DmNM$ and
$$
\DmNM|_{\C[\Mnm]}
=\Dmnm.\tag6.16
$$
\endproclaim 

\demo{Proof} This follows at once from the fact that for each summand
occurring in the definition of the operator $\Dmnm$ (see \tht{6.6}) the
coefficient and the derivatives depend on the same set of indices. \qed
\enddemo

In other words, the expression \tht{6.6} is {\it stable\/} as
$n,m\to\infty$.

It is well known that the space $\C[\Mnm]$  as
$GL(n)\times GL(m)$ bi--module decomposes into a multiplicity free sum,
$$
\C[\Mnm]=
\bigoplus_{\l} V_{\l|n}\otimes V_{\l|m}\,,
\quad  \ell(\l)\le \min(n,m) \,,\tag6.17
$$
and the highest vectors are
$$
v_\l=\prod_i 
\det
\left[
\matrix
x_{11}&\hdots&x_{1i}\\
\vdots&&\vdots\\
x_{i1}&\hdots&x_{ii}
\endmatrix
\right]^{\l_i-\l_{i+1}}\,. \tag6.18
$$
By virtue of Proposition 6.1 each higher Capelli operator $\Dmnm$
acts in the subspace $V_{\l|n}\otimes V_{\l|m}$ as a scalar
operator, say $f_\mu(\l)$. We can obtain the scalar $f_\mu(\l)$ by
applying $\Dmnm$ to the highest vector $v_\l$. Now remark that the polynomial
function $v_\l$ does not depend on $n$ and $m$. It follows that the scalar
$f_\mu(\l)$ does not depend on $n,m$, too. 

\proclaim{\smc Proposition 6.5} For each $\mu$ the eigenvalue $f_\mu(\l)$
of the higher Capelli operator $\Dmnm$ in $V_{\l|n}\otimes
V_{\l|m}$ is a shifted symmetric polynomial in $\l$ of degree $\le|\mu|$.
\endproclaim

\demo{Proof} By a well--known result (see, e.g., Howe \cite{H}), 
each of the two subalgebras
$$
L(\Ugln)\subset D(\Mnm), \qquad R(\Uglm)\subset D(\Mnm)
$$ 
coincides with the centralizer of the other subalgebra in $D(\Mnm)$. 
Hence $\Dmnm$ belongs
both to the center of $L(\Ugln)$ and $R(\Uglm)$. Take the minimal of the
numbers $n,m$; assume this is $n$. Then the morphism $L$ is an embedding,
so that $\Dmnm$ is the image of a central element
of $\Ugln$. 
\footnote {In fact the latter claim holds even if $n>m$, see Remark 6.10
below.} 
But we know 
that the eigenvalue of each central element in $V_{\l|n}$ is a shifted
symmetric polynomial in $\l$ (see section 2). 

Further, we have $f_\mu(\l)=(\Dmnm v_\l)(1)$. Since $\Dmnm$ is of order
$|\mu|$, it follows from the form of $v_\l$ that the polynomial
$f_\mu(\l)$ has degree $\le|\mu|$. This completes the proof.
\qed 
\enddemo

Propositions 6.4 and 6.5 imply the following
\proclaim{\smc Corollary 6.6} The eigenvalue $f_\mu(\l)$ is a shifted
symmetric function in $\l$. \qed
\endproclaim

\proclaim{\smc Theorem 6.7} For each partition $\mu$ the eigenvalue
$f_\mu(\l)$ of the higher Capelli operator $\Dmnm$ in the irreducible
component $V_{\l|n}\otimes V_{\l|m}$ of the decomposition \tht{6.17}
coincides with $s_\mu^*(\l)$.
\endproclaim

(It is tacitly assumed that $n,m\ge\ell(\mu)$.)

\demo{Proof} We will apply Characterization Theorem I (Theorem 3.2). 
Observe that the component $V_{\l|n}\otimes V_{\l|m}$
is contained in the subspace of polynomials of degree $|\l|$.
{F}rom \tht{6.6} it is clear that $\Dmnm$ annihilates all 
polynomials of degree $<|\mu|$. Therefore
$$
f_\mu(\l)=0, \qquad |\l|<|\mu| \,.
$$
Suppose $|\l|=|\mu|=k$. By virtue of Proposition 6.4 we may assume $m=k$. 
Let $\{e_i\}$, $i=1,\dots,n$, be the standard basis of $\C^n$. 
Embed $(\C^n)^{\otimes k}$ in $\C[\M(n,k)]$ as follows: 
$$
e_{i_1}\otimes\dots\otimes e_{i_k} \to
x_{i_1 1}\dots x_{i_k k} \,. \tag6.19
$$
This embedding is $GL(n)$-equivariant. Hence it takes the isotypic
component of $V_{\l|n}$ in the tensor space $(\C^n)^{\otimes k}$ to
$V_{\l|n}\otimes V_{\l|k}$. 
 
On the other hand, it follows from \tht{6.6} (and the fact that
$\chm(s^{-1})=\chm(s)$) that the restriction of the operator $\Dm^{(n,k)}$ to the
image of the tensor space is reduced to the following action:
$$
e_{i_1}\otimes\dots\otimes e_{i_k} \mapsto
\sum_{s\in S(k)} \chm(s) 
e_{i_{s(1)}}\otimes\dots\otimes e_{i_{s(k)}} \,.
$$
I.e., the action of $\Dm^{(n,k)}$ turns into the action of 
the central element $\chm\in\C[S(k)]$ in $(\C^n)^{\otimes k}$. By the
Schur--Weyl duality, the eigenvalue of $\chm$ in the isotypic component of
$V_{\l|n}$ is the same as in the irreducible $S(k)$-module, indexed by
$\l$. Thus, this eigenvalue (which coincides with $f_\mu(\l)$) is equal to
$$
\dt_{\l\mu} k!\ov \dim\mu = \dt_{\l\mu} H(\mu) \,.
$$
By Characterization Theorem I we conclude that $f_\mu(\l)=s_\mu^*(\l)$.
\qed 
\enddemo  

\proclaim{\smc Corollary 6.8} The identity \tht{6.5} holds \qed.
\endproclaim

Note that in the proof of Theorem 6.7, in place of Characterization
Theorem I one could use 
Characterization Theorem II (Theorem 3.4). Then one has to check that the
highest term of $f_\mu$ equals $s_\mu$. This is equivalent to the fact that 
$\Dmnm$ and $L(\S_{\mu|n})$ coincide up to lower terms, which can be
deduced from \tht{2.10} and \tht{6.6}, by making use of \tht{6.3}.

The identity \tht{6.5} plus an explicit expression for the quantum immanants
(which is obtained in \cite{Ok1}, \cite{N2}, and \cite{Ok2}) can be called 
the {\it higher analogue of the classical Capelli 
identity\/}. (The well--known classical Capelli identity corresponds
to the case $\mu=(1^k)$.)

\example{\smc Remark 6.9} Let 
$$
\M(\infty,\infty)=\varprojlim\Mnm
$$
be the space of all $\infty\times\infty$ matrices and let
$$
\C[\M(\infty,\infty)]=\varinjlim \C[\Mnm]
$$
be the space of (cylindrical) polynomial functions on $\M(\infty,\infty)$.
If $\mu$ is fixed and $n,m\to\infty$ then the collection
$\Delta_\mu=\{\Dmnm\}$ can be viewed as an `infinite--dimensional
differential operator', acting in $\C[\M(\infty,\infty)]$. It was shown in
Olshanski's paper \cite{O2} that the algebra $\Ls$ can be realized as
the algebra of all differential operators on $\M(\infty,\infty)$ commuting
with the action of the group 
$$
GL(\infty)\times GL(\infty)=\varinjlim GL(n)\times GL(m).
$$
In this realization, the operators $\Delta_\mu$ just correspond to the
shifted Schur functions $s^*_\mu\in\Ls$.
\endexample

\example{\smc Remark 6.10} Let us return to the proof of Proposition 6.5
and note that the center of the algebra $L(\Ugln)$ always coincides 
with $L(\Zgln)$, even if $L$ is not an embedding. 

To see this we use the canonical projection 
$$
\Ugln\to\Zgln, \tag6.20
$$
the unique projection commuting with the adjoint action of $\gln$; it
exists because $\Ugln$ is a semisimple $\gln$-module (cf. \cite{Bou}, Ch.
VIII, 8.5, Proposition 7). Next, consider 
$D(\Mnm)$ as a
$\gln$-module with respect to commutation with elements of $L(\gln)$; this
module also is semisimple, whence there exists a unique
equivariant projection of $D(\Mnm)$ onto the subspace of $\gln$-invariants
(i.e., onto differential operators commuting with $L(\Ugln)$). The map
$L:\Ugln\to D(\Mnm)$ is clearly compatible with these two projections. 

Now take 
an arbitrary element $A$ in the center of $L(\Ugln)$ and choose
$B\in\Ugln$ such that $L(B)=A$. Applying to $B$ the first projection we
obtain a central element $B'$, and $L(B')$ coincides with the image $A'$
of $A$ under the second projection. But we have $A'=A$, whence $L(B')=A$,
so that $A$ lies in the image of the center.

Finally, note that the both projections can also be defined as averaging over
the compact form $U(n)\subset GL(n)$.
\endexample

\head 7. Capelli--type identity for Schur--Weyl duality \endhead

Recall that the Schur--Weyl duality between the general linear and
symmetric groups is based on the following decomposition of the tensor
space $(\C^n)^{\otimes l}$ as a bi--module over $(GL(n),S(l))$
($n,l=1,2,\ldots$) 
$$
(\C^n)^{\otimes l}=\sum_\l (V_{\l|n}\otimes W_\l),
\qquad \l\vdash l,\quad\ell(\l)\le n.\tag7.1
$$
Here $W_\l$ denotes the irreducible $S(l)$-module, indexed by $\l$, and
$V_{\l|n}$, as above, is the irreducible polynomial $GL(n)$-module (or
$\Ugln$-module), also indexed by $\l$.

Let 
$$
\tau_{GL(n)}: \Ugln\to \End((\C^n)^{\otimes l}) \tag7.2
$$
and
$$
\tau_{S(l)}: \C[S(l)]\to \End((\C^n)^{\otimes l}) \tag7.3
$$
denote the algebra morphisms defined by the respective actions.  The
decomposition \tht{7.1} implies that each of the subalgebras
$$
\tau_{GL(n)}(\Ugln)\subset \End((\C^n)^{\otimes l}),\qquad
\tau_{S(l)}(\C[S(l)])\subset \End((\C^n)^{\otimes l})\tag7.4
$$
is the centralizer of the other. It follows that the centers of the both
subalgebras \tht{7.4} are the same. But these two centers coincide with the
images of the centers of the algebras $\Ugln$ and $\C[S(l)]$, respectively
(this can be shown as in Remark 6.10 above), whence we obtain that
$$
\tau_{GL(n)}(\Zgln)=\tau_{S(l)}(Z(S(l))),\tag7.5
$$
where $Z(S(l))$ stands for the center of the group algebra $\C[S(l)]$.

A natural problem arising from \tht{7.5} is to produce `sufficiently many'
couples $(A,a)$ of central elements $A\in\Zgln$, $a\in Z(S(l))$
satisfying the identity 
$$
\tau_{GL(n)}(A)=\tau_{S(l)}(a).\tag7.6
$$
By analogy with the duality between $GL(n)$ and $GL(m)$ acting in
$\C[\Mnm]$ we call such identities the {\it Capelli--type identities for the
Schur--Weyl duality\/}.

A family of identities of type \tht{7.6} was found in \cite{KO} (see also
section 15 below). Now we present another family.

\proclaim{\smc Theorem 7.1} Let $\mu\vdash k$ be a partition such that
$k\le l$ and $\ell(\mu)\le n$, and let
$$
\Ind\chm =\sum_{t\in S(l)/S(k)} t\cdot\chm\cdot t^{-1}\in Z(S(l)) \tag7.7
$$
be the induced character. 
Then
$$
\tau_{GL(n)}(\S_{\mu|n})=\tau_{S(l)}(\Ind\chm/(l-k)!).\tag7.8
$$
\endproclaim

\demo{Proof} The idea is to reduce the problem to the duality between
two general linear groups, studied in section 6.

As in the proof of Theorem 6.7, let us embed the tensor space
$(\C^n)^{\otimes l}$ into the space of polynomials 
$\C[\M(n,l)]$: 
$$
\Phi: e_{r_1}\otimes\dots\otimes e_{r_l} \to
x_{r_1 1}\dots x_{r_l l} \,. \tag7.9
$$
The map $\Phi$ commutes with the action of $\Ugln$. By Corollary 6.8,
$\Phi$ intertwines the operators $\tau_{GL(n)}(\S_{\mu|n})$ and
$\Delta_\mu^{(n,l)}$. Hence it suffices to prove that $\Phi$ also
intertwines $\tau_{S(l)}$ and $\Delta_\mu^{(n,l)}$, and this can be 
done by a direct computation.

Indeed, write
$$
\Delta_\mu^{(n,l)}=(1/k!)
\sum_{i_1,\dots,i_k=1\,}^n 
\sum_{j_1,\dots,j_k=1\,}^l
\sum_{s\in S(k)} \chm(s)\cdot
x_{i_{s(1)} j_1} \dots x_{i_{s(k)} j_k}
\p_{i_1 j_1} \dots \p_{i_k j_k}
$$
and apply this operator to the monomial
$$
x_{r_1 1}\ldots x_{r_l l}. \tag7.10
$$
Then the result can be described as follows. We choose an arbitrary subset
$J=\{j(1)<\ldots<j(k)\}\subset\{1,\ldots,l\}$, next replace 
in \tht{7.10} each
of the letters $x_{r_{j(1)}j(1)}$, $\ldots$, $x_{r_{j(k)}j(k)}$ by 
$x_{r_{j(s(1))}j(1)}$, $\ldots$, $x_{r_{j(s(k))}j(k)}$, respectively, 
where $s$ is an arbitrary 
permutation of $1,\ldots,k$, then multiply by $\chm(s)$ and sum over all $s$
and all $J$.

On the other hand, the same result is obtained if we first apply to the
tensor $e_{r_1}\otimes\dots\otimes e_{r_l}$ the operator
$\tau_{S(l)}(\Ind\chm/(l-k)!)$ and then apply the map $\Phi$. This
completes the proof. \qed
\enddemo

\head 8. Dimension of skew Young diagrams \endhead

Let $\mu\vdash k$ and $\l\vdash l$ be two Young diagrams such that $k\le
l$ and $\mu\subset\l$, and let $\l/\mu$ be the corresponding skew diagram.
The number $\dim\l/\mu$, the dimension of $\l/\mu$, equals the number of
standard tableaux of shape $\l/\mu$ or, that is the same, the number of
paths from $\mu$ to $\l$ in the Young graph (see Vershik--Kerov
\cite{VK1}, \cite{VK3}). Recall that the
Young graph is the oriented graph whose vertices are Young diagrams and
two diagrams $\mu$ and $\nu$ are connected by an edge if $\nu$ is obtained
from $\mu$ by adding a single box. 

Equivalently, in terms of the irreducible characters $\chi^\l$ and $\chm$,
indexed by $\l$ and $\mu$,
$$
\dim \l/\mu = \la\Res\chi^\l,\chm\ra_{S(k)}\,, \tag8.1
$$
where $\Res\chi^\l$ is the restriction of $\chi^\l$ to
$S(k)$ and $\la\cdot,\cdot\ra_{S(k)}$ is the standard scalar
product of functions on $S(k)$
$$
\la\phi,\psi\ra_{S(k)}=
(1/k!)\sum_{s\in S(k)} \phi(s) \overline{\psi(s)} \,. \tag8.2
$$

\proclaim{\smc Theorem 8.1} We have
$$
\frac{\dim \l/\mu}{\dim \l} =
\frac{s^*_\mu(\l)}{(l\f k)} \,. \tag8.3
$$
\endproclaim 

Note that this is an explicit formula for $\dim\l/\mu$, because for
$\dim\l$ there are nice formulas (see \tht{3.7}, \tht{3.8}). 
We will prove \tht{8.3} in
three different ways. The first proof is an immediate application of the
Capelli--type identity \tht{7.8}. The second proof is a direct calculation which
uses only formula \tht{8.1} and the definition \tht{1.6} of shifted Schur
polynomials. \footnote{In fact this was the first proof we found, and it
was inspired by Macdonald's letter \cite{M3}, see a discussion in section
15.} The third proof is given in section 9 below as a corollary of
a Pieri--type formula for the $s^*$-functions.

\demo{First Proof} Choose an arbitrary natural $n\ge\ell(\l)$ and let us 
calculate the eigenvalues of the both sides of
the identity \tht{7.8} in the irreducible component $V_{\l|n}\otimes W_\l$ of
the decomposition \tht{7.1}. The eigenvalue of $\tau_{GL(n)}(\S_{\mu|n})$
equals the eigenvalue of $\S_{\mu|n}$ in $V_{\l|n}$, which equals 
$s^*_\mu(\l)$ by the very definition of $\S_{\mu|n}$. 

The eigenvalue of $\tau_{S(l)}(\Ind\chm)$  is equal to
$$
\gather
\dfrac 1{\dim\l}\tr_{W_\l}(\Ind\chm)=
\dfrac{l!}{\dim\l}\la\chi^\l,\Ind\chm\ra_{S(l)}\\
=\dfrac{l!}{\dim\l}\la\Res\chi^\l,\chm\ra_{S(k)}=
\dfrac{l!\dim\l/\mu}{\dim\l},\tag8.4
\endgather
$$
where the two last equalities are given by Frobenius reciprocity and
formula \tht{8.1}, respectively.

By virtue of the identity \tht{7.8}, we have
$$
s^*_\mu(\l)=\dfrac{\dim\l/\mu}{\dim\l}\;\dfrac{l!}{(l-k)!},
$$
which is equivalent to \tht{8.3}. \qed
\enddemo

\demo{Second Proof} By using formula \tht{8.1}, Frobenius reciprocity and the
characteristic map (see Macdonald's book \cite{M1}, Ch. I, section 7) 
we express $\dim\l/\mu$
in terms of Schur functions,
$$
\dim\l/\mu=\la\Res\chi^\l,\chm\ra_{S(k)}=
\la\chi^\l,\Ind\chm\ra_{S(l)}=
(s_\l,s_\mu p_1^{l-k}),\tag8.5
$$
where $(\cdot,\cdot)$ denotes the canonical inner product in the algebra
$\L$.

Further, we fix $n\ge\ell(\l)$ and remark that by a standard argument (see
e.g., \cite{M1}, Ch. I, section 7, Example 7), \tht{8.5} equals the
coefficient of 
$$
x_1^{\l_1+n-1}\ldots x_n^{\l_n}
$$
in the expansion of 
$$
\det[x_i^{\mu_j+n-j}]_{1\le i,j\le n}\;\cdot\; 
(x_1+\ldots+x_n)^{l-k}.\tag8.6
$$

By expanding the both factors of \tht{8.6} we find that \tht{8.6} 
can be written
as
$$
\sum_{s\in S(n)} \sgn(s)\;
\sum_{r_1+\ldots+r_n=l-k}\; \dfrac{(l-k)!}{r_1!\ldots r_n!}
\prod_{i=1}^n x_i^{\mu_{s(i)}+n-s(i)+r_i},
$$
so that the the desired coefficient is equal to
$$
\align
\sum_{s\in S(n)}&\sgn(s)\;\dfrac{(l-k)!}
{\prod_{i=1}^n (\l_i-\mu_{s(i)}-i+s(i))!}\\
&=(l-k)!\det\left[\dfrac 1{(\l_i-\mu_j-i+j)!}\right]\\
&=(l-k)!\det\left[\dfrac 1{((\l_i+n-i)-(\mu_j+n-j))!}\right]\\
&=\dfrac{(l-k)!}{\prod(\l_i+n-i)!}
\det\left[\dfrac{(\l_i+n-i)!}{((\l_i+n-i)-(\mu_j+n-j))!}\right]\\
&=\dfrac{(l-k)!}{l!}\cdot
\dfrac{l!\prod_{p<q}(\l_p-\l_q+q-p)}{\prod(\l_i+n-i)!}
\cdot\dfrac{\det[(\l_i+n-i\f\mu_j+n-j)]}
{\prod_{p<q}(\l_p-\l_q+q-p)}\\
&=\dfrac{\dim\l}{(l\f k)}s^*_\mu(\l),
\endalign
$$
where we have used the definition \tht{1.6} of the shifted Schur polynomials
and formula \tht{3.8} for $\dim\l$. This completes the proof. \qed
\enddemo

\head 9. Pieri--type formula for $s^*$-functions \endhead

In this section we present one more proof of Theorem 8.1. It is based on
a Pieri--type formula \tht{9.2} which is of independent interest.

For two Young diagrams $\mu$ and $\nu$, let the symbol $\mu\nearrow\nu$ mean
that $|\nu|=|\mu|+1$ and $\mu\subset\nu$ (then $\mu$ and $\nu$ are
connected by an edge in the Young graph). 

Recall that 
$$
s^*_{(1)}(x)=s_{(1)}(x)=\sum_i x_i \,.\tag9.1
$$

\proclaim{\smc Theorem 9.1} For any Young diagram $\mu$
$$
s^*_\mu (p_1-|\mu|) = \sum_{\nu,\;\mu\nearrow\nu} s^*_{\nu} \,. \tag9.2 
$$
\endproclaim

\demo{Proof} We shall apply Characterization Theorem II (Theorem 3.4).
Let $k=|\mu|$. The highest terms of the both sides of \tht{9.2} (which are 
of degree $k+1$) agree by the classical Pieri formula. So it suffices to
verify that the left--hand side of \tht{9.2} vanishes at each diagram
$\eta$ with $|\eta|\le k$. We have
$$
s^*_\mu(\eta)(p_1(\eta)-|\mu|)=0,\quad\text{if}\quad 
\eta\ne\mu \,,
$$
because $s^*_\mu(\eta)=0$ then (Theorem 3.1). Further, the same is also
true for $\eta=\mu$, because then
$$
s^*_{(1)}(\eta)-|\mu|=s^*_{(1)}(\mu)-|\mu|=0,
$$
by \tht{9.1}. This completes the proof. \qed
\enddemo

Now we shall deduce Theorem 8.1 from Theorem 9.1. 
\demo{ Third Proof of Theorem {\rm 8.1}}
Recall the equality to be proved: if $\mu\vdash k$, $\l\vdash l$, and
$\mu\subset\l$ then
$$
\frac{\dim \l/\mu}{\dim \l} =
\frac{s^*_\mu(\l)}{(l\f k)} \,. \tag9.3
$$  
Using the identity \tht{9.2} $l-k$ times we get
$$
s^*_\mu(p_1-k)\dots(p_1-l+1)=
\sum_{\nu\supset\mu,\;|\nu|=l} \dim\nu/\mu\cdot s^*_{\nu}\,. \tag9.4 
$$

Let us evaluate the both sides of \tht{9.4} at $\l$. We obtain
$$
s^*_\mu(\l) (l-k)! = \dim\l/\mu\cdot s^*_{\l}(\l)\,. \tag9.5 
$$
In particular for $\mu=\varnothing$ we obtain
$$
l! = \dim\l\cdot s^*_{\l}(\l)\,.
$$
(This is a different proof of \tht{3.3}.) By substituting 
this into \tht{9.5} we
arrive to \tht{9.3}. \qed
\enddemo

\head 10. Coherence property of quantum immanants and shifted Schur
polynomials \endhead

Recall (Proposition 1.3) that the shifted Schur polynomials (just as the
ordinary ones) possess the stability property
$$
s^*_{\mu|n+1}(x_1,\ldots,x_n,0)=s^*_{\mu|n}(x_1,\ldots,x_n),\tag10.1
$$
which can also be formulated in terms of bi--invariant differential
operators on matrix spaces (Remark 6.9). The aim of this section is
to establish another stability property of the polynomials $s^*_{\mu|n}$,
which we call the coherence property. This new property is best stated in
terms of the quantum immanants $\S_{\mu|n}$.

Let $n=1,2,\ldots$. Recall (Remark 6.10) that since the adjoint
representation of the Lie algebra 
$\gln$ in its enveloping algebra $\Ugln$ is completely reducible, there
exists a unique $\gln$-equivariant projection
$$
\Ugln\to\Zgln \tag10.2
$$
on the center, which can also be defined as averaging over the compact
form $U(n)\subset GL(n)$ with respect to its adjoint action,
$$
X\mapsto \int_{U(n)}\operatorname{Ad}(u)\cdot X\;du,
\qquad X\in\Ugln,\tag10.3
$$
where $du$ is the normalized Haar measure.

For each couple $n<N$ define the {\it averaging operator\/} $\Av$,
$$
\Av: \Zgln\to\ZglN,\tag10.4
$$
as the composition of three maps,
$$
\Zgln\to\Ugln\to\UglN\to\ZglN,\tag10.5
$$
where the first and second arrows are natural embeddings and the third
arrow is the projection \tht{10.2} for $\glN$. 

For each $n$ we consider the canonical pairing \tht{5.17} between the
enveloping algebra $\Ugln$ and the space $\C[GL(n)]$ of regular functions
on the group $GL(n)$. For $N>n$ the embedding $GL(n)\to GL(N)$ induces 
the restriction map
$$
\Res_{Nn}: \C[GL(N)]\to \C[GL(n)], \tag10.6
$$
which is dual to the embedding 
$$
\Ugln\to \U(\gl(N)). \tag10.7
$$

Similarly, for each $n$ the embedding
$$
I(GL(n))\to \C[GL(n)] \tag10.8
$$
is dual to the projection \tht{10.2} (recall that by $I(GL(n))$ we denoted the
subspace of central functions in $\C[GL(n)]$).

Therefore, given $Z\in\Zgln$, its image $\Av Z$ under \tht{10.4} can be
characterized as the unique element of $\ZglN$ satisfying
$$
\la\Av Z ,F \ra = 
\la Z, \Res_{Nn} F \ra,\qquad F\in I(GL(N)).  \tag10.9
$$

\proclaim{\smc Theorem 10.1 (Coherence Theorem)} Let $\mu$ be a partition
and let $N>n\ge\ell(\mu)$. Then
$$
\Av \S_{\mu|n} = \frac{\,(n\u \mu)}{(N\u \mu)} \,\, \S_{\mu|N}\;.\tag10.10
$$
\endproclaim

(Recall that by $\S_{\mu|n}$ we denote the quantum immanant, see
Definition 2.3); the symbol $(n\u\mu)$ is defined in \tht{5.2}.)

We call the relation \tht{10.10} the {\it coherence property} of the quantum
immanants (or equivalently, of the shifted Schur polynomials). Three
different proofs of this property will be given. The first and second ones
are obtained by making use of the Binomial Theorem and the
Characterization Theorem, respectively. The third proof is a direct
computation which uses only the initial definition \tht{1.6} of the shifted
Schur polynomials; it leads to an interesting identity involving these
polynomials.  

\demo{First Proof} For arbitrary $n$ and $\mu$ ($\ell(\mu)\le n$) set
$$
\Smnt=\frac 1{(n\u\mu)}\Smn.\tag10.11
$$
We have to prove that
$$
\Av\Smnt=\SmNt, \qquad N>n.\tag10.12
$$
By duality, this can be restated as the relation
$$
\la\SmNt,F\ra=\la\Smnt,\Res_{Nn}F\ra\qquad
\text{for each $F\in I(GL(N))$.}\tag10.13
$$ 

By virtue of the Binomial Theorem (Theorem 5.1) 
and the argument used in its second
proof (see, in particular, formula  \tht{5.21}), we have the expansion
$$
F(1+X)=\sum_{\mu,\;\ell(\mu)\le n}\; 
\la\Smnt,F\ra S_{\mu|n}(X),\qquad X\in\gln,\tag10.14
$$
where the invariant polynomials $S_{\mu|n}(X)$ were introduced in section 2. 

Let us show that the desired relation \tht{10.10} is a formal consequence of
the relation \tht{10.14} and the stability property:
$$
X\in\gl(n)\subset\glN\;\Rightarrow\; 
S_{\mu|N}(X)=\cases S_{\mu|n}(X), &\ell(\mu)\le n,\\
0, &\ell(\mu)>n. \endcases \tag10.15
$$
Indeed, let us write the $N$-th relation of the form \tht{10.14},
$$
F(1+X)=\sum_{\mu,\;\ell(\mu)\le N}\; 
\la\SmNt,F\ra S_{\mu|N}(X),\qquad X\in\glN,\quad F\in I(GL(N)),\tag10.16
$$
and then take $X\in\gl(n)$. By virtue of \tht{10.15} we obtain then
$$
\Res_{Nn}F(1+X)=\sum_{\mu,\;\ell(\mu)\le n}\; 
\la\SmNt,F\ra S_{\mu|n}(X),\qquad X\in\gln.\tag10.17
$$
On the other hand, the $n$-th relation \tht{10.14}, applied to the function
$\Res_{Nn}F$, gives
$$
\Res_{Nn}F(1+X)=\sum_{\mu,\;\ell(\mu)\le n}\; 
\la\Smnt,\Res_{Nn}F\ra S_{\mu|n}(X),\qquad X\in\gln.\tag10.18
$$
Since the polynomials $S_{\mu|n}$, $\ell(\mu)\le n$, are linearly
independent, the comparison of \tht{10.17} and \tht{10.18} 
implies \tht{10.13}. This
completes the proof. \qed
\enddemo

Below we shall need the well--known 
{\it branching rule\/} for the general linear
groups. Let $\l$ be a signature for $GL(n+1)$ and let $V_{\l|n+1}$ be the
corresponding irreducible $GL(n+1)$-module; then its decomposition under
the action of the subgroup $GL(n)\subset GL(n+1)$ looks as follows
$$
V_{\l|n+1}\bigm|_{GL(n)}\;\sim\;\bigoplus_{\nu\prec\l}V_{\nu|n}\;,\tag10.19
$$
where $\nu\prec\l$ means
$$
\l_1\ge\nu_1\ge\l_2\ge\ldots\l_n\ge\nu_n\ge\l_{n+1}\tag10.20
$$
(see, e.g., {\v Z}elobenko \cite{Z}). In terms of characters the branching rule
looks as follows
$$
gl(n+1)_\l\bigm|_{GL(n)}=\sum_{\nu\prec\l}gl(n)_\nu.\tag10.21
$$

\demo{Second Proof} We apply Characterization Theorem I (Theorem 3.2).
Let $k=|\mu|$. The central element $\SmNt$ has degree $k$. Since the
averaging operator does not raise degree, the central element $\Av\Smnt$
has degree $\le k$. Hence, by Characterization Theorem, it suffices to
prove that the both central elements have the same eigenvalues in any
irreducible polynomial module $V_{\l|N}$ such that $|\l|\le k$.

The eigenvalue of $\SmNt$ is equal to $s^*_\mu(\l)/(N\u\l)$. For $\l\ne\mu$
this is zero, and for $\l=\mu$ this number equals
$$
\frac{s^*_\mu(\mu)}{(N\u\mu)}=\frac{H(\mu)}{(N\u\mu)}=
\frac 1{\dim_{GL(N)}\mu}\;.\tag10.22
$$

The eigenvalue of $\Av\Smnt$ can be written as
$$
\dfrac{\la\Av\Smnt,gl(N)_\l\ra}{\dim_{GL(N)}\l}=
\dfrac{\la\Smn,\Res_{Nn} gl(N)_\l\ra}{(n\u\mu)\dim_{GL(N)}\l}\;.\tag10.23
$$
By the branching rule \tht{10.19}
$$
\Res_{Nn}gl(N)_\l=\cases
gl(n)_{\lambda}+\text{lower terms}, &\text{if $\ell(\lambda)\le n$,}\\
\text{lower terms}, &\text{if $\ell(\lambda)>n$,}
\endcases\tag10.24
$$
where `lower terms' means `a linear combination of characters $gl(n)_\nu$
with $|\nu|<k$'.

Recall that
$$
\la\Smn,gl(n)_\nu\ra=\dim_{GL(n)}\nu\cdot s^*_\mu(\nu).
$$
Hence \tht{10.24} implies that the expression \tht{10.23 }
vanishes when $\l\ne\mu$.
If $\l=\mu$ then \tht{10.23} equals
$$
\dfrac{s^*_\mu(\mu)\dim_{GL(n)}\mu}
{(n\u\mu)\dim_{GL(N)}\mu}=
\dfrac{H(\mu)\dim_{GL(n)}\mu}
{(n\u\mu)\dim_{GL(N)}\mu}=
\frac 1{\dim_{GL(N)}\mu}\;.
$$
Thus, we obtain the same result as for $\SmNt$, which completes the proof.
\qed
\enddemo

For the third proof of Theorem 10.1 we need the following claim, which is
of independent interest.

\proclaim{\smc Proposition 10.2} Fix $n=1,2,\ldots$ and introduce a linear
map $\Ls(n)\to\Ls(n+1)$ making the following diagram commutative:
$$
\CD
\Zgln @>\operatorname{Avr}_{n,n+1}>> \frak{Z}(\gl(n+1))\\
@VVV @VVV \\
\Ls(n) @>>> \Ls(n+1)
\endCD \tag10.25
$$
{\rm(}here the vertical arrows are the canonical algebra isomorphisms
discussed in section {\rm 2)}. Given $f\in\Ls(n)$, let $f'\in\Ls(n+1)$
denote its image under that map. Then for each signature $\l$ for
$GL(n+1)$ the following relation holds
$$
f'(\l_1,\ldots,\l_{n+1})=\sum_{\nu\prec\l}
\frac{\dim_{GL(n)}\nu}{\dim_{GL(n+1)}\l}
f(\nu_1,\ldots,\nu_n).\tag10.26
$$
Moreover, $f'\in\Ls(n+1)$ is uniquely determined by these relations.
\endproclaim

\demo{Proof} Let $Z\in\Zgln$ and $Z'\in\frak{Z}(\gl(n+1))$ be the central
elements corresponding to $f$ and $f'$. Then, by virtue of \tht{10.9}, $Z'$ is
characterized by the relations
$$
\la Z',F\ra=\la Z,\Res_{n+1,n}F\ra,
\qquad F\in I(GL(n+1)).\tag10.27
$$
We may assume $F$ is an irreducible character of $GL(n+1)$, i.e., 
$F=gl(n+1)_\l$, where $\l$ is
an arbitrary signature. By the branching rule \tht{10.19}, 
$$
\Res_{n+1,n}gl(n+1)_\l=\sum_{\nu\prec\l}gl(n)_\nu.\tag10.28
$$
Hence \tht{10.27} can be rewritten as 
$$
\dim_{GL(n+1)}\l\cdot f'(\l_1,\ldots,\l_{n+1})=
\sum_{\nu\prec\l}\dim_{GL(n)}\nu\cdot f(\nu_1,\ldots,\nu_n),\tag10.29
$$
which is equivalent to \tht{10.26}. \qed
\enddemo

Note that in Theorem 10.1 we may assume without loss of generality that
$N=n+1$. It follows that the claim of Theorem 10.1 is equivalent to the
following family of relations: 
$$
\frac{s^*_{\mu|n+1}(\l_1,\ldots,\l_{n+1})}
{(n+1\u\mu)}=
\sum_{\nu\prec\l}\frac{\dim_{GL(n)}\nu}{\dim_{GL(n+1)}\l}
\cdot\frac{s^*_{\mu|n}(\nu_1,\ldots,\nu_n)}
{(n\u\mu)},\tag10.30
$$
where $n\ge\ell(\mu)$ and $\l$ is an arbitrary signature for $GL(n+1)$.

We call these the {\it coherence relations\/} for the shifted Schur
polynomials. 

\demo{Third Proof of Theorem 10.1} We aim to give a direct proof of
the coherence relations \tht{10.30}. We have 
$$
\gather
s^*_{\mu|n+1}(\l_1,\ldots,\l_{n+1})=
\frac{\det[(\l_i+n+1-i\f\mu_j+n+1-j)]}
{\prod_{i<j}(\l_i-\l_j+j-i)},
\qquad 1\le i,j\le n+1,\\
s^*_{\mu|n}(\nu_1,\ldots,\nu_n)=
\frac{\det[(\nu_i+n-i\f\mu_j+n-j)]}
{\prod_{i<j}(\nu_i-\nu_j+j-i)},
\qquad 1\le i,j\le n,\\
\dim_{GL(n+1)}\l=\frac{\prod_{i<j}(\l_i-\l_j+j-i)}
{n!(n-1)!\ldots0!},\\
\dim_{GL(n)}\nu=\frac{\prod_{i<j}(\nu_i-\nu_j+j-i)}
{(n-1)!\ldots0!},\\
\frac{(n+1\u\mu)}{(n\u\mu)}=
\prod_{i=1}^n\frac{\mu_i+n+1-i}{n-i+1}=
\frac 1{n!}\prod_{i=1}^n(\mu_i+n+1-i)
\endgather
$$
(in the last formula we have used the assumption $\ell(\mu)\le n$).

It follows that the coherence relations \tht{10.30} can be rewritten as 
$$
\gather
\det[(\l_i+n+1-i\f\mu_j+n+1-j)]_{1\le i,j\le n+1}=
\left(\prod_{j=1}^n(\mu_j+n+1-j)\right)\\
\cdot \sum_{\nu\prec\l}
\det[(\nu_i+n-i\f\mu_j+n-j)]_{1\le i,j\le n}. \tag10.31
\endgather
$$

Now consider the first determinant in \tht{10.31} and remark that its last
column is equal to $(1,\ldots,1)$, because $\mu_{n+1}=0$. 
For $i=1,\ldots,n$ let us subtract from
the $i$-th row the $i+1$-th one. Then we get the determinant of a
$n\times n$ matrix $A=[A_{ij}]$, where
$$
A_{ij}=(\l_i+n+1-i\f\mu_j+n+1-j)-(\l_{i+1}+n-i\f\mu_j+n+1-j).\tag10.32
$$

Note the following identity (which we already used above, see \tht{4.16})
$$
(b+1\f k)-(a\f k)=k\sum_{c=a}^b(c\f k-1),
\qquad b\ge a,\quad a,b\in\Z.\tag10.33
$$
By applying this identity with
$$
a=\l_{i+1}+n-i,\;\; b=\l_i+n-i,\;\; c=\nu_i+n-i,\;\; k=\mu_j+n+1-j,
$$
we get
$$
A_{ij}=(\mu_j+n+1-j)\sum_{\l_i\ge\nu_i\ge\l_{i+1}}
(\nu_i+n-i\f\mu_j+n-j).\tag10.34
$$
This gives an expansion of $\det A$ which coincides with the right--hand
side of \tht{10.31}. \qed
\enddemo

\proclaim{\smc Corollary 10.3} Let $\mu$ be a partition, $\ell(\mu)\le n$,
and let $\l$ be a signature for $GL(N)$, where $N>n$. Let $V_{\l|N}$ be
the irreducible $GL(N)$-module, indexed by $\l$. Then
$$
\frac{\tr_{V_{\l|N}}\Smn}{\dim_{GL(N)}\l}
=\frac{(n\u\mu)}{(N\u\mu)}s^*_\mu(\l).\tag10.35
$$
\qed
\endproclaim

\head 11. Combinatorial formula for $s^*$-functions \endhead

Let $\mu/\nu$ be a skew diagram. By definition, put
$$
(x\f \mu/\nu) = \prod_{\a\in\mu/\nu} (x-c(\a))\,. \tag11.1
$$
This is a generalization of the falling factorial
powers. It is clear that 
$$
(x\f\mu/\nu)=(x\f\mu)(x\f\nu)^{-1}.
$$
We have $(x\f\mu)=(x\f k)$ if $\mu=(k)$. According to \tht{10.20} we write 
$\mu\succ\nu$ if $\mu_i\ge\nu_i\ge\mu_{i+1}$, $i=1,2,\dots$.
Fix some $n$ and denote by $\RTab(\mu,n)$ the set of sequences  
$$
\mu=\mu^{(1)}\succ\mu^{(2)}\succ\dots\succ\mu^{(i)}\succ\dots
\succ\mu^{(n+1)}=\varnothing\,.\tag11.2
$$
Equivalently, elements of $\RTab(\mu,n)$ can be considered
as tableaux $T$ of shape $\mu$ whose entries $T(\a)$ (where $\a\in\mu$) 
belong
to $\{1,\dots,n\}$ and weakly decrease along each row and 
strictly decrease down each column. 
We call such tableaux {\it reverse tableaux} of shape 
$\mu$. 

There is a natural inclusion 
$$
\RTab(\mu,n)\subset\RTab(\mu,n+1)\,,\tag11.3
$$
where $\RTab(\mu,n)$ is identified with the set of sequences \tht{11.2} such
that $\mu^{(n+2)}=\varnothing$. 
Denote by $\RTab(\mu)$ the union of all the sets $\RTab(\mu,n)$,
$n=1,2,\ldots$. This set consists of all infinite sequences
$$
\mu=\mu^{(1)}\succ\mu^{(2)}\succ\dots\succ\mu^{(i)}\succ\dots\,, \tag11.4
$$
such that $\mu^{(i)}=\varnothing$ for all sufficiently large $i$.

\proclaim{\smc Theorem 11.1} We have
$$
\align
s^*_\mu(x_1,x_2,\ldots)& =
\sum_{\RTab(\mu)}
\prod(x_i\f\mu^{(i)}/\mu^{(i+1)})\\
&=\sum_{T\in\RTab(\mu)}
\prod_{\a\in\mu}(x_{T(\a)}-c(\a)) \,, \tag11.5
\endalign
$$
that is,
$$
\align
s^*_{\mu|n}(x_1,\ldots,x_n)& =
\sum_{\RTab(\mu,n)}
\prod(x_i\f\mu^{(i)}/\mu^{(i+1)})\\
&=\sum_{T\in\RTab(\mu,n)}
\prod_{\a\in\mu}(x_{T(\a)}-c(\a)) \tag11.6
\endalign
$$
for all $n\ge\ell(\mu)$.
\endproclaim

Formula \tht{11.5} is equivalent to the combinatorial presentation of
factorial Schur polynomials \tht{0.4},
$$
t_\mu(x_1,\ldots,x_n)=
\sum_{T\in\operatorname{Tab}(\mu,n)}
\prod_{\alpha\in\mu}(x_T(\alpha)-T(\alpha)-c(\alpha)+1),\tag11.7
$$
where $\operatorname{Tab}(\mu,n)$ denotes the set of ordinary tableaux
of shape $\mu$, see (6.16) in Macdonald \cite{M2} and Theorem 3.2 in
Chen--Louck \cite{CL}. To establish the equivalence of \tht{11.6} and
\tht{11.7}
it suffices to use the obvious relation  between both kinds of polynomials 
and reverse in \tht{11.7} the order of
variables (which is possible due to symmetry of the factorial Schur 
polynomials). 

Below we present two proofs of formula \tht{11.6}. The first one is 
essentially 
the proof of Chen and Louck \cite{CL}, rewritten in terms of the shifted
Schur polynomials; we give it for sake of
completeness. The second proof is based on the Characterization Theorem
and closely follows the argument of \cite{Ok1},
Proposition 3.8; the only difference is that we directly verify that the
right--hand side of \tht{11.6} is shifted symmetric.

Let us consider first the simplest case when $\mu=(k)$ and we have only
two variables $x$ and $y$. By virtue of the definition \tht{1.6} of 
the shifted
Schur polynomials formula \tht{11.6} reduces in this case to
the following elementary claim: 

\proclaim{\smc Lemma 11.2}
$$
\frac{(x+1\f k+1)-(y\f k+1)}{x-y+1} =
\sum_{l=0}^k (y\f l) (x-l\f k-l)  \,.\tag11.8
$$
\endproclaim

This is Lemma 2.1 in Chen--Louck \cite{CL}; we will use it in the
both proofs of the theorem. 

\demo{Proof of Lemma}
We have
$$
\multline
(x-y+1) \sum_{l=0}^k (y\f l) (x-l\f k-l) \\
=\sum_{l=0}^k \big((y\f l) (x-l\f k-l)(x+1-l) -
(y\f l) (x-l\f k-l)(y-l)\big) \\
=\sum_{l=0}^k(y\f l)(x+1\f k-l+1)-
\sum_{l=0}^k(y\f l+1)(x-l\f k-l)\,.
\endmultline
$$
It is easy to see that all summands cancel each other
except 
$$
(x+1\f k+1)-(y\f k+1)\,. \qed
$$
\enddemo 

\demo{First Proof of Theorem 11.1}
It is clear that the theorem is equivalent to the following
{\it branching rule\/} for $\s$-functions:
$$
\s_\mu(x_1,x_2,\dots)=\sum_{\nu\prec\mu} 
(x_1\f\mu/\nu)\, \s_\nu(x_2,x_3,\dots) \,, \tag11.9
$$
or equivalently,
$$
\s_\mu(x_1,\dots,x_n)=\sum_{\nu\prec\mu} 
(x_1\f\mu/\nu)\, \s_\nu(x_2,\dots,x_n) \,,\qquad n\ge\ell(\mu). \tag11.10
$$

Let us check \tht{11.10} using the definition \tht{1.6} of the polynomials
$s^*_{\mu|n}$. Consider the numerator of \tht{1.6},
$$
\det
\left[
\matrix
(x_1+n-1\f\mu_1+n-1)&\hdots&(x_1+n-1\f\mu_n)\\
\vdots&&\vdots\\
(x_n\f\mu_1+n-1)&\hdots&(x_n\f\mu_n)
\endmatrix
\right] \,. \tag11.11
$$
For all $j=1,\dots,n-1$ subtract from the $j$-th column
of \tht{11.11} the $(j+1)$-th column, multiplied by
$$
(x_1-\mu_{j+1}+j\f\mu_j-\mu_{j+1}+1) \,.
$$
Then for all $j<n$ the $(i,j)$-th entry of \tht{11.11} becomes
$$
\multline
(x_i+n-i\f\mu_{j+1}+n-j-1)
\big[
(x_i-\mu_{j+1}+j+1-i\f\mu_j-\mu_{j+1}+1) \\
- (x_1-\mu_{j+1}+j\f\mu_j-\mu_{j+1}+1) 
\big] \,. 
\endmultline \tag11.12
$$
In particular the first row of \tht{11.11} becomes
$$
(0,\dots,0,(x_1+n-1\f\mu_n)) \,.
$$

Let us apply the lemma, where we substitute
$$
\alignat2
x&=x_1-\mu_{j+1}+j-1,&\qquad k&=\mu_j-\mu_{j+1},\\
y&=x_i-\mu_{j+1}+j+1-i,&\qquad l&=\nu_j-\mu_{j+1}.
\endalignat
$$
Then we obtain
$$
\gather
x-y+1=x_1-x_i+i-1\\
(x-l\f k-l)=(x_1-\nu_j+j-1\f\mu_j-\nu_j)\\
(y\f l)=(x_i-\mu_{j+1}+j+1-i\f \nu_j-\mu_{j+1})\\
(x_i+n-i\f\mu_{j+1}+n-j-1)(y\f l)=(x_i+n-i\f \nu_j+n-j-1),
\endgather
$$
whence the entry \tht{11.12} equals
$$
-(x_1-x_i+i-1)
\sum_{\nu_j=\mu_{j+1}}^{\mu_j}
(x_1-\nu_j+j-1\f\mu_j-\nu_j)\,
(x_i+n-i\f\nu_j+n-j-1)\,.
$$
Therefore the determinant \tht{11.11} equals
$$
\prod_{1<i}(x_1-x_i+i-1) 
\sum_{\nu\prec\mu} 
(x_1\f\mu/\nu) \,
\det
\big[
(x_{i+1}+n-i-1\f\nu_j+n-j-1)
\big]_{i,j=1}^{n-1} \,, 
$$
which implies \tht{11.10}. \qed
\enddemo

\demo{Second Proof of Theorem 11.1} Denote the right--hand side of
\tht{11.6} by $\Sigma_{\mu|n}$ or by $\Sigma_\mu(x_1,\ldots,x_n)$ (below we
will also need an 
evident generalization of this notation to skew diagrams).

To apply Theorem 3.4 we have to check that, first, $\Sigma_{\mu|n}$ is a
shifted symmetric polynomial and, second, that it vanishes at certain
partitions and differs from $s^*_{\mu|n}$ in lower terms only.

For the first claim we shall suitably modify the well--known argument of 
Bender and Knuth \cite{BK}, which proves  symmetry of the combinatorial
formula for ordinary Schur polynomials. 

Let us show that $\Sigma_{\mu|n}$ is a shifted symmetric polynomial, i.e.,
$$
\Sigma_\mu(x_1,\ldots,x_i,x_{i+1},\ldots,x_n)=
\Sigma_\mu(x_1,\ldots,x_{i+1}-1,x_i+1,\ldots,x_n)\tag11.13
$$
for $i=1,\ldots,n-1$. 

Let us fix some $i$ and put $y_1=x_i$, $y_2=x_{i+1}$.
Consider the skew diagram 
$\nu=\mu^{(i)}/\mu^{(i+2)}$. It suffices to prove that the expression
$$
\Sigma_{\nu|2}=\Sigma_\nu(y_1,y_2)=
\sum_{T\in\RTab(\nu,2)}\prod_{\alpha\in\nu}
(y_{T(\alpha)}-c(\alpha))\tag11.14
$$
satisfies the shifted symmetry condition
$$
\Sigma_\nu(y_1,y_2)=\Sigma_\nu(y_2-1,y_1+1).\tag11.15
$$

Note that each column in $\nu$ may contain either 2 or 1 boxes. For
each column with 2 boxes, the entries of the boxes do not depend on $T$:
by the definition of reverse tableaux we always have to put 2 over 1.
Hence the contribution of any such column in the 
product in \tht{11.14} has the form $(y_2-c)(y_1-c+1)$, which is clearly a
shifted symmetric expression. Thus, we may strike out from $\nu$ all the
columns with 2 boxes. After this operation the shape $\nu$ becomes a
horizontal strip, and the property \tht{11.15} can be verified separately 
for each
row of $\nu$. This essentially means that we have reduced the problem to
the case $\nu=(k)$ when \tht{11.15} follows from Lemma 11.2.

Let us turn to the second claim. We have to show that if $\l$ is a
partition with $\ell(\l)\le n$ then
$\Sigma_{\mu|n}(\l)=0$ unless $\mu\subset\l$. 
In fact we will check a stronger fact: for each $T\in\RTab(\mu,n)$
$$
\prod_{\a\in\mu}(\l_{T(\a)}-c(\a)) = 0 \tag11.16
$$
unless $\mu\subset\l$. 

Let $\Pi_T$ denote the product in \tht{11.16}. Put
$\l_{(i,j)}=\l_{T(i,j)}$.
Since $T\in\RTab(\mu)$ we have
$$
\l_{(1,1)}\le\l_{(1,2)}\le\dots\le\l_{(1,\mu_1)} \,. \tag11.17
$$
If $\Pi_T(\l)\ne 0$ then
$$
\l_{(1,1)}\ne 0,\quad\l_{(1,2)}\ne1,\quad\l_{(1,3)}\ne2,\dots \tag11.18
$$ 
By \tht{11.17} and \tht{11.18} we have
$$
\l_{(1,1)}\ge1,\quad\l_{(1,2)}\ge2,\quad\l_{(1,3)}\ge3,\dots \tag11.19
$$ 
Again since $T\in\RTab(\mu)$ we have for each $i$
$$
T(1,i)<T(2,i)<\dots<T(\mu'_i,i), \tag11.20
$$
whence by \tht{11.19}
$$
i\le\l_{(1,i)}\le\l_{(2,i)}\le\dots\le\l_{(\mu'_i,i)} \tag11.21
$$
By virtue of \tht{11.20} and \tht{11.21}, for each $i$, in the diagram
$\l$, there are at least $\mu'_i$ rows of length $\ge i$, so that
$\l'_i\ge\mu'_i$. 
Thus,  $\Pi_T(\l)\ne 0$ implies $\mu\subset\l$. By 
Theorem 3.4, $\Sigma_{\mu|n}$ equals $s^*_{\mu}$
up to a constant factor. In order to see that this
factor equals 1 we can either compare the highest
terms  or calculate 
explicitly the unique non-vanishing summand in $\Sigma_{\mu|n}(\mu)$.
\qed
\enddemo

\proclaim{\smc Corollary 11.3} The complete shifted functions
$h^*_r=s^*_{(k)}$ and the elementary shifted functions $e^*_r=s^*_{(1^r)}$
($r=1,2,\ldots$) can be written as
$$
\gather
h^*_r(x_1,x_2,\ldots)=\sum_{1\le i_1\le\ldots\le i_r<\infty}
(x_{i_1}-r+1)(x_{i_2}-r+2)\ldots x_{i_r},\tag11.22\\
e^*_r(x_1,x_2,\ldots)=\sum_{1\le i_1<\ldots< i_r<\infty}
(x_{i_1}+r-1)(x_{i_2}+r-2)\ldots x_{i_r}.\tag11.23
\endgather
$$
\qed
\endproclaim

Our next aim is to interpret formula \tht{11.6} in terms
of Gessel-Viennot
nonintersecting lattice paths (see Gessel--Viennot \cite{GV}, Sagan
\cite{S}). It will be  
convenient for us to consider south-western paths instead of
more customary north-eastern ones. We define  a {\it path\/}
in $\Z^2$ as a sequence
$$
p=(p(l))\in\Z^2,\quad l=1,2,\dots \,,
$$
such that 
$$
p(l+1)-p(l)\in\{(-1,0),(0,-1)\} \,.
$$
According to these two possibilities we speak of a {\it horizontal\/}
step or a {\it vertical\/} step, respectively. We shall
say that the path $p$ ends at $(a,-\infty)$ if the 
$x$-coordinate of $p(l)$ equals $a$ for all
sufficiently large $l$. 

For a partition $\mu$, we denote by $P(\mu)$ the set of sequences $P=(p_i)$,
$i=1,2,\dots$, such that $p_i$ starts at $(\mu_i-i,-1)$,
ends at $(-i,-\infty)$, and $p_i\cap p_j=
\varnothing$ if $i\ne j$. Remark that only finitely many
paths $p_i$ are not purely vertical.

To each element of $\RTab(\mu)$
$$
\mu=\mu^{(1)}\succ\mu^{(2)}\succ\dots\succ\mu^{(i)}\succ\dots \tag11.24
$$
assign an element $P=(p_i)$ of $P(\mu)$
such that the vertical steps of the path $p_i$ are 
$$
\{(\mu^{(j)}_i-i,1-j),
(\mu^{(j)}_i-i,-j)\},\quad j=2,3,\dots \,.
$$
By the standard arguments of the Gessel-Viennot theory (see \cite{GV} or
\cite{S}) 
this is a bijection between the set
of sequences \tht{11.24} (or, equivalently, the set $\RTab(\mu)$)
and the set $P(\mu)$.

Next, the product over all boxes of $T\in\RTab$ in \tht{11.5}
is equal to the product 
over all horizontal steps $\{(a,-b),(a-1,-b)\}$ of $P$
of the following factors
$$
x_b-a\,. \tag11.25
$$
Denote this product by $\Pi_P(x)$. Then Theorem 11.1 can be restated as
follows: 

\proclaim{\smc Proposition 11.4} In the above notation
$$
\s_\mu(x)=\sum_{P\in P(\mu)} \Pi_P(x)\,. \tag11.26
$$
\qed
\endproclaim 

\example{\smc Example}
Suppose $\mu=(3,2)$ and $T$ is the following 
element of $\RTab(\mu)$
\picture 10 {} 
{
\hl 2 1 6
\hl 2 3 6
\hl 2 5 4
\vl 2 5 4
\vl 4 5 4
\vl 6 5 4
\vl 8 3 2
\boxit 2 3 {$3$}
\boxit 4 3 {$2$}
\boxit 6 3 {$2$}
\boxit 2 5 {$1$}
\boxit 4 5 {$1$}
\boxit 0 4 {$T=$}
}
\endpicture
\noindent
Then the corresponding nonintersecting paths are drawn
thick in the following picture
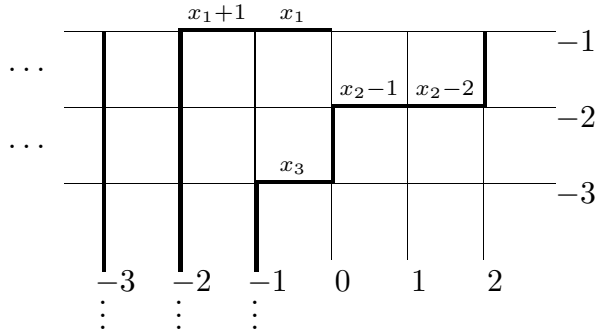
\picture 14 {}
{
\hl -1 1 12.9
\hl -1 3 12.9
\hl -1 5 12.9
\vl 0 7 6
\vl 2 7 6
\vl 4 7 6
\vl 6 7 6
\vl 8 7 6
\vl 10 7 6
\bvl 0 7.3 6.3
\bvl 2 7.3 6.3
\bhl 2 1 4
\bvl 4 7.3 2.3
\bhl 4 5 2
\bvl 6 5 2
\bhl 6 3 4
\bvl 10 3 2
\boxit -1 9.2 {$\vdots$}
\boxit 1 9.2 {$\vdots$}
\boxit 3 9.2 {$\vdots$}
\botit 2 1 {$\scriptstyle{x_1+1}$} 
\botit 4 1 {$\scriptstyle{x_1}$} 
\botit 4 5 {$\scriptstyle{x_3}$} 
\botit 6 3 {$\scriptstyle{x_2-1}$} 
\botit 8 3 {$\scriptstyle{x_2-2}$} 
\wr -0.2 7.8 {$-3$}
\wr 1.8 7.8 {$-2$}
\wr 3.8 7.8 {$-1$}
\wr 6.1 7.8 {$0$}
\wr 8.1 7.8 {$1$}
\wr 10.1 7.8 {$2$}
\wr 11.9 1.5 {$-1$}
\wr 11.9 3.5 {$-2$}
\wr 11.9 5.5 {$-3$}
\boxit -3 3 {$\hdots$}
\boxit -3 5 {$\hdots$}
}
\endpicture
\noindent 
On the horizontal steps are written the corresponding
factors \tht{11.25}.
\endexample

Denote by $P(\mu,n)$ 
the set of $n$-tuples $P=(p_i)$ of non-intersecting
lattice paths $p_i$ such that $p_i$ starts at $(\mu_i-i,-1)$,
ends at $(-i,-n)$, and $p_i\cap p_j=
\varnothing$ if $i\ne j$. Then it follows from \tht{11.6} 
that
$$
\s_\mu(x_1,\dots,x_n)=\sum_{P\in P(\mu,n)} \Pi_P(x)\,. \tag11.27
$$

\proclaim{\smc Proposition 11.5} We have the following Jacobi--Trudy--type
formula 
$$
\s_\mu(x_1,\dots,x_n)=
\det\big[h^*_{\mu_i-i+j}(x_1+j-1,\dots,x_n+j-1)\big]\,. \tag11.28
$$
\endproclaim

\demo{Proof} The standard arguments of the Gessel--Viennot theory, 
see \cite{GV} and \cite{S}, imply 
$$
\sum_{P\in P(\mu,n)} \Pi_P(x)=
\det\big[\Sigma(i,j)\big]\,, \tag11.29
$$
where 
$$
\Sigma(i,j) = \sum_p \Pi_p(x)\,,
$$
summed over all paths $p$ from $(\mu_i-i,-1)$ to
$(-j,-n)$. It follows from \tht{11.27} that
$$
\Sigma(i,j) = h^*_{\mu_i-i+j}(x_1+j-1,\dots,x_n+j-1)\,,
$$
provided $\mu_i-i+j\ge 0$ and $\Sigma(i,j)=0$ 
otherwise, which implies \tht{11.29}. \qed
\enddemo

Another possible way to prove \tht{11.28} is to derive it directly
from the determinantal formula \tht{1.6} for shifted Schur polynomials: this
can be done using the transformations of determinants, described (in 
case of ordinary Schur functions) in Littlewood \cite{L}, 6.3. 
Note that formula \tht{11.28} can easily be rewritten in terms of the
factorial Schur polynomials $t_\mu(x_1,\ldots,x_n)$. For these
polynomials, more general results are contained in Chen--Louck \cite{CL},
Theorem 4.1; Goulden--Hamel \cite{GH}, Theorem 4.2; Macdonald \cite{M2},
(6.7), and \cite{M1}, 2nd edition, Ch. I, section 3, Example 20 (c).

However, formula \tht{11.28} 
does not express the $\s$-functions in terms of the 
generators $h^*_k$ of the algebra $\Ls$, because the shift of
arguments is not defined in $\Ls$ (in other words, formula \tht{11.28} is not
stable as $n\to\infty$). In section 13 below we shall derive from \tht{11.28} 
a true Jacobi-Trudi formula for $\s$-functions.

\head 12. Generating series for $h^*$- and $e^*$-functions \endhead

Let $u$ be a formal variable. We shall deal
with series like
$$
c_0+
\frac{c_1}{u}+
\frac{c_2}{u(u-1)}+
\frac{c_3}{u(u-1)(u-2)}+\dots\,,\tag12.1
$$
which make sense because they can be rewritten as
formal power series in $u^{-1}$. 
Introduce the following generating series for the $h^*$- and
$e^*$-functions: 
$$ 
H^*(u)=\sum_{r\ge0}\frac{h^*_r(x_1,x_2,\dots)} 
{(u\f r)}\,,\qquad 
E^*(u)=\sum_{r\ge0}\frac{e^*_r(x_1,x_2,\dots)} 
{(u\f r)} \,.\tag12.2
$$ 

\proclaim{\smc Theorem 12.1}
$$ 
\align
H^*(u)&=
\prod_{i=1}^{\infty}\frac{u+i}{u+i-x_i}\,,\tag12.3\\ 
E^*(u)&=
\prod_{i=1}^{\infty}\frac{u-i+1+x_i}{u-i+1}\,.\tag12.4 
\endalign
$$ 
\endproclaim 

We present two proofs: one is a direct computation and another one is
based on the Characterization Theorem.

\demo{First Proof} Let us prove \tht{12.3}.  
If $n=1$ then \tht{12.3} takes the form (we write $x$ instead of $x_1$) 
$$ 
\sum_{r\ge0}\frac{(x\f r)}{(u\f r)}=\frac{u+1}{u+1-x}\,. \tag12.5
$$ 
This is a particular case of Gauss' formula for the value of the 
hypergeometric function $F(a,b;c;z)$ at $z=1$ 
(see, e.g., Whittaker--Watson \cite{WW}, Part II, Section 14.11); in our case
$a=-x$, $b=1$, $c=-u$. One also can check \tht{12.5} directly. Indeed, 
let us show by induction on $m=1,2,\dots$ that 
$$ 
\frac{u+1}{u+1-x}=\sum_{r=0}^m\frac{(x\f r)}{(u\f r)} 
+O\left(\frac 1{u^{m+1}}\right)\,. \tag12.6
$$ 
For $m=1$ this is trivial. Next for $m>1$ 
$$ 
\alignat2
\sum_{r=0}^m\frac{(x\f u)}{(u\f r)}&=1+\frac xu 
\sum_{r=0}^{m-1}\frac{(x-1\f r-1)}{(u-1\f r-1)}&&\\ 
&=1+\frac xu\left(\frac u{u+1-x}+O\left(\frac 1{u^m}\right)\right)
&\quad &\text{by assumption}\\ 
&=1+\frac x{u+1-x}+O\left(\frac 1{u^{m+1}}\right)&&\\
&=\frac{u+1}{u+1-x}+O\left(\frac 1{u^{m+1}}\right)\,.
\endalignat
$$ 
Thus we have checked \tht{12.3} for $n=1$. 
 
Now we use induction on $n$.
It follows from \tht{11.22} that 
$$ 
h^*_r(x_1,\ldots,x_n)=\sum\Sb p,q\ge0\\p+q=r\endSb  
h^*_p(x_1-q,\ldots,x_{n-1}-q)(x_n\f q) \tag12.7
$$ 
Therefore, 
$$ 
\align 
\sum_{r\ge0}\frac{h^*_r(x_1,\ldots,x_n)}{(u\f r)}
&= 
\sum_{p,q\ge0}\frac{h^*_p(x_1-q,\ldots,x_{n-1}-q)(x_n\f q)}{(u\f p+q)}\\ 
&=\sum_{q\ge0}\frac{(x_n\f q)}{(u\f q)} 
\sum_{p\ge0}\frac{h^*_p(x_1-q,\ldots,x_{n-1}-q)}{(u-q\f p)}\\ 
&=\sum_{q\ge0}\frac{(x_n\f q)}{(u\f q)}\prod_{i=1}^{n-1} 
\frac{u-q+i}{u+i-x_i}\,, \tag12.8
\endalign
$$ 
where the last equality holds by the inductive assumption.
It is easy to see that
$$
\frac{(u-q+n-1)\dots(u-q+1)}
{u(u-1)\dots(u-q+1)} =
\frac{(u+n-1)\dots(u+1)}
{(u+n-1)\dots(u-q+n)} \,.
$$
Therefore, \tht{12.8} equals
$$
\multline
\left(\prod_{i=1}^{n-1}\frac{u+i}{u+i-x_i}\right) 
\sum_{q\ge0}\frac{(x_n\f q)}{(u+n-1\f q)} 
\\
=\left(\prod_{i=1}^{n-1}\frac{u+i}{u+i-x_i}\right) 
\frac{u+n}{u+n-x_n} 
=\prod_{i=1}^n\frac{u+i}{u+i-x_i}\;,
\endmultline 
$$ 
as desired. This proves \tht{12.3}. 

The proof of \tht{12.4} is similar and even simpler.  It suffices to 
check the identity 
$$ 
\sum_{r=0}^n\frac{e^*_r(x_1,\ldots,x_n)}{(u\f r)}= 
\prod_{i=1}^n\frac{u-i+1+x_i}{u-i+1}.\tag12.9 
$$ 
We again use induction on $n$. For $n=1$ the identity is trivial. 
By virtue of \tht{11.23}, 
$$ 
e^*_r(x_1,\ldots,x_n)= 
e^*_r(x_1,\ldots,x_{n-1})+e^*_{r-1}(x_1+1,\ldots,x_{n-1}+1)x_n.\tag12.10 
$$ 
It follows 
$$ 
\sum_{r=0}^n\frac{e^*_r(x_1,\ldots,x_n)}{(u\f r)}= 
\sum_{r=0}^{n-1}\frac{e^*_r(x_1,\ldots,x_{n-1})}{(u\f r)}+ 
x_n\sum_{r=1}^n\frac{e^*_{r-1}(x_1+1,\ldots,x_{n-1}+1)}{(u\f r)}.\tag12.11 
$$ 
Using the relation $(u\f r)=u(u-1\f r-1)$, $r\ge1$, we can rewrite this as 
$$ 
\sum_{r=0}^{n-1}\frac{e^*_r(x_1,\ldots,x_{n-1})}{(u\f r)}+ 
\frac{x_n}u\sum_{r=0}^{n-1}\frac{e^*_r(x_1+1,\ldots,x_{n-1}+1)}{(u-
1\f r)}. 
\tag12.12
$$ 
By the inductive assumption, \tht{12.12} is equal to 
$$ 
\multline 
\prod_{i=1}^{n-1}\frac{u-i+1+x_i}{u-i+1}+\frac{x_n}u 
\prod_{i=1}^{n-1}\frac{u-i+1+x_i}{u-i}\\ 
=\left(1+\frac{x_n}{u-n+1}\right)\prod_{i=1}^{n-1} 
\frac{u-i+1+x_i}{u-i+1}=\prod_{i=1}^n\frac{u-i+1+x_i}{u-i+1}. 
\endmultline 
$$ 
\qed 
\enddemo 

Note that \tht{12.4} also follows from a relation proved by Macdonald, see
\cite{M2}, (6.5), or \cite{M1}, 2nd edition, Chapter I, Section 3, 
Example 20 (a).

For the second proof of Theorem 12.1 we need the following evident lemma.

\proclaim{\smc Lemma 12.2}
Let $f(u)$ be a polynomial of degree $\le k+1$ and write
$$
F(u)=\frac{f(u)}
{u(u-1)\dots(u-k)}\;. \tag12.13
$$
Then
$$
F(u)=
c_0+
\frac{c_1}{u}+
\frac{c_2}{u(u-1)}+ \dots +
\frac{c_{k+1}}{u(u-1)\dots(u-k)}\;, \tag12.14
$$
where $c_i$ are some coefficients.
\endproclaim

\demo{Proof} It suffices to write
$$
\multline
f(u)=c_{k+1}+c_k(u-k)+c_{k-1}(u-k)(u-k+1)+
\dots\\
+c_0 (u-k)(u-k+1)\dots(u-1)u\,.
\endmultline
$$
\qed
\enddemo

\demo{Second Proof of Theorem 12.1} Let us prove \tht{12.3}. 
It is clear that 
$$
\prod_{i=1}^\infty\frac{u+i}{u+i-x_i}=
c_0(x)+
\frac{c_1(x)}{u}+
\frac{c_2(x)}{u(u-1)}+
\frac{c_3(x)}{u(u-1)(u-2)}+\dots\,, \tag12.15
$$
where $c_i(x)$ are certain shifted symmetric functions in $x$. By Theorem
3.4, it suffices to check two claims: first, the highest term of $c_k$ is
equal to $h_k$, $k=0,1,2,\ldots$, and, second,
$$
c_k(\l)=0\qquad \text{for a partition $\l$ such that $|\l|<k$.}\tag12.16
$$

Put $x_i=u t \xi_i$ and let $u\to\infty$. Then the 
left--hand side  of \tht{12.1.53} turns into
$$
\prod_i (1-t\xi_i)^{-1} \,,
$$
which implies
the first claim.

Next, substitute $x=\l$ into \tht{12.15}. Fix an 
arbitrary $n\ge\ell(\l)$. Then
the left--hand side of \tht{12.15} can be written as
$$
\frac{(u+1)\ldots(u+n)}
{(u-(\l_n-n))\ldots(u-(\l_1-1))}.\tag12.17
$$
Note that  the factors in the denominator are of the form $u-r$, where
the $r$'s are pairwise distinct integers from $\{-n,-n+1,\ldots,
\l_1-1\}$.  If $r<0$ then the corresponding factor cancels with a certain
factor in the numerator. After all cancellations only factors with
$r\in\{0,\ldots,\l_1-1\}$ will remain (in the denominator). Hence the
product \tht{12.17} can be rewritten in the form \tht{12.13} 
with $k\le\l_1-1$, so
that, by Lemma 12.2, 
$$
c_{\l_1+1}(\l) = c_{\l_1+2}(\l) = \ldots = 0 \,.
$$
Thus, $c_k(\l)=0$ when $\l_1<k$, which implies our second claim.

For \tht{12.4} the argument is similar and even a bit simpler. The key
observation is that if $\l$ is a partition and $n=\ell(\l)$ then
$$
\prod_{i=1}^\infty\frac{u-i+1+\l_i}{u-i+1}=
\frac{(u+\l_1)\ldots(u-n+1+\l_n)}
{u(u-1)\ldots(u-n+1)},
$$
so that, by Lemma 12.2, the corresponding coefficients $c_k$ vanish at
$\l$ if 
$k>\ell(\l)$ and hence if $|\l|<k$.
\qed
\enddemo

\proclaim{\smc Corollary 12.3} The generating series for the $h^*$- and
$e^*$-functions satisfy the relation
$$
H^*(u) E^*(-u-1) = 1 \,. \tag12.18
$$
\qed
\endproclaim

Recall that in section 1 we have introduced certain generators $p^*_k$
of the algebra $\Ls$ (see \tht{1.12}). We have also introduced their
generating series $P^*(u)$, see \tht{4.3} and \tht{4.4}. By 
comparing \tht{4.4} and \tht{12.3} we see that
$$
P^*(u)=\frac{d}{du}\log H^*(u).\tag12.19
$$
Similarly, \tht{4.4} and \tht{12.4} imply 
$$
P^*(u)=-\frac{d}{du}\log E^*(-u-1).\tag12.20
$$
Now we would like to remark that this agrees with the following facts, 
proved in section 4 :
$$
\omega P^*(u)=P^*(-u-1),\qquad \omega H^*(u)=E^*(u).\tag12.21
$$

Thus, each of identities \tht{12.3}, \tht{12.4}, 
combined with the relation $\o H^*(u)=E^*(u)$, implies
the other. 

Besides $\{p^*_k\}$, there is one more family  of elements
in $\Ls$ which 
also may be viewed as an analogue of the Newton power sums. 
We denote these new elements by $p^\circ_k$, $k=1,2,\ldots$, and 
define them in terms of their generating series
$$
P^\circ(u)=\sum_{k\ge 1}p^\circ_k u^{-k}\tag12.22
$$
as follows
$$
1+u^{-1}P^\circ(u)=\frac{H^*(u)}{H^*(u-1)}=
\frac{E^*(-u)}{E^*(-u-1)}.\tag12.23
$$
Clearly, the highest term of $p^\circ_k$ (as that of $p^*_k$) is the
Newton power sum $p_k$. An advantage of the elements $p^\circ_k$  
is that we have from \tht{12.23} and the relation $\o H^*(u)=E^*(u)$
$$
\omega P^\circ(u)=-P^\circ(-u),\tag12.24
$$
so that
$$
\omega p^\circ_k=(-1)^{k-1}p^\circ_k, \qquad k=1,2,\ldots,\tag12.25
$$
just as for the customary power sums.

\example{\smc Remark 12.4} The generating series $E^*(u)$ is closely
related to the classical Capelli identity and the so called quantum
determinant for the Yangian of $\gln$, see Howe--Umeda \cite{HU} and
Nazarov \cite{N1}. As was shown by Nazarov \cite{N1}, the series $H^*(u)$
also appears in a Capelli  identity; it arises when  inverting the
quantum determinant. It should be mentioned that the form of the
generating series $H^*(u)$ and $E^*(u)$ was suggested  us by Nazarov's
work \cite{N1}. Finally, note that  
the series $P^\circ(u)$ appears in the computation of the square of 
the antipode in the Yangian, see section 6 in  
Molev--Nazarov--Olshanski \cite{MNO}.
\endexample

\head 13. Jacobi--Trudi formula for $\s$-functions \endhead
 
Define an automorphism $\phi$ of the algebra $\Ls$ by 
$$ 
\phi H^*(u)=H^*(u-1),\tag13.1  
$$ 
where  $H^*(u)$ is the generating series \tht{12.2} and 
$$ 
\phi H^*(u)=
\sum_{k\ge0} \phi(h^*_k)\,(u\f k)^{-1} \,. \tag13.2
$$ 
{F}rom the the identity 
$$ 
\frac 1{(u-1\f k)}=\frac 1{(u\f k)}+\frac k{(u\f k+1)} \tag13.3
$$ 
it follows that
$$ 
\phi(h^*_k)=h^*_k+(k-1)h^*_{k-1},\quad k=1,2,\dots\,.\tag13.4 
$$ 
Iterating \tht{13.4} we obtain
$$ 
\phi^r(h^*_k)=\sum_{i=0}^r \binom{r}{i} (k-1\f i)\, h^*_{k-i}\;,
\qquad r=1,2,\ldots\tag13.5
$$

{F}rom \tht{12.18} it follows that 
$$
\phi E^*(u) = E^*(u+1) \,. \tag13.6
$$ 
Therefore, by the very definition of $E^*(u)$, see \tht{12.2}, 
$$
\phi^{-1}(e^*_k)=e^*_k+(k-1)e^*_{k-1},\quad r=1,2,\dots\,.\tag13.7 
$$
and
$$
\phi^{-r}(e^*_k)=\sum_{i=0}^r \binom{r}{i} (k-1\f i)\, e^*_{k-i}\,.\tag13.8
$$

We have the following analogues of the classical Jacobi--Trudi and
N\"agelsbach--Kostka formulas.

\proclaim{\smc Theorem 13.1}
$$
\align
\s_\mu &= \det\left[\phi^{j-1} h^*_{\mu_i-i+j} 
\right]_{1\le i,j\le l} \tag13.9\\
\s_\mu &= \det\left[\phi^{1-j} e^*_{\mu'_i-i+j} 
\right]_{1\le i,j\le m}\,\tag13.10
\endalign
$$
for arbitrary $l,m$ such that $l\ge\ell(\mu)$ and $m\ge\mu_1$.
\endproclaim

\demo{Proof} By virtue of \tht{13.4} and \tht{13.7}, 
the action of $\phi$ on the
$h^*_k$'s is the same as that of $\phi^{-1}$ on the $e^*_k$'s. Hence, by
duality \tht{4.17}, it suffices to check one of these two formulas. We shall
deduce formula \tht{13.9} from \tht{11.28}. 

We can assume that we have a finite number of variables $x_1,\dots,x_n$.
Define in the space $\Ls(n)$ the `shift operator' $\tau$, 
$$ 
[\tau f](x_1,\ldots,x_n)=f(x_1+1,\ldots,x_n+1). \tag13.11
$$ 
Then 
$$ 
\align 
\tau^rH^*(u)&=
\prod_{i=1}^n\frac{u+i}{u+i-x_i-r}\\ 
&=\prod_{i=1}^n\frac{u+i}{u+i-r}\prod_{i=1}^n\frac{u-r+i}{u-r+i-x_i}\\ 
&=\frac{(u+n\f r)}{(u\f r)}H^*(u-r)\,.\tag13.12
\endalign
$$ 
Next consider the following `lowering' operator, acting in the linear span
of the polynomials $h^*_k\in\Ls(n)$, 
$$ 
\si : h^*_k\mapsto h^*_{k-1}\,.\tag13.13
$$ 
It is readily seen that 
$$ 
\si^r H^*(u)=\frac 1{(u\f r)}H^*(u-r)\,.\tag13.14
$$ 
Let us view $\tau$ as an operator in the linear span of the $h^*_k$. 
Clearly $\tau$ and $\si$ commute. 

Now we claim that for any $r=1,2,\dots$ 
$$ 
H^*(u-r)=\tau^rH^*(u)+\sum_{i=1}^r c_i\tau^{r-i}\si^iH^*(u)\,,\tag13.15
$$ 
where $c_i$ are certain number coefficients
which depend on $r$ and $n$. Indeed, by \tht{13.12} and \tht{13.14} 
$$ 
\tau^{r-i}\si^i H^*(u)=
\frac{(u+n\f r-i)}{(u\f r)}H^*(u-r) 
$$ 
and hence \tht{13.15} reduces to the evident identity
$$ 
(u\f r)=(u+n\f r)+\sum_{i=1}^r c_i\,(u+n\f r-i)\,.
$$ 
 
Now we are in a position to derive \tht{13.9} from \tht{11.28}. 
Abbreviate $h^*_r=h^*_r(x_1,\ldots,x_n)$.
We have to prove that 
$$ 
\det[\phi^{j-1}h^*_{\mu_i-i+j}]_{1\le i,j\le l} =
\det[\tau^{j-1}h^*_{\mu_i-i+j}]_{1\le i,j\le l}\,. \tag13.16
$$ 
To do this we shall show that for any $j=2,3,\ldots,l$, the $j$-th column 
of the matrix in the left--hand side of \tht{13.16} is equal to 
the $j$-th column of the 
matrix in the right--hand side of \tht{13.16} plus a linear combination 
of the preceding columns.  
 
Fix some $j$ and denote by $H_j$ the $j$-th column of the 
matrix $[h^*_{\mu_i-i+j}]$. Then the 
first $j$ columns of the matrix 
$[\tau^{j-1}h^*_{\mu_i-i+j}]$ equal
$$
\si^{j-1}H_j, \tau\si^{j-2} H_j,\dots,\tau^{j-1} H_j\,.
$$ 
By \tht{13.15}
$$ 
\phi^{j-1}H_j=
\tau^{j-1}H_j+\sum_{i=1}^{j-1} c_i\tau^{j-i-1}\si^i H_j \,.
$$ 
This concludes the proof. \qed 
\enddemo 

\example{\smc Remark 13.2} 
The equivalence of \tht{13.9} and \tht{13.10} can also be 
deduced from a general result due to Macdonald (\cite{M2}, Section 9, or
\cite{M1}, 2nd edition, Chapter I, Section 3, Example 21). Macdonald's
argument also 
implies that for the $s^*$-functions there exists a precise analogue of
the classical Giambelli formula: write $\mu$ in Frobenius notation, 
$$
\mu=(\alpha_1,\ldots,\alpha_r\;\vert\;\beta_1,\ldots,\beta_r);
$$
then
$$
s^*_\mu=\det[s^*_{(\alpha_i\;\vert\;\beta_j)}]\;.\tag13.17
$$
\endexample

\head 14. Special symmetrization\endhead

In this section we discuss a relation between the results of
section 6 and a certain linear
isomorphism
$$
\si: S(\gln)\to\Ugln\,.\tag14.1
$$
This isomorphism was introduced by Olshanski \cite{O2} in
connection with infinite--dimensional Lie algebras and Yangians. Then
$\si$ was used in Kerov--Olshanski \cite{KO}; there it
was called the {\it special symmetrization\/}. We also present some new
results about the map $\si$. 

The construction of the special symmetrization is based on the fact that
the Lie algebra structure of $\gln$ is obtained from the associative
algebra structure on the space of $n\times n$ matrices.  To define $\si$
let us view the enveloping algebra $\Ugln$ as the algebra of the (complex)
distributions on the real Lie group $GL(n,\R)$, supported at the unity,
and, similarly, identify the symmetric algebra $S(\gln)$ with the algebra
of (complex) distribution on the additive Lie group $\gln$, supported at
zero. Further, introduce a local chart on $GL(n,\R)$, centered at the
unity, by using the map $X\to 1+X$ from a neighborhood of $0\in\gln$ to
$GL(n,\R)$. This chart allows us to identify the both algebras of
distributions as vector spaces, and we take this identification as the map
$\si$. (This definition is equivalent to that given in \cite{O2}, 2.2.11.)

Note that $\si$ differs from the customary symmetrization map
$S(\gln)\to\Ugln$ (see, e.g., Dixmier \cite{D}, 2.4.6) but shares several
its properties. 
In particular, $\si$ preserves highest terms and commutes with the adjoint
action of the group $GL(n)$. 

We maintain the notation of sections 2 and 6. Let us use the map $R$ to
identify 
$\Ugln$ with the algebra of left--invariant differential operators on the
space of $n\times n$ matrices (see \tht{6.4}, where we shall assume $m=n$).
Then we have  a surprisingly simple formula 
$$
R(\si(E_{i_1j_1}\ldots E_{i_kj_k}))=
\sum_{\alpha_1,\ldots,\alpha_k=1}^n
x_{\alpha_1i_1}\ldots x_{\alpha_ki_k}
\p_{\a_1j_1}\ldots\p_{\a_kj_k}\tag14.2
$$
(see \cite{O2}, Lemma 2.2.12).

Since $\si$ is $GL(n)$-equivariant, it establishes a linear isomorphism
$$
\si: I(\gln)\to\Zgln\,. \tag14.3
$$
We recall that $I(\gln)$ stands for the subalgebra of $GL(n)$-invariants
in $S(\gln)$. Recall also that we have defined natural bases
$\{S_{\mu|n}\}\subset I(\gln)$ and $\{\Smn\}\subset\Zgln$, see \tht{2.9} and
Definition 2.3; here $\mu$ ranges over the set of partitions with
$\ell(\mu)\le n$.

\proclaim {\smc Theorem 14.1} The special symmetrization $\si$ takes each
$S_{\mu|n}\in I(\gln)$ to $\Smn\in\Zgln$.
\endproclaim

\demo{Proof} It suffices to prove that $R(\si(S_{\mu|n}))$ coincides with
$R(\Smn)$.  

Let $k=|\mu|$. By Proposition 2.4 and formula \tht{14.2},
$R(\si(S_{\mu|n}))$ is equal to
to
$$
(k!)^{-1}\sum_{\a_1,\ldots,\a_k=1}^n\sum_{i_1,\ldots,i_k=1}^n
\sum_{s\in S(k)} \chm(s)\cdot 
x_{\a_1i_1}\ldots x_{\a_ki_k}\p_{\a_1i_{s(1)}}\ldots
\p_{\a_ki_{s(k)}}.\tag14.4
$$

Next, by the identity \tht{6.5}, proved in Corollary 6.8, $R(\Smn)$ is the
differential operator $\Delta^{(n,n)}_\mu$, which is given by \tht{6.6},
where we have to set $m=n$:
$$
\Delta^{(n,n)}_\mu=
(k!)^{-1}\sum_{i_1,\ldots,i_k=1}^n\sum_{j_1,\ldots,j_k=1}^n
\sum_{s\in S(k)} \chm(s)\cdot 
x_{i_1j_1}\ldots x_{i_kj_k}\p_{i_{s(1)}j_1}\ldots\p_{i_{s(k)}j_k}.\tag14.5
$$
Now it is easy to see that the both expressions, \tht{14.4} and 
\tht{14.5}, are
the same (it suffices to reorder the derivatives in \tht{14.5} and then change
the notation of the subscripts). This completes the proof. \qed
\enddemo

There is a formula for the inverse map $\si^{-1}:\Ugln\to S(\gln)$,
obtained in \cite{O2}, 2.2.13. Let us fix $k=1,2,\ldots$ and denote by
$\Xi(k)$ the set of 
all partitions $\xi$ of the set $\{1,\dots,k\}$ into disjoint subsets
$\{i_1<i_2<\dots\}$, $\{j_1<j_2<\dots\}$, $\dots$, called the {\it
clusters\/} of $\xi$. Suppose we have a monomial 
$A_1\circ\dots\circ A_k\in\Ugln$ (we write ``$\circ$'' for the
multiplication in $\Ugln$ to distinguish it from the multiplication in
$S(\gln)$). Put 
$$
\Pi_\xi(A_1,\ldots, A_k)=
\la A_{i_1} A_{i_2} \ldots \ra \cdot
\la A_{j_1} A_{j_2} \ldots \ra \cdot 
\ldots \; \in S(\gln),\tag14.6
$$
where
$\la A_{i_1} A_{i_2} \dots \ra\in\gln$ is the ordinary matrix product of
$A_{i_1}$, $A_{i_2}$, $\dots$ (that is, the product in the associative
algebra of $n\times n$ matrices). Then
$$
\si^{-1} (A_1\circ\ldots\circ A_k)
=\sum_{\xi\in\Xi(k)} \Pi_\xi(A_1,\dots,A_k) \,. \tag14.7 
$$

For instance,
$$
\alignat2
&\si^{-1} (A_1)&=& A_1\,,\\
&\si^{-1} (A_1\circ A_2)&=& A_1 A_2+ \la A_1 A_2\ra\,,\\
&\si^{-1} (A_1\circ A_2\circ A_3)&=& 
A_1 A_2 A_3+ \la A_1 A_2\ra A_3+
\la A_1 A_3\ra A_2 +\\
&&&\la A_2 A_3\ra A_1+
\la A_1 A_2 A_3\ra\,.
\endalignat
$$

We call \tht{14.7} the {\it cluster formula\/}. Note that this formula is
used by Okounkov \cite{Ok2} 
in computation of the quantum immanants. There
it is 
also remarked that the cluster formula can be derived from the classical
Wick formula. 

Our next aim is to inverse the cluster formula and so find an explicit
expression for the special symmetrization.

For $A\in \gln$ we write $A^k$ for the $k$-th power of $A$ in $S(\gln)$
and $\la A\ra^k$ for its $k$-th power in the matrix algebra (so that $\la
A\ra^k\in\gln$). Recall a standard notation: if $\l$ is a partition,
$\l=(1^{r_1}2^{r_2}\ldots)$, then
$$
z_\l=1^{r_1}r_1!2^{r_2}r_2!\ldots\;. \tag14.8
$$

\proclaim{\smc Theorem 14.2} For any $A\in\gln$
$$
\si(A^k)=\sum_{\l\vdash k}(-1)^{k-\ell(\l)}\frac{k!}{z_\l}
\la A\ra^{\l_1}\circ\la A\ra^{\l_2}\circ\ldots\;.\tag14.9
$$
\endproclaim

Using polarization we easily obtain from \tht{14.9} the following {\it 
inversion formula\/} (with respect to \tht{14.7}). 

\proclaim{\smc Theorem 14.3} For any $A_1,\ldots,A_k\in\gln$
$$
\multline
\si(A_1\ldots A_k)\\
=\sum_{s\in S(k)}\sum_{\l\vdash k}
(-1)^{k-\ell(\l)}z_\l^{-1}\la A_{s(1)}\ldots A_{s(\l_1)}\ra\circ
\la A_{s(\l_1+1)}\ldots A_{s(\l_1+\l_2)}\ra\circ\ldots\;.
\endmultline\tag14.10
$$
\qed
\endproclaim

Note that the first version of the inversion formula was obtained by
Postnikov \cite{P}. He did not use the polarization argument, and his
formula involved an extra symmetrization.

\demo{Proof of Theorem 14.2} We remark first that our inversion problem 
makes sense for
any associative algebra. Indeed, let $M$ be such an algebra, also viewed
as a Lie algebra with bracket $[A,B]=AB-BA$. Then formula \tht{14.7} defines
a linear isomorphism 
$$
\si^{-1}: \U(M)\to S(M).\tag14.11
$$
Indeed, it is readily verified that for $i=1,\ldots,k-1$
$$
\gather
\si^{-1}(A_1\circ\ldots\circ A_i\circ A_{i+1}\circ\ldots\circ A_k)-
\si^{-1}(A_1\circ\ldots\circ A_{i+1}\circ A_i\circ\ldots\circ A_k)\\
=
\si^{-1}(A_1\circ\ldots\circ [A_i, A_{i+1}]\circ\ldots\circ A_k),\tag14.12
\endgather
$$
so that the map $\si^{-1}$ is correctly defined.

The next observation is that without loss of generality we may assume $M$
to be commutative. Indeed, given $A\in M$ we consider the associative
subalgebra $M_A$ in $M$, spanned by the powers of $A$, and remark that
$\si^{-1}$ maps $\U(M_A)$ onto $S(M_A)$. So we may replace $M$ by $M_A$
and reduce the problem to the case when $M$ is a commutative associative
algebra and the corresponding Lie structure is trivial. Then $\U(M)$ is
canonically identified with $S(M)$, so we may view $\si$ as a linear
isomorphism of the vector space $S(M)=\U(M)$.

The set $\Xi(k)$ is partially ordered: 
$\xi_1 \le \xi_2$ if each cluster of $\xi_2$ is a subset of
a cluster of $\xi_1$. Note that this order is inverse to
the order usually used in combinatorics. The partition
$$
\1=\{\{1\},\{2\},\dots,\{k\}\}
$$
is the maximal element of $\Xi(k)$. 

It follows from \tht{14.7}
that
$$
\si^{-1} \big(\Pi_\xi(A,\dots,A)\big) = \sum_{\zeta \le \xi}
\Pi_\zeta (A,\dots,A) \,. \tag14.13
$$
Hence we can use the M\"obius inversion in the poset $\Xi(k)$
to obtain a formula for $\si$. The M\"obius function $\mu$ of
the poset $\Xi(k)$ is well-known (see Stanley \cite{St}, section 3.10.4).
In particular,
$$
\mu(\xi,\1)=(-1)^{k-\ell(\lambda)} \prod (\lambda_i-1)!\,, \tag14.14
$$ 
where $\lambda_i$ are the cardinalities of the clusters of $\xi$.
Therefore
$$
\align
\si \big(A\cdot \dots \cdot A\big) &=
\si \big(\Pi_{\1}(A,\dots,A)\big) \\
&=
\sum_{\xi \in \Xi(k)} 
(-1)^{k-\ell(\lambda)} \prod (\lambda_i-1)!\,
\Pi_\xi (A,\dots,A) \,. \tag14.15
\endalign
$$
It is easy to see that there are 
$$
\frac {k!}{z_\lambda \prod (\lambda_i-1)!}
$$
elements of $\Xi(k)$ with cardinalities of clusters
equal to $\lambda$, which implies \tht{14.9}. \qed
\enddemo

\example{\smc Example 14.4} Suppose $A\in\gln$ is a matrix projection,
$\la A\ra^2=A$. Then
$$
\si(A^k)=A\circ(A-1)\circ\ldots\circ(A-k+1).\tag14.16
$$
\endexample

\example{\smc Remark 14.5} In the right--hand side 
of \tht{14.9}, for any $\l$,
the factors 
$$
\la A\ra^{\l_1}, \la A\ra^{\l_2}, \ldots
$$ 
can be rewritten in
an arbitrary order. After polarization we can obtain in this way a number
of various formulas
which differ from \tht{14.10} but are equivalent to it.
\endexample

\example{\smc Remark 14.6} By applying formula \tht{14.10} to the expression
\tht{2.10} for the elements $S_{\mu|n}$  we obtain, by virtue of Theorem 
14.1, a certain explicit expression for the quantum immanants
$\S_{\mu|n}$. However, this expression seems to be more complicate than
the expression obtained in the papers \cite{Ok1}, \cite{N2}, \cite{Ok2}.
\endexample

We conclude with one more interpretation of the special symmetrization.
Note that, by virtue of \tht{14.2}, the linear isomorphism $\si:
S(\gln)\to\Ugln$ can be obtained via realizing both $S(\gln)$ and $\Ugln$ by
differential operators on $n\times n$ matrices. Now we aim to describe
a different realization which leads to the same result.

As in the proof of Theorem 14.2 we shall work with an arbitrary
associative algebra $M$, also viewed as a Lie algebra.

First, we form the algebra $\tM=\C1\oplus M$, where 1 is a formally
adjoint unity. For each $k=1,2,\ldots$ we consider the algebra
$\tM^{\otimes k}$. Its elements are linear combinations of the tensors
$A_1\otimes\ldots\otimes A_k$, where each $A_i$ belongs either to $\C1$ or
to $M$. We equip $\tM$ with a filtration by setting
$$
\deg(A_1\otimes\ldots\otimes A_k)=\operatorname{card}\{i\mi A_i\in
M\}.\tag14.17 
$$
Note that this filtration is compatible with the algebra structure. 

Next, for each $k$ we define a projection $\tM^{\otimes
(k+1)}\to\tM^{\otimes k}$ by setting
$$
A_1\otimes\ldots\otimes A_{k+1}\mapsto
\cases \alpha A_1\otimes\ldots\otimes A_k, &\text{if $A_{k+1}=\a1$,}\\
0,&\text{if $A_{k+1}\in M$.}
\endcases\tag14.18
$$
Let $\Omega(M)$ be the projective limit of the algebras $\tM^{\otimes k}$,
taken with respect to the projections \tht{14.18} in the category of filtered
algebras. 

Finally, we define the embeddings
$$
\U(M)\to\Omega(M),\qquad S(M)\to\Omega(M)\tag14.19
$$
as follows.

Given $A\in M$,  put
$$
A^{(i)}=1^{\otimes(i-1)}\otimes A\otimes 1^{\otimes\infty},
i=1,2,\ldots\;.\tag14.20 
$$
Then the first arrow in \tht{14.19} is an algebra morphism, 
defined on elements
$A\in M$ by
$$
A\mapsto \sum_{i=1}^\infty A^{(i)}\tag14.21
$$
(note that the right--hand side of \tht{14.11} is a well--defined element of
$\Omega(M)$). 

As for the second arrow in \tht{14.19}, we define it on monomials in the
symmetric algebra as follows:
$$
A_1\ldots A_k\mapsto \sum_{i_1,\ldots,i_k}
A_1^{(i_1)}\ldots A_k^{(i_k)},\tag14.22
$$
summed over all $k$-tuples of {\it pairwise distinct\/} indices. We again
note that the right--hand side is well--defined in $\Omega(M)$.

\proclaim{\smc Proposition 14.7} The mappings \tht{14.19} defined above are
embeddings with the same image, so that they define a linear isomorphism
$S(M)\to\U(M)$. This isomorphism coincides with the map $\si$, defined by
the cluster formula \tht{14.7}. 
\endproclaim

\demo{Proof} This follows at once from \tht{14.7}. \qed
\enddemo

\head 15. Concluding remarks \endhead

\subhead 1. The map $\varphi: \L\to\Ls$\endsubhead By definition, this is
a linear isomorphism such that
$$
\varphi(s_\mu)=s^*_\mu\qquad \text{for each partition $\mu$.}\tag15.1
$$
In case of finitely many variables $\varphi$ turns into the linear
isomorphisms 
$$
\varphi:\L(n)\to\Ls(n),\qquad n=1,2,\ldots,\tag15.2
$$ 
such that 
$$
\varphi(s_{\mu|n})=s^*_{\mu|n},\qquad \ell(\mu)\le n.\tag15.3
$$

Equivalently, by Theorem 14.1, the map $\varphi$ can be defined via the
commutative diagram
$$
\CD
I(\gln) @>\si>> \Zgln\\
@VVV     @VVV\\
\L(n) @>\varphi>> \Ls(n)
\endCD \tag15.4
$$
where $\si$ is the special symmetrization and the vertical arrows are the
canonical isomorphisms discussed in section 2. Thus, the special 
symmetrization is a natural extension of
the map $\varphi$. 

\subhead 2. The basis $\{p^\#_\mu\}\subset\Ls$\endsubhead 
Let $p_\mu\in\L$ denote the ordinary power sum functions, where $\mu$ is
an arbitrary partition. We define the elements $p^\#_\mu\in\Ls$ by
$$
p^\#_\mu=\varphi(p_\mu).\tag15.5
$$
In particular,
$$
p^\#_k=\varphi(p_k),\qquad k=1,2,\ldots.\tag15.6
$$
Note that
$$
\omega(p^\#_k)=(-1)^{k-1}p^\#_k, \tag15.7
$$
where $\omega:\Ls\to\Ls$ is the involution of section 4. 

Thus, in the algebra $\Ls$, we have three different families of generators,
$\{p^*_k\}$, $\{p^\circ_k\}$, and $\{p^\#_k\}$; each of them can be
regarded as an analogue of Newton power sums $p_k$.

Recall the following fundamental identity in the theory of 
symmetric functions,
$$
p_\mu=\sum_{\l\vdash k}\chi^\l_\mu s_\l, \tag15.8
$$
where $k=|\mu|$, $\chi^\l$ is the irreducible character of the symmetric
group $S(k)$, indexed by $\l$, and $\chi^\l_\mu$ is its value on any
permutation with cycle--type $\mu$. It follows
$$
p^\#_\mu=\sum_{\l\vdash k}\chi^\l_\mu s^*_\l.\tag15.9
$$

\subhead 3. The functions $p^\#_\mu$ and a Capelli--type identity for the
Schur--Weyl duality \endsubhead  
Let $\rho=(\rho_1,\ldots,\rho_m)$ be a partition of $k$ and let $l\ge k$
and $n$ be natural 
numbers. As in section 7 we consider the tensor space 
$(\C^n)^{\otimes l}$ as a bimodule over $GL(n)$ and $S(l)$. In the paper
\cite{KO} there were introduced central elements
$$
a_{\rho,l}\in Z(S(l)), \qquad A_{\rho,n}\in\Zgln, \tag15.10
$$
such that, in notation \tht{7.2} and \tht{7.3},
$$
\tau_{GL(n)}(A_{\rho,n})=\tau_{S(l)}(a_{\rho,l}) \tag15.11
$$
(\cite{KO}, Theorem 2).

The element $a_{\rho,l}$ is defined as a sum of products of cycles,
$$
a_{\rho,l}=\sum_{i_1,\ldots,i_k}
(i_1,\ldots,i_{\rho_1})(i_{\rho_1+1},\ldots,i_{\rho_1+\rho_2})\ldots
(i_{\rho_i+\ldots+\rho_{m-1}+1},\ldots,i_k),\tag15.12
$$
where summation is taken over all $k$-tuples of pairwise distinct
indices from the set $\{1,\ldots,n\}$. That is, $a_{\rho,l}$ is a scalar 
multiple
of the conjugacy class in $S(l)$ corresponding to the partition $\rho\cup
1^{l-k}$.

The element $A_{\rho,n}$ is defined via the special
symmetrization map  as follows:
$$
A_{\rho,n}=\sigma(P_{\rho|n}),\tag15.13
$$
where the element $P_{\rho|n}\in I(\gln)$, viewed as an invariant
polynomial function 
on $\gln$, is defined by
$$
P_{\rho|n}(X)=\tr X^{\rho_1}\cdot\tr X^{\rho_2}
\cdot\ldots\cdot\tr X^{\rho_m}, \qquad X\in\gln. \tag15.14
$$
In other words, under the identification $I(\gln)=\L(n)$, $P_{\rho|n}$
becomes the power sum function $p_\rho$ in $n$ variables.

For a partition $\l\vdash l$, let $f_\rho(\l)$ be the
eigenvalue of the 
central element $a_{\rho,l}$ in $W_\l$ (the irreducible $S(l)$-module,
indexed by $\l$). By virtue of \tht{15.11} and the Schur--Weyl duality
\tht{7.1}, $f_\rho(\l)$ also is the eigenvalue of $A_{\rho,n}$ in $V_{\l|n}$,
the irreducible polynomial $\gln$-module corresponding to $\l$ (we assume
$\ell(\l)\le n$). It follows (\cite{KO}, Proposition 3) that
$f_\rho(\cdot)$ is a shifted symmetric function with highest term
$p_\rho$. 

Now we remark that in notation of section 7
$$
a_{\rho,l}=\sum_{\mu\vdash k}\chi^\mu_\rho \Ind\chi^\mu/(l-k)! \tag15.15
$$
and
$$
A_{\rho,n}=\sum_{\mu\vdash k}\chi^\mu_\rho S_{\mu|n}. \tag15.16
$$
Therefore, the relation \tht{15.11} is equivalent to the relation \tht{7.8} of
Theorem 7.1. 

An important consequence of this is that
$$
f_\rho=\varphi(p_\rho),\tag15.17
$$
or, in notation \tht{15.5}, 
$$
f_\rho=p^\#_\rho. \tag15.18
$$

Formula \tht{15.17} first appeared in Macdonald's letter \cite{M3}
addressed to Olshanski (this letter was written as a comment to
the papers \cite{KO} and \cite{O3}). There Macdonald introduced the map
$\varphi$ and gave a simple direct proof of \tht{15.17}. 
At that moment we were
already aware of the Characterization Theorem, of formula \tht{6.5} for
the differential operators corresponding to the quantum immanants, 
and of Theorem 14.1, which 
is essentially equivalent to the relation \tht{15.17}. However,
the elegant argument of Macdonald 
suggested us formula \tht{8.3} for the dimension of skew Young diagrams
(Theorem 8.1).

Theorem 8.1 is one more result  equivalent 
to \tht{15.17}. As we already noted in section 8, the second
proof of Theorem 8.1 is nothing but a slight modification of
Macdonald's proof of the relation \tht{15.17}.

Note also that the eigenvalue $f_\rho(\l)$ can be written in 
the `Gibbsian form'
$$
f_\rho(\l)=\frac{(s_\l,p_\rho e^{p_1})}{(s_\l,e^{p_1})}.\tag15.19
$$
This fact was communicated by Kerov to Olshanski when working 
on the paper \cite{KO}; it is close to expressions used by 
Macdonald \cite{M3}. 

\subhead 4. Character values on small cycle--types \endsubhead
By the very definition of $f_\rho(\l)$, we have
$$
f_\rho(\l)=\frac{(l\f k)}{\dim\l}\chi^\l_{\rho\cup 1^{l-k}}. \tag15.20
$$
Therefore, in view of \tht{15.18} and
\tht{15.9} we obtain
$$
\chi^\l_{\rho\cup 1^{l-k}}=
\frac{\dim\l}{(l\f k)} 
\sum_{\mu\vdash k}\chi^\mu_\rho s^*_\mu(\l), 
\qquad \rho\vdash k, \; \l\vdash l,\; l\ge k.\tag15.21
$$

This implies that the values of an
irreducible character of $S(l)$ on cycle--types of the form $\rho\cup
1^{l-k}$, $k<l$, can be expressed in terms of the character table of the
smaller group $S(k)$ and the $s^*$-functions. Hence, for small $k$ and
arbitrary $l$ and $\l\vdash l$ we can obtain from \tht{15.21} explicit
formulas for the character values. Such formulas were derived 
by several authors: Frobenius, Murnaghan, and Ingram
(see \cite{I} and references therein), Kerov (unpublished), Wassermann
\cite{W}.

\subhead 4. Explicit formula for quantum immanants \endsubhead
The following formula for the quantum immanant $\Sm$
was found by Okounkov \cite{Ok1}.

Denote by  $E=(E_{ij})$ the $n\times n$ matrix formed
by the standard generators $E_{ij}$ of $\Ugln$.
(In section 6 we already used matrices with non-commutative 
entries.)

Let $T$ be a Young tableau of shape $\mu$. Put $k=|\mu|$. Let 
$$
P_T\in\C[S(k)]
$$
be the orthogonal projection onto the corresponding Young
basis vector in the irreducible $S(k)$-module indexed by $\mu$; 
this projection is proportional to the corresponding diagonal matrix
element. Denote by $c_T(i)$ the content of the $i$-th
box in the tableau $T$. We have the following 

\proclaim{\smc Theorem\cite{Ok1}}
$$
\Sm=\tr\, ((E-c_T(1))\otimes\dots\otimes(E-c_T(k))\cdot P_T) \,.
$$
In particular, the right--hand side does not depend on the
particular choice of the tableau $T$.
\endproclaim

Further discussions of this topic 
can be found in \cite{Ok1}, \cite{N2}, and \cite{Ok2}.

\subhead 5. A connection between the binomial formula \tht{5.7} and the
dimension formula \tht{8.3} for skew shapes \endsubhead 
Here we aim to explain why the coefficients in the binomial formula
\tht{5.7} are equal, within simple factors, to the numbers
$\dim\lambda/\mu$. Our argument will follow an idea due to Bingham
\cite{Bi}, Theorem I (see also Lassalle \cite{La}, Macdonald \cite{M4}). 

\proclaim{\smc Lemma} Consider the algebra $\Lambda(n)$ of symmetric
polynomials in $x_1,\ldots,x_n$ and endow it with the inner product
$\la\;,\;\ra$ such that
$$
\la s_{\mu|n},s_{\nu|n}\ra=\delta_{\mu\nu}(n\u\mu).
$$
Further, introduce two operators in $\Lambda(n)$, $D$ and $M$, where
$$
D=\partial/\partial x_1+\ldots+\partial/\partial x_n, \qquad
M=\text{multiplication by $x_1+\ldots+x_n$}.
$$
Then $D$ and $M$ are mutually adjoint with respect to this inner product.
\endproclaim

\demo{\smc Proof} The simplest way to check this is to use the basic
formula \tht{0.1} for the Schur polynomials. \qed
\enddemo

An explanation of this fact is as follows. Consider the reproducing kernel
of the inner product,
$$
{\Cal F}(x,y)=\sum_{\ell(\mu)\le n}\frac{s_{\mu|n}(x)s_{\mu|n}(y)}
{\la s_{\mu|n},s_{\mu|n}\ra},
$$
where $x=(x_1,\ldots,x_n)$, $y=(y_1,\ldots,y_n)$, and let 
$$
X=\diag(x_1,\ldots,x_n), \qquad Y=\diag(y_1,\ldots,y_n)
$$
be the diagonal matrices corresponding to $x$ and $y$. Then it turns out
that 
$$
{\Cal F}(x,y)=\int_{u\in U(n)} e^{\tr XuYu^{-1}}du,
$$
where $du$ is the normalized Haar measure on the unitary group $U(n)$.
{F}rom this formula on easily obtains the relation
$$
(\partial/\partial x_1+\ldots+\partial/\partial x_n){\Cal F}(x,y)=
(y_1+\ldots+y_n){\Cal F}(x,y),
$$
which is equivalent to the claim of the lemma. Finally, note that the
formula 
$$
\int_{u\in U(n)} e^{\tr XuYu^{-1}}du=
\sum_{\ell(\mu)\le n}\frac{s_{\mu|n}(x)s_{\mu|n}(y)}
{\la s_{\mu|n},s_{\mu|n}\ra}
$$
can be obtained from the binomial formula by a limit transition.

Now we apply the lemma to establish a close connection between \tht{5.7}
and \tht{8.3}. Let $\lambda\vdash l$ and $\mu\vdash k$ be two partitions
such that $l>k$, $\ell(\lambda)\le n$, $\ell(\mu)\le n$. We remark that
$$
s_{\lambda|n}(1+x_1,\ldots,1+x_n)=(e^Ds_{\lambda|n})(x_1,\ldots,x_n),
$$
hence the coefficient of $s_{\mu|n}$ in the decomposition of the left--hand
side is equal to
$$
\gather
\frac{\la e^Ds_{\lambda|n},s_{\mu|n}\ra}{\la s_{\mu|n},s_{\mu|n}\ra}=
\frac{\la s_{\lambda|n},e^M s_{\mu|n}\ra}{\la s_{\mu|n},s_{\mu|n}\ra}\\
=\frac{\la s_{\lambda|n},p_1^{l-k} s_{\mu|n}\ra}
{(l-k)!\la s_{\mu|n},s_{\mu|n}\ra}=
\frac{\la s_{\lambda|n},s_{\lambda|n}\ra}
{(l-k)!\la s_{\mu|n},s_{\mu|n}\ra}\dim\lambda/\mu.
\endgather
$$
Thus,
$$
s_{\lambda|n}(1+x_1,\ldots,1+x_n)=
\sum_\mu \frac{(n\u\lambda)}{(|\lambda|-|\mu|)!(n\u\mu)}
\dim\lambda/\mu \cdot s_{\mu|n}(x_1,\ldots,x_n).
$$

\subhead 6. Factorial Schur polynomials from Schubert polynomials 
\endsubhead The following remark is due to Alain Lascoux.

Consider double Schubert polynomials: these are polynomials 
in two sets of variables, $X$ and $Y$, which are indexed by permutations,
see Lascoux \cite{Lasc}, Macdonald 
\cite{M5}. Some special permutations (Grassmannian
permutations) give Schubert polynomials which are symmetric in 
$x_1,\ldots,x_n$. In this case, when one specializes $Y$ to
$\{0,1,2,\ldots\}$, these Schubert polynomials become the factorial Schur
polynomials in $x_1,\ldots,x_n$.

One can also recover factorial Schur polynomials 
(or rather generalized binomial coefficients $\binom{\lambda}{\mu}$, which are
closely connected to them) by specializing $X=\{1,1,\ldots\}$,
$Y=\{0,0,\ldots\}$ and considering the specialized Schur polynomial as a
function of its indexing permutation. The permutations that one has to take are
defined, in a certain way,  by a pair $\lambda,\mu$ of partitions. Then one
obtains the coefficients $\binom{\lambda}{\mu}$.

\subhead 7. Further developments \endsubhead 
Recently Ivanov and Okounkov \cite{Iv}, \cite{IO} 
studied factorial 
analogues of Schur $Q$-functions. Factorial analogues of super Schur
functions are studied in Molev's paper \cite{Mo}. Certain factorial
functions play an important role in Molev--Nazarov work \cite{MN} on
Capelli operators for the orthogonal and symplectic groups. 
The study of general shifted Macdonald polynomials was started in
\cite{KnS,Kn,Sa2} and \cite{Ok3}.

\head References \endhead

\item{[BK]} {\smc E.~A.~Bender and D.~E.~Knuth,} Enumeration of plane
partitions, {\it J. Combin. Theor.\/} A 13 (1972), 40--54.
\item{[BL1]} {\smc L.~C. Biedenharn and J.~D.~Louck}, A new class of 
symmetric 
polynomials defined in terms of tableaux, {\it Advances in Appl. Math.\/} 
10 
(1989), 396--438. 
\item{[BL2]} {\smc L.~C.~Biedenharn and J.~D.~Louck,} Inhomogeneous 
basis set of symmetric polynomials defined by tableaux, {\it Proc. Nat. 
Acad.
Sci.   U.S.A.\/} 87 (1990), 1441--1445.
\item{[Bi]} {\smc C.~Bingham,}  An identity involving partitional 
generalized  binomial coefficients, {\it J. Mutivariate Analysis\/} 4
(1974), 210--223.
\item{[Bou]} {\smc N.~Bourbaki,} Groupes et alg\`ebres de Lie, Chapitres
VII--VIII, Hermann, Paris, 1975.
\item{[CL]} {\smc W.~Y.~C.~Chen and J.~D.~Louck,}  
The factorial Schur function,  
{\it J. Math. Phys.\/} 34 (1993), 4144--4160. 
\item{[D]} {\smc J.~Dixmier,} Alg\`ebres enveloppantes, Gauthier--Villars,
Paris/Bruxelles/Mont\-r\'eal, 1974.
\item{[GG]} {\smc I. Goulden and C. Greene,} A new tableau 
representation for supersymmetric Schur functions, {\it J. Algebra\/} 170
(1994), 687--704. 
\item{[GH]} {\smc I. P. Goulden and A. M. Hamel,} 
Shift operators and factorial symmetric functions, 
{\it J. Comb. Theor. A.\/} 69 (1995), 51--60. 
\item{[H]} {\smc R. Howe,} Remarks on classical invariant theory,  
{\it Trans. Amer. Math. Soc.\/} 313 (1989), 539--570. 
\item{[HU]} {\smc R. Howe and T. Umeda,}  
The Capelli identity, the double commutant theorem, and  
multiplicity--free actions, {\it Math. Ann.\/} 290 (1991), 
569--619.
\item{[I]} {\smc R.~E.~Ingram,} Some characters of the symmetric group,
Proc. Amer. Math. Soc. 1 (1950), 358--369.  
\item{[Iv]} {\smc V.~Ivanov,} Dimension of skew shifted Young diagrams
and projective characters of the infinite symmetric group, to appear.
\item{[IO]} {\smc V.~Ivanov and A.~Okounkov,} in preparation.
\item{[K]} {\smc S.~V.~Kerov,} Combinatorial examples in the theory of 
$AF$-algebras, in: Differential geometry, Lie groups and mechanics X ({\it 
Zapiski Nauchnyh Seminarov LOMI,\/} vol. 172), Nauka, Leningrad, 1989, 
pp. 55--67 (Russian); English translation in {\it J. Soviet Math.} 59, no.
5  (1992). 
\item{[KO]} {\smc S.~Kerov and G.~Olshanski,}  
Polynomial functions on the set of Young diagrams,  
{\it Comptes Rendus Acad. Sci. Paris, S\'er.\/} I, 319 (1994), 
121--126.
\item{[Kn]} {\smc F.~Knop,} Symmetric and non-symmetric quantum
Capelli polynomials, Preprint, q-alg/9603028.
\item{[KnS]} {\smc F.~Knop and S.~Sahi,} Difference equations and
symmetric polynomials defined by their zeros, Preprint, February 1996.
\item{[Lasc]} {\smc A.~Lascoux,} Classes de Chern des vari\'et\'es de
drapeaux, {\it Comptes Rendus Acad. Sci. Paris, S\'er. I, \/} 295 (1982),
393--398.
\item{[La]} {\smc M.~Lassalle,} Une formule de bin\^ome
g\'en\'eralis\'ee pour les polyn\^omes de Jack,  {\it Comptes Rendus
Acad. Sci. Paris, S\'er. I,\/} 310 (1990), 253--256.
\item{[L]} {\smc D. E. Littlewood,} The theory of group characters and  
matrix representations of groups, 2nd edition, Clarendon Press,  Oxford, 
1958. 
\item{[M1]} {\smc I. G. Macdonald,}  
Symmetric functions and Hall polynomials, Oxford University Press, 1979;
2nd edition, 1995.
\item{[M2]} {\smc I. G. Macdonald,}  
Schur functions: theme and variations,  
{\it Publ. I.R.M.A. Strasbourg,\/} 1992, 498/S--27.  
{\it Actes 28-e S\'eminaire Lotharingien,\/} pp. 5--39. 
\item{[M3]} {\smc I.~G.~Macdonald,} Letter of February 22, 1995.
\item{[M4]} {\smc I. G. Macdonald,} Hypergeometric functions, unpublished
manuscript.
\item{[M5]} {\smc I. G. Macdonald,} Notes on Schubert polynomials, Publ.
LACIM, Univ\'ersit\'e du Qu\'ebec, Montr\'eal, 1991.
\item{[Mo]} {\smc A.~Molev,} Factorial supersymmetric functions, 
Research Report. Australian National University, Canberra, November 1995.
\item{[MN]} {\smc A.~Molev and M.~Nazarov,} Capelli identities 
for classical groups, Mathematical Research Report Series 95-21,
University of Wales, Swansea, November 1995.
\item{[MNO]} {\smc A.~Molev, M.~Nazarov, and G.~Olshanski,\/} Yangians and
classical Lie algebras. Research Report CMA-MR53-93. Centre for Mathematics 
and its Applications. Australian National University. November 1993, 
105 pp.; hep-th/9409025. To appear in {\it Russian Math. Surveys\/} 51,
No. 2 (1996).
\item{[N1]} {\smc M. L. Nazarov,}  
Quantum Berezinian and the classical Capelli identity,  
{\it Letters in Math. Phys.\/} 21 (1991), 123--131. 
\item{[N2]} {\smc M.~L.~Nazarov,} Yangians and Capelli identities,
preprint, December 1995; q-alg/9601027. 
\item{[Ok1]} {\smc A.~Yu.~Okounkov,} Quantum immanants and higher Capelli
identities, {\it Transformation Groups\/} 1 (1996), 99--126;
q-alg/9602028.
\item{[Ok2]} {\smc A.~Yu.~Okounkov,} Young basis, Wick formula, and higher
Capelli identities, Preprint, October 1995; q-alg/9602027.
\item{[Ok3]} {\smc A.~Yu.~Okounkov,} (Shifted) Macdonald Polynomials:
$q$-Integral Representation and Combinatorial Formula,
Preprint, April 1996; q-alg/9605013.
\item{[O1]} {\smc G. I. Olshanskii,}  
Unitary representations of infinite--dimensional classical grou\-ps  
$U(p,\infty)$, $SO_0(p,\infty)$, 
$Sp(p,\infty)$ and the corresponding motion groups,  
{\it Funkts. Analiz i 
Pril.\/} 12, no. 3 (1978), 32--44 (Russian); English translation in  
{\it Funct. Anal Appl.\/} (1978), 185--195. 
\item{[O2]} {\smc G. I. Olshanskii,}  
Representations of infinite--dimensional 
classical groups, limits of enveloping algebras, and Yangians, in: Topics 
in representation theory (A.~A.~Kirillov, Editor), Advances in Soviet 
Mathematics, vol. 2, Amer. Math. Soc., Providence, RI, 1991, pp. 1--66. 
\item{[O3]} {\smc G.~Olshanski,} Quasi--symmetric functions and factorial
Schur functions, unpublished paper, January 1995, 21 pp.
\item{[S]} {\smc B.~E.~Sagan,} The symmetric group: representations,  
combinatorial algorithms, and symmetric functions,  Wadsworth \& Brooks  
/ Cole mathematics series, Pacific Grove, CA, 1991.
\item{[Sa1]} {\smc S.~Sahi,} The spectrum of certain invariant differential
operators associated to a Hermitian symmetric space, in: Lie theory and
geometry (J.-L.~Brylinski et al., Editors), Progress in Math. 123 (1994),
569--576. 
\item{[Sa2]} {\smc S.~Sahi,} Interpolation, integrality, and a 
generalization of Macdonald's polynomials, Preprint.
\item{[St]} {\smc R.~P.~Stanley,} Enumerative combinatorics, Vol. 1,
Wadsworth \& Brooks/ Cole  mathematics series, Pacific Grove, CA, 1986.
\item{[VK1]} {\smc A. M. Vershik and S. V. Kerov,}  
Asymptotic theory of 
characters of the infinite symmetric group,  
{\it Funkts. Analiz i Pril.\/} 15, no. 
4 (1981), 17--27 (Russian); English translation in {\it Funct. Anal. 
Appl.\/} 15 (1981), 246--255. 
\item{[VK2]} {\smc A.~M.~Vershik and S.~V.~Kerov,}  
Characters and factor representations of the infinite unitary group,  
{\it Dokl. Akad. Nauk SSSR\/} 267 (1982), 272--276 (Russian);  
English translation in {\it Soviet Math. Dokl.\/} 26 
(1982), 570--574. 
\item{[VK3]} {\smc A.~M.~Vershik and S.~V.~Kerov,} 
Locally semisimple algebras: combinatorial theory and $K_0$-functor, 
{\it J. Soviet Math.\/} 38 (1987), 1700-1750. 
\item{[W]} {\smc A.~J.~Wassermann,} Automorphic actions of compact groups
on operator algebras, Dissertation, University of Pennsylvania, 1981.
\item{[WW]} {\smc E.~T.~Whittaker and G.~N.~Watson,} A course of modern
analysis, Cambridge Univ. Press, 1927.
\item{[Z]} {\smc D.~P.~\v{Z}elobenko,} Compact Lie groups and their
representations, Nauka, Moscow, 1970 (Russian); English translation: 
Transl.\ Math.\ Monogr., Amer.\ Math.\ Soc.\ Providence, RI, 1973.

\bye